# First Families of Regular Polygons

## G.H. Hughes

**Introduction**

Every regular polygon can be regarded as a member of a well-defined 'family' of related regular polygons. These families arise naturally in the study of piecewise rotations such as outer billiards. In some cases they exist at all scales and can be used to define the fractal dimension of the 'singularity set'. This is well-documented for regular N-gons such as the pentagon, octagon and dodecagon, whose algebraic complexity is 'quadratic' ($\varphi(N)/2 = 2$ where $\varphi$ is the Euler totient function).

Recent evidence suggests that the geometry of these families is intrinsic to the 'parent' polygon and can be derived independently of any mapping. It is the purpose of this paper to show how the First Family for any regular polygon arises naturally from the geometry of the 'star polygons' first studied by Thomas Bradwardine (1290-1349).

Any star polygon based on a regular N-gon should share the same algebraic complexity as the embedded N-gon. This complexity is traditionally defined to be the degree of the minimal polynomial of $2\cos(2\pi/N)$ which is always $\varphi(N)/2$. We will show that there are advantages to basing this complexity on the tangent function, because when N is at the origin with apothem h, the vertices of the nested star polygons are of the form $\pm h\{s_k, 1\}$ where $s_k = \tan(k\pi/N)$. Each of these 'star points' defines a scale and the star points and the scales together define the First Family.

The $s_k$ with $\gcd(k,N) = 1$ will be called 'primitive' because they play a role similar to the primitive roots of unity. Based on a 1949 result by Carl Siegel, communicated to Sarvadaman Chowla [Ch], these primitive star points and matching primitive scales are independent and form a basis for all the star points and scales. This implies that the algebraic 'complexity' of the star points and scales matches the complexity of N - and the First Family is a window on this complexity.

In studies of piecewise rotations based on regular N-gons, the rotational parameters and scaling parameters typically lie in $\mathbb{Q}_N^+$ - the maximal real subfield of the cyclotomic field $\mathbb{Q}_N$. This is the subfield of order $\varphi(N)/2$ which can be generated by $\cos(2\pi/N)$. We will show that the primitive scales are a unit basis for this 'scaling' subfield, so they are compatible with a wide range of studies such as Ashwin [A], Goetz [Go1], Lowenstein[L], Kahng [K], Vivaldi [LV] and Poggiaspalla [GP]. The investigations in Appendix C, support our conjecture that all 'tiles' that arise in the outer billiards map can be defined by scaling in $\mathbb{Q}_N^+$. This means that these tiles are known to arbitrary accuracy, and working within $\mathbb{Q}_N^+$ is very computationally efficient.

The First Families to be defined here are a natural extension of the 'parent' polygon, so it is hard to explain why they have not been studied earlier. One explanation is that they are really a 'dynamic' extension of the regular polygons, and this is an area that is still under development. Since these are dynamic issues, there is natural path toward recursion, and in some cases there are families at all scales. This is well-documented for the pentagon [T], heptagon [H2] and octagon [S2]. These extended 'families' and their symmetric counterparts can be seen below.

| First Families and extended families for the regular pentagon, heptagon and octagon |
|---|
| 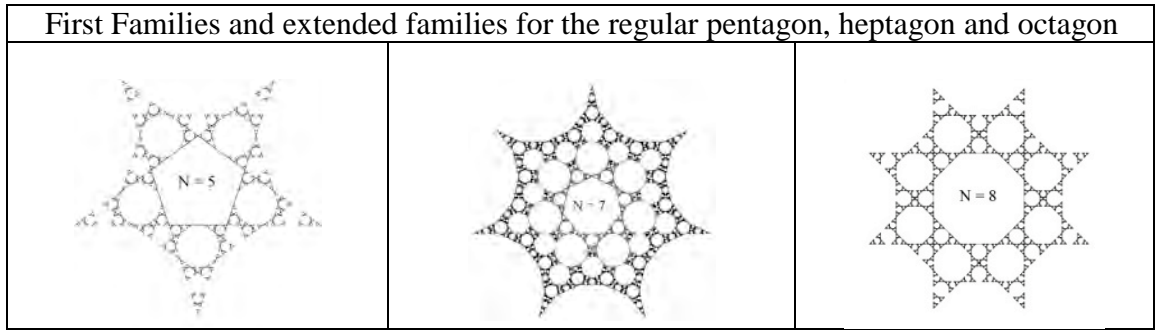 |

In 1978 Jurgen Moser [M2] suggested that the 'outer billiards' map for a convex polygon, could be used as a 'toy model' for the study of stability for 'conservative' systems such as the Solar System. This follows because the outer billiards map has properties similar to a discrete Hamiltonian system. One of the issues raised by Moser has already been settled – namely the existence of unbounded orbits. In the regular case, Vivaldi and Shaidenko [VS] showed in 1989 that all orbits were bounded, but in 2006, Richard Schwartz [S1] showed that non-regular polygons can support unbounded orbits. See chronology for a history of the outer billiards map.

Below are examples of the outer billiards map $\tau$ for a non-regular polygon and for a regular pentagon. For an initial point p each iteration is a reflection (central symmetry) about one of two possible 'support' vertices (clockwise here), so the formula is $\tau(p) = 2c - p$ where c is the support vertex.

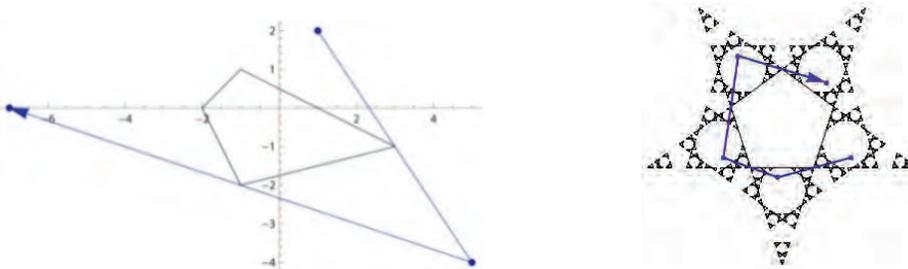

There is a growing consensus that the outer billiards map and related piecewise affine maps on a polygon, can uncover geometric and algebraic relationships that are inherent in the polygon itself. In [H2] the author defines the First Family of a regular polygon as the major 'resonances' of the outer billiards map – and it is the purpose of this paper to show that these families can be derived independently of the outer billiards map. This would give evidence that these structures are an inherent part of the geometry of any regular polygon. Algebraic evidence is presented in Section 4.

## Section 1: Star points and scales of regular polygons

As a first step, we will look at the 'star polygons' first studied by Thomas Bradwardine (1290-1349), and later by Johannes Kepler (1571-1630).

The vertices of a regular n-gon with radius r are **Table[{rCos[2πk/n],rSin[2πk/n]},{k,1,n}].** A 'star polygon' {p,q} generalizes this by allowing n to be rational of the form p/q so the vertices are given by:
$$\{p,q\} = \text{Table[\{rCos[2πkq/p],rSin[2πkq/p]\},\{k,1,p\}]};$$

Using the notation of H.S.Coxeter [Co] a regular heptagon can be written as {7,1} (or just {7}) and {7,3} is a 'step-3' heptagon formed by joining every third vertex of {7} so the exterior angles are 2π/(7/3) instead of 2π/7.

By the definition above, {14,6} would be the same as {7,3}, but there are two heptagons embedded in N = 14 and a different starting vertex would yield another copy of {7,3} - so a common convention is to define {14,6} using both copies of {7,3} as shown below. This convention guarantees that all the star polygons for {N} will have N vertices.

| {7,1} (a.k.a. N = 7) | {7,3} | {14,6} | {10,2} |
|---|---|---|---|
| 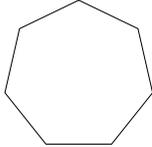 | 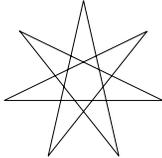 | 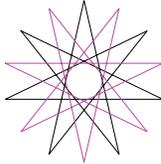 | 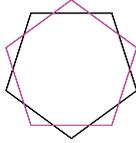 |

The number of 'distinct' star polygons for {N} is the number of integers less than N/2 – which we call 'HalfN' and write as ⟨N/2⟩. So for a regular N-gon, the 'maximal' star polygon is {N, ⟨N/2⟩}. (Some authors would also allow the 'boundary' cases such as {14,0} and {14,7} – which are isolated points or 'asterisks'.)

Our default convention for the 'parent' N-gon will be centered at the origin with 'base' edge horizontal, and the matching {N,1} will be assumed equal to N. In general sN, rN and hN will denote the side, radius and height (apothem). Typically we will use hN as the lone parameter so the 'cyclotomic' case of r = 1 will have h = Cos[π/N] – which is always in $\mathbb{Q}_N^+$.

**Definition**: The *star points* of a regular N-gon are the intersections of the edges of {N, ⟨N/2⟩} with a single extended edge of the N-gon (which will be assumed horizontal).

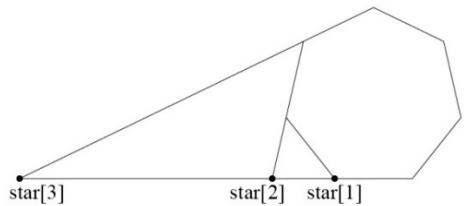

By convention the star points are numbered from star[1] (a vertex of N) outwards to star[⟨N/2⟩] – which is a vertex of {N, ⟨N/2⟩}. So every star[k] is a vertex of {N,k} embedded in {N,⟨N/2⟩}.

(It would seem natural to define the star points to be on the 'positive' side of N, but over the years we have adopted a clockwise rotation for outer-billiards, which makes it convenient to use negative star points. The symmetry between these choices makes it irrelevant which one is used.)

**Lemma 1**: The star points of a regular N-gon with apothem h are

$$\text{star}[k] = -h\{s_k, 1\} \text{ where } s_k = \tan(k\pi/N) \text{ for } 1 \leq k < N/2$$

**Proof**: Since $\tan(k\pi/N) = \tan(2k\pi/2N)$, the indices divide N into 2N segments centered at the origin with slopes given by $s_k$. □ (In particular $2s_1 = 2\tan[\pi/N]$ is always $sN/hN$).

**Example**: The two star points of N = 6 are defined using rotation by $\pi/6$. N = 6 and N = 12 are shown here in 'standard position' so that the star polygons can be compared.

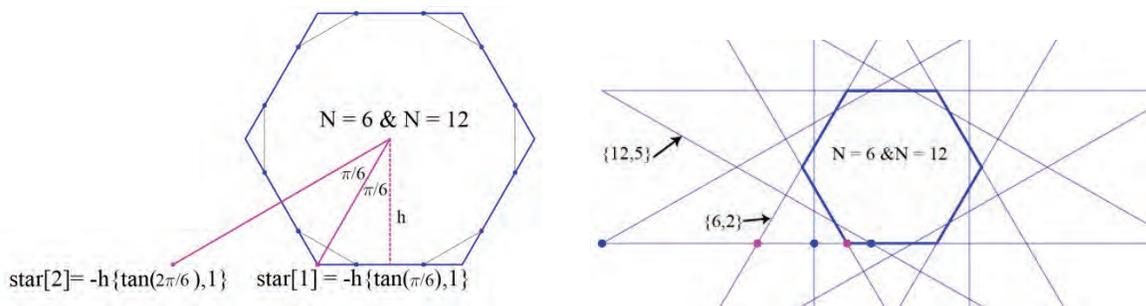

**Example**: Below is the same comparison with N = 14 and N = 7 and now there is a perfect match between the 'even' star points of N = 14 and the three star points of N = 7. This will imply that N = 14 and N = 7 have equivalent cyclotomic fields. From the standpoint of the outer billiards map, the 'rotationally equivalent' star polygons will imply identical singularity sets, and 'conjugate' dynamics. This cannot occur when N is twice-even, because <N/2> will be odd. Therefore from the standpoint of outer billiards, N = 6 and N = 12 will have very different dynamics – and in addition, the members of the $2^k$ family appear to be virtually unrelated.

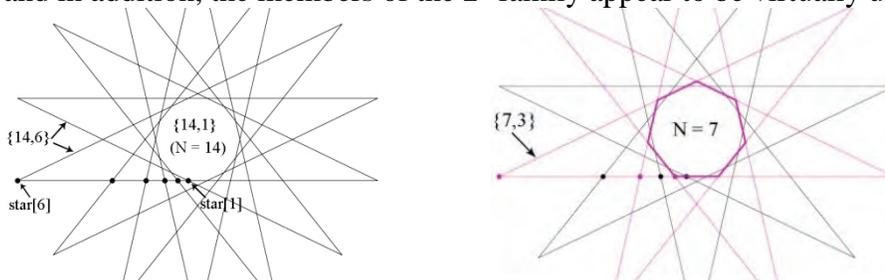

Since $-s_k = -\tan(k\pi/N) = \tan(-k\pi/N)$ the 'positive' star points can be treated as in the same fashion as the 'negative' star points by swapping k and -k. For the most part we will only assume that k is an integer, so the results will apply equally to all star points. We will sometimes refer to the $s_k$ as the 'star' points.

The set $S = \{s_k: (k,N) = 1\}$ will play an important part in what follows. These are called the 'primitive' star points because they play a role similar to the primitive roots of unity for N. The scales formed from these primitive star points will form a unit basis for the maximal real subfield of $\mathbb{Q}_N$ so they will have the same rank as the units in $\mathbb{Z}[\zeta_N]$.

There is a long history of interest in trigonometric functions of rational multiples of $\pi$. In 1949 C.L.Siegel (see chronology) communicated to S.Chowla [Ch] a proof that the primitive $s_k$ for the cotangent function are linearly independent. This was a non-trivial result in algebraic number theory and Siegel only proved it for prime N. In 1970 Chowla generalized this result using character theory and Dirichlet's $\mathcal{L}$-series, but the general tangent case was only settled recently in [G2] (1997). Therefore the set $S = \{s_k = \tan(k\pi/N): 1 \le k < N/2, (k,N) = 1\}$ is linearly independent and it follows (Lemma 9) that any non-primitive $s_k$ must be a $\mathbb{Q}$-linear combination of the primitive $s_k$, although there may be no elementary formulas which yield the coefficients of these 'degenerate' $s_k$. This means that there may be non-trivial relationships between the star points of many star polygons. All of this parallels traditional Galois theory - which is outlined in Section 4.

We will define $T_N[x]$ = **MinimalPolynomial[Tan[Pi/N]][x]**. A portion of $T_{14}[x]$ and $T_7[x]$ are shown below in blue and red. Both polynomials are symmetric and order 6 but only $T_7[x]$ is monic because, $\tan(\pi/N)$ is monic iff $N \ne 2p^k$ for an odd prime p. The six roots of $N = 14$ are shown in black and four of the six roots of $N = 7$ are shown in magenta.

Both polynomials are needed to generate the full set of star points of $N = 14$, because each polynomial only generates the (positive and negative) **primitive** roots - so $s_1$ is primitive with respect to $N = 14$ and $s_2$ is degenerate, but $s_2$ is primitive with respect to $N = 7$.

The indices of the roots of $N = 14$ are $k = \pm 1, \pm 3$ and $\pm 5$ with $(k,14) = 1$. By definition, these are the prime residue classes mod 14, so they can be used to define the Galois group $G_{14}$. Each element k defines an automorphism on the cyclotomic field $\mathbb{Q}(\zeta_{14})$ given by $\sigma_k : \zeta_{14} \to \zeta_{14}^k$ with $\zeta_{14}$ a primitive root of unity of $N = 14$. Therefore each $\sigma_k$ is in $G_{14}$ and $|G_{14}| = [\mathbb{Q}(\zeta_{14}):\mathbb{Q}] = 6$. The indices $\{\pm 1, \pm 2$ and $\pm 3\}$ of $N = 7$ define the same Galois group and in fact $\mathbb{Q}(\zeta_7) \approx \mathbb{Q}(\zeta_{14})$ because $(-\zeta_7)$ will generate all of $\mathbb{Q}(\zeta_{14})$. This is only true for 2N,N pairs with N odd.

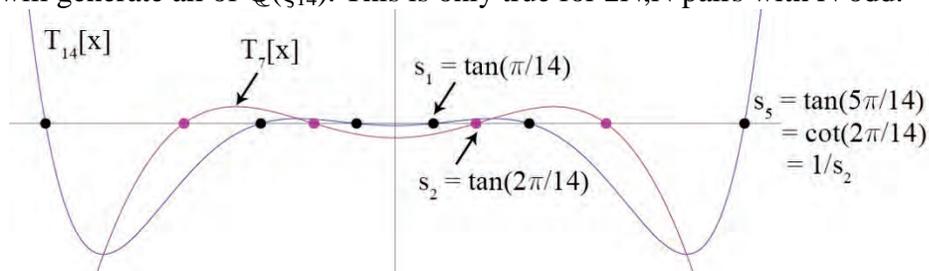

It is not obvious that the three primitive $s_k$ of $N = 7$ or $N = 14$ are linearly independent over $\mathbb{Q}$. They are Galois conjugates but that clearly does not imply independence – however it does show that they share the same number field- $\mathbb{Q}(s_1)$ (although $s_1$ may not be in $\mathbb{Q}(\zeta)$ unless $4|N$). Since $T_{14}[x]$ and $T_7[x]$ are minimal and degree 6, the degree of $\mathbb{Q}(s_1)$ is 6 for either $s_1$, and one basis is $\{1,\alpha,\alpha^2,..,\alpha^5\}$ where $\alpha = s_1$. (This basis is only integral for $N = 7$). These basis elements must be linearly independent but this does not imply that the primitive $s_k$ are linearly independent. In fact it is often easier to prove the 'dual' result about the cotangent function – although this does not directly imply the tangent result.

**Definition:** The *dual* star points are $c_k = \cot(k\pi/N)$ for $1 \le k < N/2$ and the primitive dual star points satisfy $(k,N) = 1$.

Therefore if $s_k$ is a primitive star point (w/r tan), $c_k = 1/s_k$ is always primitive (w/r cot) because the primitive k indices do not change. This does not imply that $1/s_k$ is always primitive (w/r to tan). For example, $s_5$ for N = 14 has inverse $s_2$ as shown in the illustration above.

Since $\tan[\pi/2 - \theta] = \cot[\theta]$, if Mathematica is asked to list the star points for a regular N-gon it will list the first 'half' using tangents and the second 'half' using cotangents. When N is even these two half-lists will match up perfectly (in reverse order) because symmetry demands that
**star[N] = Reverse[1/star[N]].**

**Definition**: For a regular N-gon, *scale[k]* = $s_1/s_k$ with the primitive scales satisfying (k,N) = 1

So scale[1] is always 1 and the scales of any N-gon are strictly decreasing. Of course these scales are independent of size or orientation so they are fundamental parameters of any regular N-gon. Since the dual scales are just the inverse of the scales, there is a useful tan-cot duality. The smallest scale, scale[<N/2>], is the matching scale for GenStar[N] so it is called GenScale[N].

For N-gon scaling, any of the basic trig functions would work, but clearly $\tan(\pi/N)$ mediates height/side, while $\cos(\pi/N)$ gives height/radius, and $\sin(\pi/N)$ yields side/radius. We will show that there are many advantages to using tan as the primary function. One reason why J. Stillwell [St2] felt that the " *tan function can be considered more fundamental than either cos or sin.*" is because of the 'half-angle' formulas below that express cos and sin as rational functions of tan:

$$\cos(2\theta) = \frac{1-\tan^2(\theta)}{1+\tan^2(\theta)} \quad \text{and} \quad \sin(2\theta) = \frac{2\tan(\theta)}{1+\tan^2(\theta)}$$

These relationships are based on the work of Diophantus in about 250 AD as part of his quest to find rational points on curves - in this case the unit circle. Geometrically they are based on the same 'half-angle' diagram used above to illustrate the star points of N = 6. This diagram is repeated here with N = 6 scaled to have radius 1.

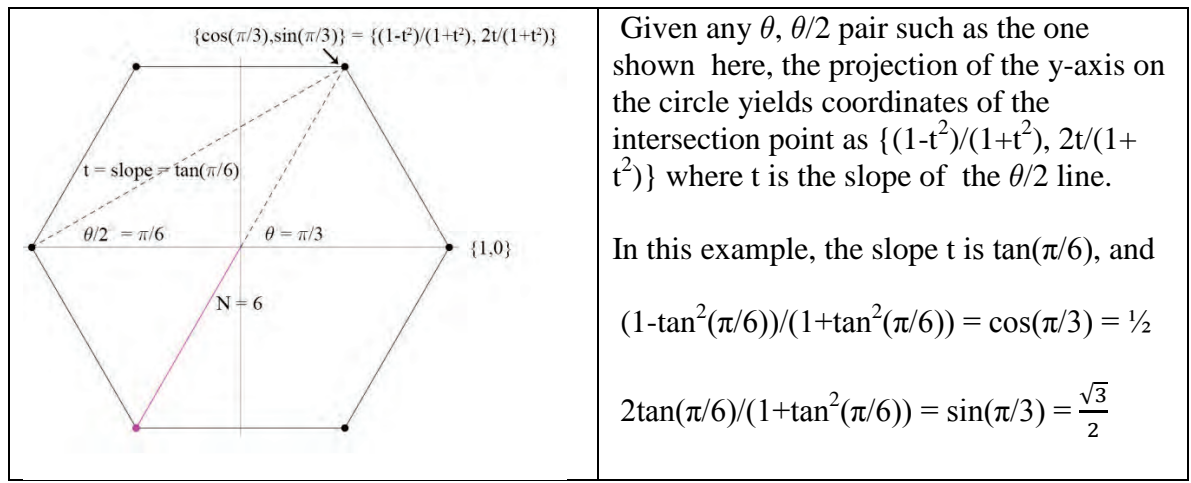

Given any $\theta$, $\theta/2$ pair such as the one shown here, the projection of the y-axis on the circle yields coordinates of the intersection point as $\{(1-t^2)/(1+t^2), 2t/(1+t^2)\}$ where t is the slope of the $\theta/2$ line.

In this example, the slope t is $\tan(\pi/6)$, and

$(1-\tan^2(\pi/6))/(1+\tan^2(\pi/6)) = \cos(\pi/3) = \frac{1}{2}$

$2\tan(\pi/6)/(1+\tan^2(\pi/6)) = \sin(\pi/3) = \frac{\sqrt{3}}{2}$

Parameters such as $\lambda_N = 2\cos(2\pi/N)$ are traditionally used for algebraic investigations in the context of the cyclotomic field $\mathbb{Q}_N$, because $\lambda_N$ is an algebraic integer and generates the maximal real subfield, $\mathbb{Q}_N^+$. Based on the half angle formula above, $s_1^2 = \tan^2(\pi/N)$ can also serve as a generator of $\mathbb{Q}_N^+$. By symmetry, $\tan^2(\pi/N)$ is GenScale[N] when N is even, but as noted above, when N is twice-odd, $\tan^2(\pi/N)$ is no longer an algebraic integer – but GenScale[N/2] can be used instead. We will see that in all cases either GenScale[N] or GenScale[N/2] can serve as integral generators for $\mathbb{Q}_N^+$ and in fact they are units – while $\lambda_N$ fails to be a unit when 4|N. Also the independence results of Siegel and Chowla fail for sin and cos.

Below is the First Family for N = 14 showing a magenta copy scaled by GenScale[7] = $\tan(\pi/7)/\tan(3\pi/7) \approx .109916$. Note that the new family fits perfectly on a side of N = 14. This is called 'generation' scaling and it is an invitation for recursion. Such families are not unusual for the outer billiards map – but they are often 'broken' families which are just a skeleton of the First Family. However hidden in these second generation 'tiles' for N = 14 are the seeds for a perfect third generation - scaled by GenScale[7]$^2$.

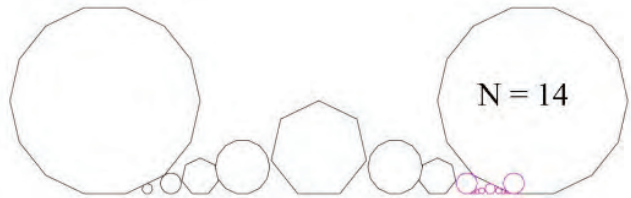

When N is even, GenStar[N] = 1/star[1] so the scales inherit the symmetry of the star points:

**scales[N] = Reverse[GenScale[N]/scales[N]]**

Usually to compare the scales of an N-gon and an M-gon it makes sense to assume that they are both in 'standard position' with the same height as show below for N = 7 and N = 14. In this case the ratio of the sides is sN/sM = $\tan(\pi/N)/\tan(\pi/M)$ - which we call **ScaleSwap[N,M].** (If instead the sides are equal then ScaleSwap[N,M] will be the ratio of the heights.)

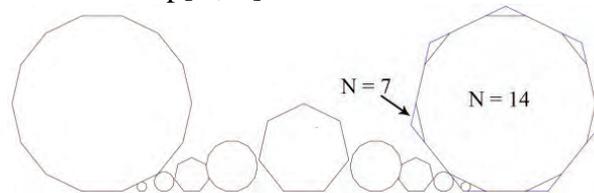

For matching heights, the ratio of sides is the same as the ratio of any matching scales, so ScaleSwap[14,7] = GenScale[14]/GenScale[7] ≈ 0.4739524. These definitions and observations imply that the scales of any twice-odd N-gon are essentially the same as those of the matching N/2-gon. This will be made precise in the scaling lemmas of Section 3. Here is a summary of these results as they apply to N = 14 & N = 7:

| N = 14 (scale[6] = GenScale[14]) | N = 7 (scale[3] = GenScale[7]) |
|---|---|
| scale[6] = scale[6]/scale[1] | scale[3] = scale[6] (of N = 14)/ ScaleSwap[14,7] |
| scale[5] = scale[6]/scale[2] | scale[2] = scale[4] (of N = 14)/ ScaleSwap[14,7] |
| scale[4] = scale[6]/scale[3] | scale[1] = scale[2](of N = 14)/ ScaleSwap[14,7] |

**Lemma 2**: GenScale[N] = $\tan(\pi/N)^2$ if N is even and $\tan(\pi/N)\cdot\tan(\pi/2N)$ if N is odd

**Proof**: GenScale[N] is defined to be $s_1/s_{<N/2>}$. When N is even $s_{<N/2>}$ is $1/s_1$ by symmetry, so GenScale[N] = $s_1^2$ = $\tan(\pi/N)^2$. When N is odd GenScale[N] = GenScale[2N]/ScaleSwap[N,2N] = $\tan(\pi/2N)^2 \cdot \tan(\pi/N)/\tan(\pi/2N) = \tan(\pi/N)\cdot\tan(\pi/2N)$ □

In the next section we will use these scales to define the 'First Family' for every regular N-gon. The 11 members of the First Family for N = 14 are shown in cyan below.

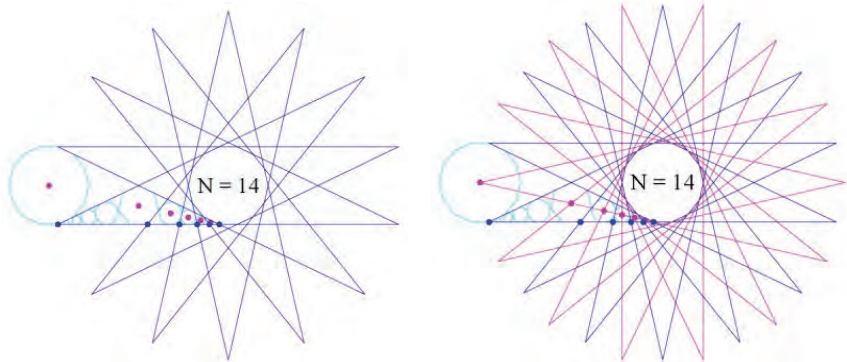

These are sometimes called 'tiles'. The largest tile in the First Family will be called $D_N$ or simply D – when N is understood. In the outer-billiards world, D is globally maximal among all tiles and there are unbroken rings of D tiles at all distances from the origin - which guarantees bounded dynamics. When N is even, we will see that $D_N$ is an identical copy of N and when N is odd, $D_N$ is a regular 2N-gon with the same side length as N. (It is easy to show that no tile can have more than 2N sides.) $D_{14}$ shown here in cyan is clearly identical to N= 14. We will define the magenta {14,6} on the right to be the 'outer-dual' of the blue {14,6}. (This looks like a simple rotation of the blue {14,6}, but there is an important scale change as well.)
.
**Finding the Centers**

The classic 'dual' of a polygon involves swapping edges and vertices. Since a regular N-gon is self-dual, we will define the 'outer dual' to be a circumscribed copy of N which is rotated by a half-turn. This will mimic the natural duality of the outer billiards map.

**OuterDual[x_, N_] := RotationTransform[-Pi/N][x*rN/hN]**

OuterDual[{14/6},14] is the magenta {14,6} on the right above. As such it is scaled by radius/height of N = 14 which is ≈1.025716863. This magenta {14,6} is the star polygon for a regular 14-gon which is circumscribed about N = 14 and rotated by a 'half-turn' (the 'outer dual' of N = 14). Therefore it is the outer-billiards orbit of the center of D, for N = 14. This is the single most important orbit for a regular polygon – because it can be used to define the parameters of all the remaining canonical resonances – which are the First Family tiles.

**Definition**: For a regular N-gon, each star[k] and matching scale[k] defines a corresponding First Family member which we call S[k]. These S[k] tiles are regular 2N-gons or N-gons with centers given by **cS[k] = OuterDual[star[k], n].** (So centers are dual to star points.)

The constraints of the star polygons provide a simple geometric relationship between the star point, scales and the S[k] tiles. This is given in the Lemma below. We will use the maximal 'D' tile as a reference for scale because it is uniform across all regular 'parent' N-gons.

**Lemma 3**: For every regular N-gon, the scales of the S[k] tiles are in direct proportion to Distance[cS[k],cS[0]] – where cS[0] is the star[1] vertex of N. Therefore hS[k] is always hD·GenScale[N]/scale[k] and all First Family tiles can be scaled relative to the maximal tile D.

**Example:** N-odd ( N = 7 - with height h ). D is a regular 14-gon with side the same as N = 7, so hD = h·ScaleSwap[7,14]. The ratio of center distances gives the three scales. Since star points map to centers, the scaled ratio of star distances would suffice also.

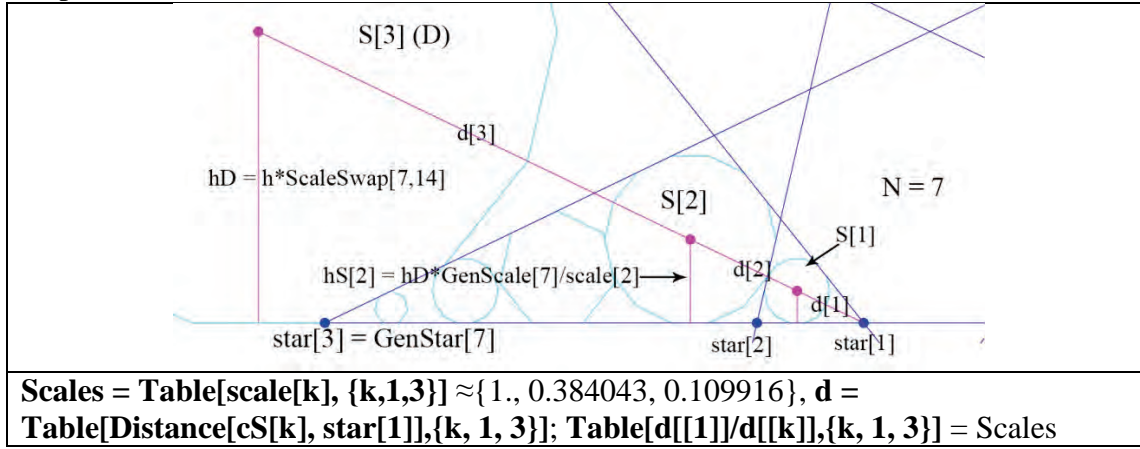

**Scales = Table[scale[k], {k,1,3}]** ≈ {1., 0.384043, 0.109916}, **d = Table[Distance[cS[k], star[1]],{k, 1, 3}]; Table[d[[1]]/d[[k]],{k, 1, 3}]** = Scales

**Example:** N twice-odd (N = 14). D is a clone of N so hD = h. This is the only case where the S[k]'s alternate regular N-gons and regular N/2 – gons, but the heights are uniform.

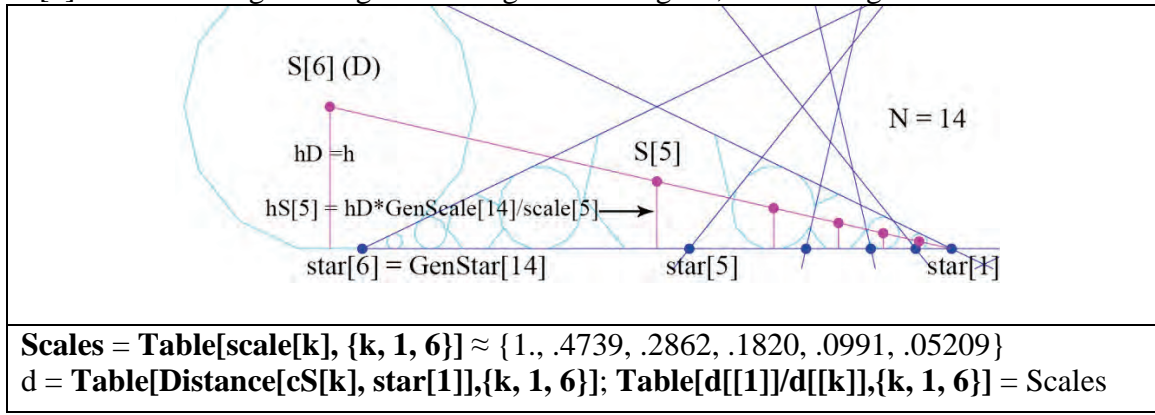

**Scales = Table[scale[k], {k, 1, 6}]** ≈ {1., .4739, .2862, .1820, .0991, .05209}
d = **Table[Distance[cS[k], star[1]],{k, 1, 6}]; Table[d[[1]]/d[[k]],{k, 1, 6}]** = Scales

**Example:** N twice-even (N = 16). D is a clone of N but now all the S[k] are regular N-gons)

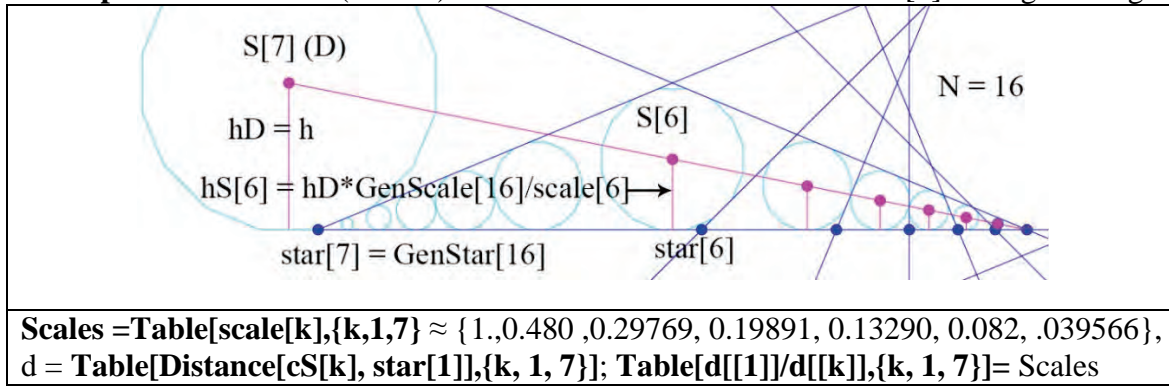

**Scales =Table[scale[k],{k,1,7}** ≈ {1.,0.480 ,0.29769, 0.19891, 0.13290, 0.082, .039566},
d = **Table[Distance[cS[k], star[1]],{k, 1, 7}]; Table[d[[1]]/d[[k]],{k, 1, 7}]**= Scales

These definitions of star points, scales, centers and heights apply to any regular N-gon, but these examples make it clear that the First Families will be slightly different depending on whether N is 1,2,3 or 4 mod 4. The odd cases of 1 and 3 mod 4 can be treated together here, but we believe that there are important dynamical distinctions between these cases Thus we will cover the three cases of odd, twice-odd and twice-even below – beginning with N twice-odd.

**Section 2: Defining the First Family from centers and scales**

**Case 1: N Twice-Odd**.

In all cases the nucleus of the First Family consists of the S[k] tiles. Based on Lemma 2 we know that hS[k] = hD*GenScale[N]/scale[k] and since the centers are known, the only issue is whether the S[k] are N-gons or 2N-gons. The answer is simple: when N is odd, the S[k] tiles will all be 2N-gons and when N is twice-even the S[k] will all be N-gons, but when N is twice-odd, the S[k] will alternate N/2- gons and N-gons - and we will show why this is true.

The default scaling is simply a scaled version of D, so the 'even' S[k] are the expected result. To explain the odd tiles, we will look at the familiar case of N = 14 – where hS[k] = hGenScale[14]/scale[k]. Starting with S[5], GenScale[14]/scale[5] = scale[2] because N is even. But only when N is twice-odd, scale[2] will be GenScale[14]/GenScale[7], so s[5] is scaled (relative to D or N = 14) by ScaleSwap[14,7] – which is by definition the scale that would convert a regular 14-gon to a regular 7-gon with the same side. Therefore hS[5] = h·tan($\pi$/14)/tan($\pi$//7) ≈ h·(.4739) . This tile will always have an 'M-D' relation with D (and N)

The scaling of the other 'odd' tiles will be similar because GenScale[N]/scale[k] will always be an 'even' scale of N and for N twice-odd these are exactly the scales shared by N and N/2. So for S[3] of N = 14, hS[3] = GenScale[14]/scale[3] = scale[4]. The table above (and the Scaling Lemma of Section 3) says that scale[4] of N = 14 is scale[2] (of N = 7)·ScaleSwap[14,7], so S[3] is first scaled to match M and then scaled by scale[2] of N = 7. This makes S[2] part of a scaled version of the N = 7 First Family, which is embedded in the N = 14 First Family.

The k = 1 case for N = 14 is very important, because hS[1] = GenScale[14]/scale[1] = scale[6] = scale[3] (of N = 7)*ScaleSwap[14,7]. Therefore S[1] is scaled by scale[3] relative to M, but scale[3] of N = 7 is GenScale[7], so S[1] is a 'second generation' M tile as shown in the examples above.

Therefore the odd tiles of N = 14 are actually scaled by N = 7 – and it is fair to say that dynamically N= 14 is just a 'different view' of N = 7. (The Twice-Odd conjecture of Section 3.)

Any N-gon and N/2-gon pair such as M and N = 14 are called an *M-D pair*. We pointed out earlier the important role that D tiles play in outer billiards. Except when N is twice-even, the M tiles evolve from matching D tiles (as shown here) and conversely D tiles evolve from M tiles, so these M-D pairs are fundamental building blocks for extended families for N odd or twice-odd.

(In 1989 when the author met with Jurgen Moser at Stanford University to discuss these canonical tiles, the issue of boundedness for regular polygons had just been settled, but he agreed that a study of these 'families' would be an interesting exercise in 'recreational mathematics'. That is exactly what it became – part of a whimsical 'fairy-tale' with M and D as matriarch and patriarch. But after all these years, there is still no theory for the extended families – except for the simple cases of N = 3, 4, 5, 6, 8,10 and 12. However the 'gender issue' between tiles with even or odd number of sides is still fundamental – because they have different dynamics relative to the outer billiards map. Maybe the astute reader can see why.)

**Definition**: The **First Family** of a twice-odd N-gon is defined as follows:

Assume N is centered at the origin with apothem h and lower edge horizontal.

(i) The ⟨N/2⟩ star points of N are **star[n_]:= Table[-h*{Tan[k*Pi/n],1},{k, 1, ⟨N/2⟩}]**

(ii) Each star point defines a scale using **scale[k] = $s_1/s_k$** (where $s_k = \tan(k\pi/N)$)

(iii) The 'centers' are defined to be
    cS[k] = **OuterDual[star[k]]= RotationTransform[-Pi/N][star[k]*rN/hN]**

(iv) For k even, S[k] will be a regular N-gon with center cS[k] and hS[5] = GenScale[N]/scale[k]

(v) For k odd, [k] will be a regular N/2-gon with center cS[k] and hS[k] = GenScale[N]/scale[k]

(S[⟨N/2⟩] will be called D and S[⟨N/2⟩-1] will be called M. N is also be known as S[0])

(vi) By symmetry, each S[k] defines a matching LeftS[k] (LS[k]) as follows:

    **LS[k_]:=ReflectionTransform[{1,0},cM][S[k]];**

(where the vector {1,0} specifies a 'horizontal' reflection about cM)

    **FFRight[N_] :=Table[S[k],{k,0,⟨N/2⟩}]; FFLeft[N_] := Table[LS[k],{k,0, ⟨N/2⟩ }];**

    **FirstFamily[N_] := Join[FFRight[N], FFLeft[N]];**

**Example: Graphics[poly/@FirstFamily[14]]**

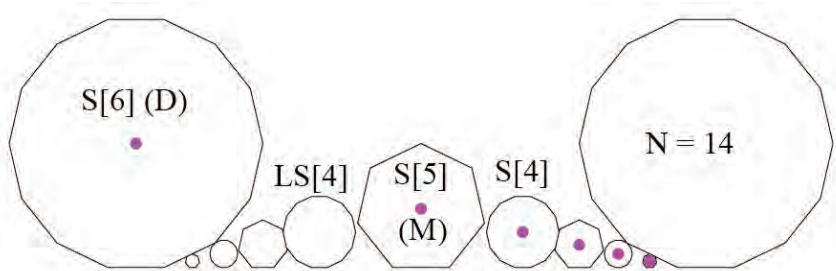

Note: To display these closed polygons in Mathematica using the Line command it is necessary to append a vertex to explicitly close the circuit. We use **poly[M_]:=Line[Append[M, M[[1]]]];**

**Example**: **Graphics[poly/@FirstFamily[26]]** (in black)

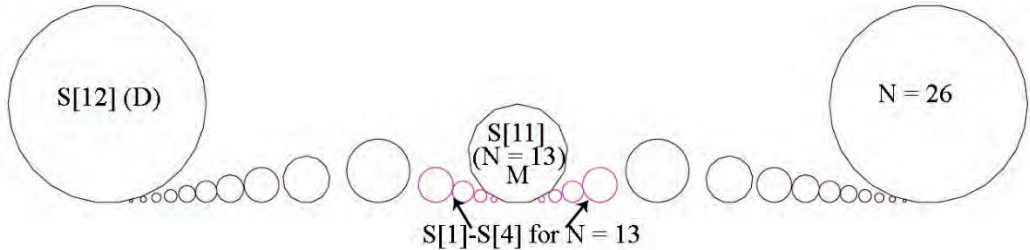

As N increases, it appears that the central 'M' tile will be further isolated from the rest of the First Family, but the Twice-Odd Lemma says that that the central S[11] tile is the (scaled) matriarch of the First Family for N = 13 – which will be defined next. The transformation that imports this First Family is **kTo2k[x_]:=TranslationTransform[cM][x*rM].** The magenta tiles above are the result of applying this transformation to just the 'missing' S[k] tiles of N = 13. If we had imported the full First Family for N = 13 it would include M and all the tiles left of M – so M and D shown here, are a canonical M-D pair. By convention we refer to a combined plot like this as the 'full' First Family, so it is FFF[26] - and as such it contains a scaled copy of the First Family for N = 13 – by simply cutting it in half – as in Parent Trap.

N = 13 (and N = 26) have 'generations' of families in a manner similar to N = 7 (and N = 14). In both cases, generation[k] is presided over by an M-D pair scaled (relative to the First Family) by GenScale[7]$^k$ and GenScale[13]$^k$ respectively – but the other 'tiles' in these families may be quite different from the S[k] tiles of the First Families. For N = 7, the difference is fairly minimal and all the even and odd generations appear to be self-similar, but N = 13 may have generations which never settle down to a recognizable pattern. The 4k+1 conjecture in [H2] predicts that all N-gons where N is prime of the form 4k+1 will have endless chains of generations. The twice-even cases have no M-D pairs, but the canonical S[1] and S[2] tiles can serve as surrogates. For N = 12 and N = 16 there appear to be chains of such generations, but the rest of the twice-even family shows little sign of generation scaling. The issue is complicated by the mutations.

**Case 2: First Family for N-odd**

The star points and centers of N= 7 are shown below. They are defined the same for all N-gons so the magenta {7,3} is the OuterDual of the blue {7,3}. Since the blue {7,3} is one component of {14,6}, these star polygons have very similar First Families and each can generate the other.

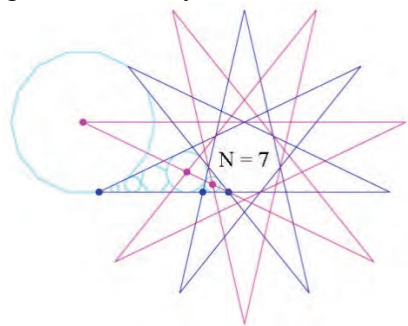

The S[k] are always scaled relative to the maximal D tile. When N is odd, D and N will share a side and form an M-D pair, so hD = ScaleSwap[N,2N] . Therefore hS[k] = ScaleSwap[N,2N]· GenScale[N]/scale[k]. For N = 7 above the list of heights for S[1], S[2] and S[3] are {0.231914,0.603875,2.10992}  with the last being hD. These are all 2N-gons by default.
The outer billiards web and related star polygons have a natural duality which allows us to view the D tiles as if they were 'generating' N-gons at the origin, so the 'First Family' for D within N = 7 is a scaled (and reflected) version of the First Family for N = 14 described in Case 1. To see this connection, we will import the LS[k] for N = 14. These tiles are shown below in magenta. Within N = 7, the imported tiles will be called DS[k] because reflection in D is their true origin.

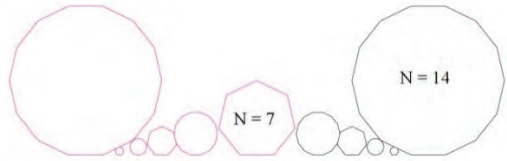

All of these will be part of the First Family for N = 7. Of course these magenta tiles can also be generated within the context of N = 7, using D as a clone of N = 14.

**Definition**: The **First Family** of an N-gon for N odd is defined as follows:
Assume that N is centered at the origin with hN = 1 and lower edge horizontal. (Items (i),(ii),(ii) below are identical for all regular N-gons.)

(i) The star points are: **star[N_]:= Table[{-Tan[k*Pi/n],-1},{k, 1, ⟨N/2⟩}]**

(ii) Each star point defines a scale using **scale[k] = $s_1/s_k$**

(iii) The 'centers' are defined to be
$$cS[k] = \textbf{OuterDual[star[k]]= RotationTransform[-Pi/N][star[k]*rN]}$$
(iv) The S[k] tiles will by regular 2N-gons with center cS[k] and height hD*GenScale[7]/scale[k] = ScaleSwap[14]/scale[k].

(v) Import FFLeft[2N] - which consists of the LS[k] for {k,0, N-2}. Recall that LS[N-2] is called M. Scale each LS[k] using the scale that would 'promote' this M tile to match N, so scale by 1/hM and translate cM to the origin. Call the scaled files DS[k] instead of LS[k]

$$DS[k] = \textbf{TwokTok[LS[k]] = TranslationTransform[\{0,0\}-cM][LS[k]]/hM;}$$

(vii) **FirstFamily[N] = Join[Table[S[k],{k,0,⟨N/2⟩, Table[DS[k], {k,0, ⟨N/2⟩ }]**

**Example: Graphics[poly/@FirstFamily[7]]**

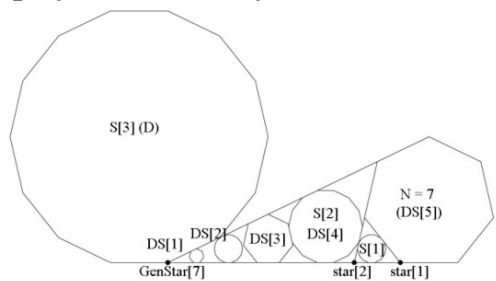

The redundant names here are useful because they highlight the connection between D and N = 7, which is the new 'M'. Basically D is playing the part of N = 14 but on the left side, so it is no surprise that N = 7 is DS[5]. The only tile that was not part of the First Family for N = 14 is S[1], but this difference grows with N, as shown earlier with N = 13.

Under the outer billiards map, the S[k] tiles of any First Family have orbits which skip k vertices on each iteration, so they have 'step-k' orbits. The DS[k] do not have constant step orbits around N = 7, but they are step-k relative to D. It would be a simple matter to extend this family with a reflection about the origin to get the First Family for N = 14. In this scenario DS[6] would be called DRight (or N = 14).

As indicated earlier, N = 7 has the potential for extended families on the edge of D – with DS[1] playing the part of the M[1] and DS[2] as D[1] and scale GenScale[7]. Using 'virtual' D[k] tiles inside of D (and matching real M[k] outside), this chain can theoretically be continued to obtain M[k]'s and D[k]'s scaled by GenScale$[7]^k$ converging to GenStar[7]. Even though N = 7 is not a 4k+1 prime, there do appear to be such chains and the families alternate between two 'templates' – where the canonical FirstFamily[7] shown above is just the template for the 'odd' generations. This alternation may be due to the fact that N = 7 has two non-trivial scales – while N = 5, 8,10 and 12 have only one non-trivial scale – and hence a well-defined fractal dimension.

**Example**: **Graphics[poly/@FirstFamily[9]]**

Since 9 is composite, the actual First Family that occurs in the outer billiards map may involve some 'mutations' of the S[k] or DS[k] tiles. This is due to shortened orbits of the center whenever (k,N) >1. On the right above, the step-3 orbit of cS[3] only 'sees' the embedded regular triangle so it has period 3 instead of 9. These shortened orbits may lead to mutations in the DS[k] or S[k] tiles. S[3] would normally be a scaled version of N = 18, but in reality it consists of two nested regular hexagons, each with center cS[3] and slightly different radii. We call these 'woven' polygons. They preserve 'half' of the dihedral symmetry of the regular S[3] and their dynamics are of interest in their own right. We will return to N = 9 in Section 4.

**Case 3- N twice-even**

The First Family for N twice-even is defined below. This is the easiest case to define because all the family members are identical except for scale and there is lateral symmetry. Of course in reality, there may be 'mutations' based on the First Family tiles as 'templates'.

**Definition**: The *First Family* of an N-gon for N twice-even is defined as follows:
Assume N is centered at the origin with height 1 and 'bottom' edge horizontal.

(i) The star points are: **star[n_]:= Table[{-Tan[k*Pi/n],-1},{k, 1, ⟨N/2⟩}]**

(ii) Each star point defines a scale using **scale[k] = $s_1/s_k$**

(iii) The 'centers' are defined to be
cS[k] = **OuterDual[star[k]]= RotationTransform[-Pi/N][star[k]*rN]**

(iv) Each S[k] will be a regular N-gon with center cS]k] and height GenScale[N]/scale[k]. (S⟨N/2⟩ will also be known as D and S[⟨N/2⟩-1] will be the 'central' tile in the First Family so we will call it C. By convention S[0] will be N)

(v) **DS[k_]:=ReflectionTransform[{1,0},cC][S[k]];**

(vi) **FFRight[N_] :=Table[S[k],{k,0, ⟨N/2⟩}]; FFLeft[N_] := Table[DS[k],{k,0, ⟨N/2⟩}];**

**FirstFamily[N_] := Join[FFRight[N], FFLeft[N]];**

**Example**: **Graphics[poly/@FirstFamily[24]]**

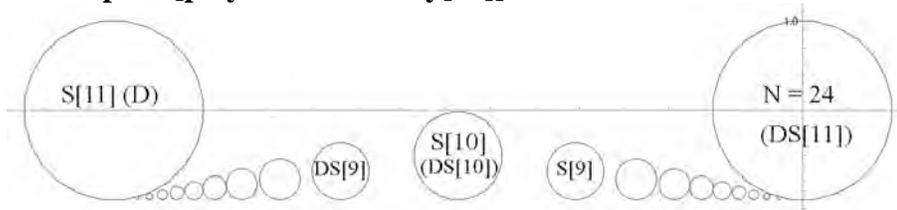

We will return to N = 24 in Section 3 and Section 4- where we find relationships between the star points and scales of inscribed polygons.

**Summary:**

The extended edges of any regular N-gon can be used to define a unique family of nested 'star-polygons'. The intersections of these extended edges are called the 'star points' of N. By symmetry there are only <N/2> 'distinct' intersection – where <N/2> is the greatest integer less than N/2. If the N-gon is in 'standard position' at the origin these 'star points' will be of the form $\pm h\{s_k,1\}$ where $s_k = \tan(k\pi/N)$ for $1 \leq k < N/2$. Each $s_k$ defines a scale and the star points and scales define a regular polygon S[k]. These S[k] form the nucleus of the First Family – which can be constructed by the algorithms in Section 2. The Families are of the form:

| N odd | N twice-odd | N twice-even |
|---|---|---|
| First Family size:⟨N/2⟩+N-3 | First Family size: N-3 | First Family size: N -3 |
| Example: N = 11 | Example: N = 18 | Example: N = 12 |
| 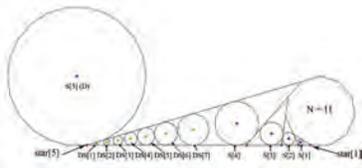 | 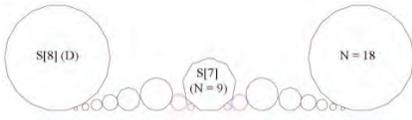 | 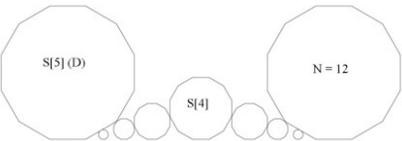 |

## Section 3: Scaling for regular polygons

For any piecewise affine map, the notion of scaling is fundamental – but very little is known about scaling under maps based on rational rotations. There may never be a uniform theory for scaling of regular polygons under the outer billiards map – because the scaling appears to be very sensitive to the algebraic complexity of N.

However every regular N-gon has a well-defined First Family of tiles with their associated scales and it is our hope that these 'canonical' scales can serve as a geometric and algebraic 'template' for the general scaling. This is true for the 'linear' and 'quadratic' cases and in general we will show that the $\varphi(N)/2$ primitive scales serve as a basis for all the 'canonical' scales. The vector space determined by these primitive scales is actually a field – namely the maximal real subfield of the cyclotomic field $\mathbb{Q}_N$. This means that there is a structure to the canonical scales – which hopefully can be exploited to obtain results about the dynamics of N.

The singularity set of a regular N-gon evolves in a recursive fashion and the star points (and associated S[k] tiles) define transitions in this evolution. Since the S[k] tiles have trivial 'step-k' dynamics, it may be possible to use this knowledge to study the evolution local to S[k].

In [H3] we have attempted to do this for N = 7 - where there are only two non-trivial primitive scales. But these scales can interact in unpredictable ways – and this potential increases with each new 'generation'. The $8^{th}$ generation for N = 7 has regions with dynamics which appear to be unrelated to any previous dynamics – but there are still embedded invariant regions with 'canonical' dynamics. For N = 11 these regions have almost vanished and there seems to be no correlation between the dynamics of 'generations'.

When N is composite, every divisor defines a 'factor' polygon and there are simple relations between the scales of N and the scales of these factors. These relations are described in the Scaling Lemma below. The factors can be considered inscribed or circumscribed as shown here for N = 24. The circumscribed form makes it clear that every factor will share its star points with N. These are classified as 'degenerate' star points of N= 24, since $s_k = \tan(k\pi/N)$ has $(k,N) >1$.

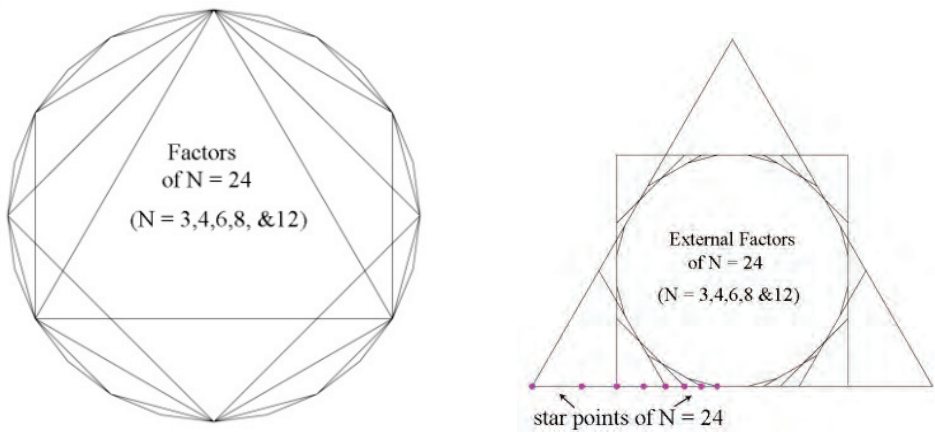

The external factors shown on the right are scaled to be in 'standard position' relative to N = 24, which means that they share the same horizontal 'base' edge and have the same height as N = 24. Since every factor shares certain edges with N = 24, they also share the corresponding star points and this implies that the corresponding scales are related by ScaleSwap.

**Lemma 4** (Scaling Lemma): *Suppose N and M are regular polygons and N/k = M, then*

*scale[j] of N/k = scale[kj] /scale[k] of N*

**Proof**: N/k can be scaled relative to N so that it is an 'external factor polygon in standard position'. To do this, center both N and N/k at the origin. N/k will have bottom edge horizontal and side s/ScaleSwap[N,N/k] where s is the side of N. The external angle of M is $2\pi k/N$ so in this position, every kth edge will coincide with an edge of N. Therefore it will share every kth star point with N and by definition the corresponding scales are related by the ratio of the sides of N and N/k, so scale[j] of N/k = (scale[kj] of N)/ScaleSwap[N,N/k]. In particular, scale[1] of M = 1 = scale[k]/ScaleSwap[N,M], so scale[k] of N is ScaleSwap[N,M]. □

**Example**: Below is N = 24 showing the external factor polygon N = 8 which corresponds to k = 3, so every third star point is aligned.

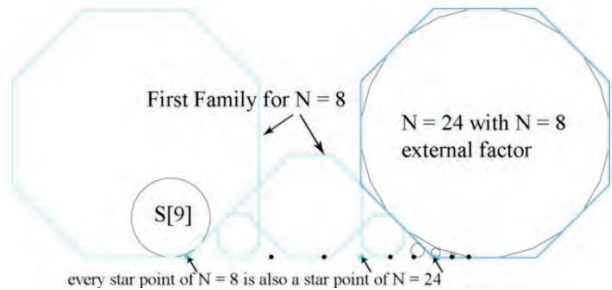

Therefore ScaleSwap[24,8] is scale[3] (of N = 24). GenScale[8] is scale[3] of N = 8 and this is equal to scale[9]/scale[3]. Because of the symmetry of the N = 8 family, scale[2] must be $\sqrt{scale[3]}$, and this implies that scale[9] of N = 24 is $scale[6]^2$/scale[3]. Because of the relationship between the scales of N = 8 and N = 24, the S[1], S[2] and S[3] tiles of N = 8 are scaled versions of S[3], S[6] and S[9] for N = 24. However there is no natural way to extend this to the rest of the First Family for N = 24, so the families appear to have very different dynamics.

Since $\varphi(24) = 8$, N = 24 has just 4 primitive scales, with indices {1,5,7,11} and we will see that all the scales are linear combinations of these 4 primitive scales. Alternatively scale[11] = GenScale[24] = $\tan^2(\pi/24)$ will generate the scales because they are all in the maximal real subfield of the cyclotomic field $\mathbb{Q}(\zeta_N)$. (See Section 4.)

**NumberFieldSignature[GenScale[24]]** = {4,0}, so the minimal polynomial for GenScale[24] has 4 real roots and this implies that every scale is of the form $a_0 + a_1 x + a_2 x^2 + a_3 x^3$ where x = GenScale[24] and the $a_i$ are rational. For example, scale[9] = $s_1/s_9$, and GenScale[24] = $s_1^2$ so

**AlgebraicNumberPolynomial[ToNumberField[Tan[Pi/24]/Tan[9*Pi/24],Tan[Pi/24]^2],x]**

yields $\dfrac{9}{64} - \dfrac{321x}{64} + \dfrac{179x^2}{64} - \dfrac{3x^3}{64} = \tan(\pi/24)/\tan(9\pi/24) = $ scale[9] $\approx 0.0545323$ when $x = s_1^2$.

The important special case of the Scaling Lemma for k = 2 will be described below. Basically this N-even Scaling Lemma says that for any regular N-gon with N-even, the scales inherit the symmetry of N.

**Lemma 5:** (N-even Scaling Lemma )

(i) In Section 1 we noted that when N is even, GenStar[N] = 1/star[1] so

**scales[N] = Reverse[GenScale[N]/scales[N]]**

Since GenScale[N], maps the scales to themselves in reverse order, the number of linearly independent scales cannot exceed N/4. We will show in Section 4 that the actual number of independent scales is always φ(N)/2, but here we are looking at the consequences of symmetry.

(ii) For N twice-odd the total number of scales is N/2 -1 which is even, so the largest number of independent scales is Floor[N/4]. The 'even' or 'odd' scales are sufficient because they map to each other under GenScale:

**EvenScales[N] = Reverse[GenScale[N]/OddScales[N]]**
**OddScales[N] = Reverse[GenScale[N]/EvenScales[N]]**

Therefore either the evens or the odds will generate the full list:
  **Scales[N] = Riffle[OddScales,EvenScales]** (Riffle interweaves them)

(iii) For N twice-even the total number of scales is still N/2 – 1 but now it is odd, so the maximum number of independent scales is N/4. The first N/4 scales (or the last N/4) form a full list of candidates for independent scales, and these sets map to each other using GenScale[N] as:

**FirstHalf[N] = Table[scale[k], {k,1, N/4}]; SecondHalf[N]= GenScale[N]/FirstHalf[N]**
  **Scales[N] = Join[FirstHalf[N], SecondHalf[N]]**

Since the 'central' S[N/4] tile of the First Family is in both lists, scale[N/4] $= \sqrt{\text{GenScale[N]}}$ .

(iv) When N is even, the corresponding N/2 –gon can be scaled to become an 'external factor polygon in standard position' where it shares all of its star points with N – as in the Scaling Lemma. It follows that ScaleSwap[N,N/2] = scale[N/(N/2)] = scale[2]. Therefore the scales of N/2 can be derived from the even scales of N as:

**Scale[N/2] = EvenScales[N]/ScaleSwap[N,N/2]**  □

**Note:** In the twice-odd case this equivalence of scales is significant because ScaleSwap[N,N/2] is GenScale[N]/GenScale[N/2] and this establishes a natural conjugacy between the First Families as shown below in Lemma 6. For N twice-even the central S[N/4] tile is an N-gon by

default, so N/2 is not embedded in the First Family of N and there is no natural connection between the First Families or the dynamics of N and N/2.

**Lemma 6:** (Twice-odd Lemma): When N is twice-odd, the scales of N/2 are related to the even scales of N by:

$$\text{scale}[k] = \text{scale}[2k]\frac{\text{GenScale}[N/2]}{\text{GenScale}[N]}$$

Therefore scale[⟨N/2⟩-1] of N is GenScale[N/2].

**Proof**: When N is even, part (iv) of the Scaling Lemma says that the scales of N/2 can be derived from the even scales of N using scale[k] = scale[2k]/scale[2] where scale[2] of N is the ratio of the sides, which is ScaleSwap[N,N/2] = tan($\pi$/N)/tan($2\pi$/N). When N is twice-odd, Lemma 2 shows that ScaleSwap[N,N/2] = GenScale[N]/GenScale[N/2].

Item (i) of the N-even Scaling Lemma says that GenScale[N]/scale[N/2-1] = scale[2] and scale[2] is GenScale[N]/GenScale[N/2]. Therefore scale[⟨N/2⟩-1] of N is GenScale[N/2]. □

The key issue in Lemma 6 is that when N is twice-odd, ScaleSwap is the ratio of the GenScales and this yields a perfect alignment of scales between N and N/2. This implies that the First Families are also related by a simple ScaleSwap as shown below.

**Corollary**: For N twice-odd, the First Families of N and N/2 are related by:

   **T[x_] := TranslationTransform[{0,0}-cS[N/2-2]][x]*ScaleSwap[N/2,N]**

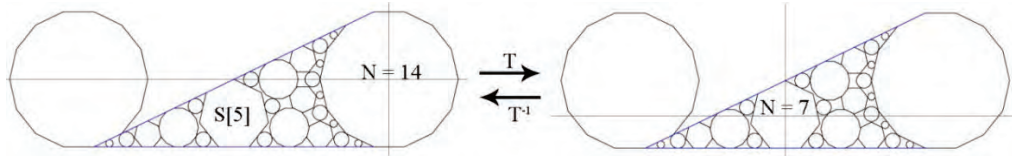

The 'half-rhombi' shown here are called the First Family 'templates' – where we have filled in the symmetric tiles for N = 14 and N = 7. It is possible to reduce all the dynamics of the outer billiards map to just these templates when N is even. This is done in Appendix B using the Digital Filter map. This implies that for N twice odd, the singularity sets of N and N/2 are related by a simple ScaleSwap and many studies of N = 7 have been done using N = 14. ScaleSwap[N/2,N] will be an algebraic integer when N is odd, so there is no loss of information. However this equivalence of 'webs' does not imply equivalent dynamics and our conjecture is that there is a conjugacy between the dynamics. Below is a comparison of the invariant regions.

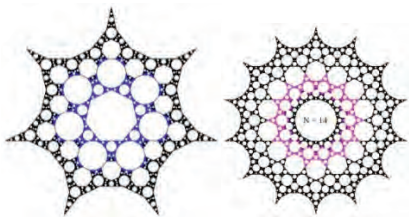

## Section 4: Algebraic complexity of star points and scales

Since the vertices of a regular n-gon with radius 1 satisfy the *cyclotomic equation* $z^n = 1$, they are algebraic and not transcendental. There are always n (complex) solutions which can be written as $\{\zeta^k : 1 \leq k \leq n\}$ where $\zeta = \exp(2\pi i/n)$. The *primitive roots* of the cyclotomic equation are those where $(k,n) = 1$, so there are $\varphi(n)$ such roots.

**Definition**: Let $\mathbb{Q}$ be the field of rational numbers. The nth *cyclotomic field* is the extension field $\mathbb{Q}(\zeta_n)$ where $\zeta_n$ is a primitive root of $z^n = 1$ – which we usually take to be $\zeta = \exp(2\pi i/n)$.

Since $\zeta$ is an algebraic number, $\mathbb{Q}(\zeta)$ is a number field. It is called the 'splitting field' of $z^n = 1$. All of the primitive roots share the same minimal polynomial, so they are called Galois conjugates. This minimal polynomial is defined below:

**Definition**: The nth cyclotomic polynomial is defined to be

$$\Phi_n(x) = \prod_{1 \leq k \leq n,\ (k,n)=1} (x - \zeta^k)$$

This polynomial has degree $\varphi(n)$ and can be shown to be irreducible, so $\mathbb{Q}(\zeta)$ has dimension $\varphi(n)$ over $\mathbb{Q}$. $\Phi_n(x)$ will be monic over $\mathbb{Z}$ and have all integer coefficients, so each primitive root is an algebraic integer and the set $\{1, \zeta, \zeta^2, \ldots, \zeta^{\varphi(n)-1}\}$ forms an integral basis. So $\mathbb{Q}(\zeta) = \mathbb{Q}[\zeta]$ since every element $\alpha \in \mathbb{Q}(\zeta)$ is a linear combination of the primitive $\zeta^k$ with coefficients $\mathbb{Q}$. If the coefficients are restricted to integers, then $\alpha \in \mathbb{Z}[\zeta]$ – the ring of integers.

For any field F and extension K, the set of all automorphisms of K that leave F fixed, form a group known as Aut(K/F). When K is the splitting field of a separable polynomial, it is called a *Galois extension* and Aut(K/F) is called the *Galois group* of K/F and written Gal(K/F). Any automorphism of $\mathbb{Q}(\zeta)$ which fixes $\mathbb{Q}$, will map the roots to themselves so $\mathrm{Gal}(\mathbb{Q}(\zeta)/\mathbb{Q})$ is $\{\sigma_k : \zeta \to \zeta^k : (k,n) = 1\}$. We will write this as $G_n$. Therefore $|G_n| = \varphi(n)$ and it is isomorphic to $\mathbb{Z}^x_n$ - the integers mod n under multiplication. $\mathbb{Z}^x_n$ consists of the 'units' of $\mathbb{Z}_n$ – where $\mathbb{Z}_n$ is the additive group of integers mod n.

For a Galois extension E/F, the fundamental theorem of Galois theory says that there is a 1-1 correspondence between the subgroups H of Gal(E/F) and the subfields K of E (which are Galois by default). Moreover the degree of E over K, [E:K], is equal to |H|. For example if K is an intermediate extension $\mathbb{Q} \subset K \subset \mathbb{Q}(\zeta)$, it will also be Galois. Then $|G_K| = [\mathbb{Q}(\zeta): K]$ and since $[\mathbb{Q}(\zeta):K][K:\mathbb{Q}] = [\mathbb{Q}(\zeta):\mathbb{Q}] = \varphi(n)$, this will yield the degree of K over $\mathbb{Q}$. Therefore to find the degree of K over $\mathbb{Q}$, it is sufficient to determine $G_K$ -which consists of the $\sigma_k$ that leave K fixed.

There are three extensions K of immediate interest – corresponding to $\cos(2\pi k/n)$, $\sin(2\pi k/n)$ and $\tan(2\pi k/n)$.

(i) $\mathrm{Cos}(2k\pi/n) = (\zeta^k + \zeta^{-k})/2$ whenever $(k,N) = 1$, so each primitive $\zeta$ can be grouped with its complex conjugate to obtain a real result in $\mathbb{Q}(\zeta)$. Setting $\xi_c = \zeta + \zeta^{-1}$, $K_c = \mathbb{Q}(\xi_c) \subset \mathbb{Q}(\zeta)$.

Complex conjugation in $\mathbb{Q}(\zeta)$ is always an automorphism of order 2, so $|G_\mathbb{K}| = [\mathbb{Q}(\zeta): \mathbb{Q}(\xi_c)] = 2$ and since $[\mathbb{Q}(\zeta):\mathbb{Q}(\xi_c)] [\mathbb{Q}(\xi_c), \mathbb{Q}] = [\mathbb{Q}(\zeta):\mathbb{Q}] = \varphi(n)$, $[\mathbb{Q}(\xi_c):\mathbb{Q}] = \varphi(n)/2$. Therefore the minimal polynomial for $\cos(2k\pi/n)$ has degree $\varphi(n)/2$ over $\mathbb{Q}$ whenever $(k,N) = 1$.

(ii) $\sin(2\pi k/n) = (\zeta^k - \zeta^{-k})/2i$ whenever $(k,N) = 1$, so the matching extension field is $K_s = \mathbb{Q}(\xi_s)$ where $\xi_s = \zeta - \zeta^{-1}$. Clearly $\mathbb{Q}(\xi_s) \subset \mathbb{Q}(\zeta)$. To find the dimension of $K_s$ over $\mathbb{Q}(\zeta)$, note that $\sigma_1$ is the only possible automorphism when n is not divisible by 4, so $|G_K| = 1$. But when n is divisible by 4, $i \in \mathbb{Q}(\zeta)$ and the extra symmetry allows for $|G_K| = 2$ when $n = 0$ mod 8 and $|G_K| = 4$ when $n = 4$ mod 8. This last case occurs because $i \in \mathbb{Q}(\xi_s)$ iff $n = 4$ mod 8. This yields the three distinct dimensions over $\mathbb{Q}$ in the theorem below.

(iii) $\tan(2\pi k/n) = (\zeta^k - \zeta^{-k})/(\zeta^k + \zeta^{-k})i$ whenever $(k,N) = 1$ so the matching extension field is $K_t = \mathbb{Q}(\xi_t)$ where $\xi_t = \zeta - \zeta^{-1}/(\zeta + \zeta^{-1})$. The Galois group of $K_t$ depends on n just like sin. When n is not divisible by 4, $|G_K|$ is again 1, but now the extra symmetry of n divisible by 4 yields $|G_K| = 2$ when $n = 4$ mod 8 and $|G_K| = 4$ when $n = 0$ mod 8. This last case occurs because now $i \in \mathbb{Q}(\xi_t)$ iff $n = 0$ mod 8.

**Complexity Theorem** ([R],[Ca]) Suppose $n > 2$ and $(k,n) = 1$. In the sine and tangent cases assume $n \neq 4$. Then the algebraic degree over $\mathbb{Q}$ of

(i) $\cos(2k\pi/n)$ is $\varphi(n)/2$

(ii) $\sin(2k\pi/n)$ is $\begin{cases} \varphi(n) & \text{if } 4 \text{ does not divide } n \\ \varphi(n)/2 & \text{if } n = 0 \text{ mod } 8 \\ \varphi(n)/4 & \text{if } n = 4 \text{ mod } 8 \end{cases}$

(iii) $\tan(2k\pi/n)$ is $\begin{cases} \varphi(n) & \text{if } 4 \text{ does not divide } n \\ \varphi(n)/2 & \text{if } n = 4 \text{ mod } 8 \\ \varphi(n)/4 & \text{if } n = 0 \text{ mod } 8 \end{cases}$

The case of $\sin(2k\pi/N)$ is credited to Lehmer [L] (1933), but he based his derivation on the cosine case using the identity $\sin(2\pi k/n) = \cos(2\pi(4k-n)/4n)$. Lehmer noted that since $(k,n) = 1$, the fraction $(4k-n)/4n$ is either in lowest terms or "may be reduced to an equal fraction in its lowest terms with the denominator 2n or n". However when n is 4 mod 8, the reduced denominator will be $(n/2)$. Setting $k = 1$ and $n = 8m + 4$, $(4 - n)/4n = - m/(n/2)$ so when $n = 12$, $(4-12)/48 = -8/48 = -1/6$. This implies that the degree of $\sin(2\pi/12)$ is 1 and not 2 as stated by Lehmer.

Niven [N] (1956) quotes a corrected version of Lehmer's result for sine and then extends it to the tangent case by showing that $\mathbb{Q}(\tan(2\pi k/n)) = \mathbb{Q}(\sin(2\pi k/n))$ when 4 does not divide n, $\mathbb{Q}(\tan(2\pi k/n)) = \mathbb{Q}(\cos(2\pi k/n))$ when $n = 4$ mod 8, and $\mathbb{Q}(\tan(2\pi k/n)) = \mathbb{Q}(\cos(4\pi k/n))$ otherwise.

Calcut [Ca] gives a thorough treatment of all three cases from first principles using Galois Theory and later reproves the tangent result using an analogue to the Chebyshev polynomials of the first kind for cosine. Stillwell [St] notes that these basic polynomials were known to John Bernoulli as early as 1702, but it was Bernoulli's student Euler who understood their significance in the development of complex logarithms and exponentials.

**Definition**: The algebraic complexity of a regular N-gon is the degree of the minimal polynomial for $\cos(2\pi/N)$ – which is $\varphi(N)/2$.

So N = 5,8,10 and 12 have 'quadratic' complexity, while N =7, N = 14 and N = 9 and N = 18 are 'cubic' and N = 15, 16 and 24 are 'quartic'.

In Mathematica, **Cyclotomic[N,x]** gives $\Phi_N(x)$ and **MinimumPolynimial[Cos[2Pi/N]][x]** gives the minimal polynomial for $\cos(2\pi/N)$.

The Complexity Theorem simplifies when $\tan(2\pi k/N)$ is replaced by $\tan(k\pi/N)$ because the only issue now is whether 4|N or not

**Corollary** (Complexity Theorem – Part 2): If n > 2 and (k,n) = 1, the degree of $\tan(k\pi/n)$ over $\mathbb{Q}$ is $\varphi(n)$, except when n is divisible by 4, the degree is $\varphi(n)/2$.

**Example**: **Table[{n, EulerPhi[n], Exponent[MinimalPolynomial[Cos[Pi/n]][x], x], Exponent[MinimalPolynomial[Tan[Pi/n]][x], x], {n, 3, 20}]** (The matching Cotangent sequence is https://oeis.org/A089929.)

| N | 3 | 4 | 5 | 6 | 7 | 8 | 9 | 10 | 11 | 12 | 13 | 14 | 15 | 16 | 17 | 18 | 19 | 20 |
|---|---|---|---|---|---|---|---|---|---|---|---|---|---|---|---|---|---|---|
| Φ(N) | 2 | 2 | 4 | 2 | 6 | 4 | 6 | 4 | 10 | 4 | 12 | 6 | 8 | 8 | 16 | 6 | 18 | 8 |
| Degree of Cos[2Pi/N] | 1 | 1 | 2 | 1 | 3 | 2 | 3 | 2 | 5 | 2 | 6 | 3 | 4 | 4 | 8 | 3 | 9 | 4 |
| Degree of Tan[Pi/N] | 2 | 1 | 4 | 2 | 6 | 2 | 6 | 4 | 10 | 2 | 12 | 6 | 8 | 4 | 16 | 6 | 18 | 4 |

For a given regular N-gon, there are always $\varphi(N)/2$ primitive $s_k$ and when 4 does not divide N, $T_N[x]$ will be an even function with $\varphi(N)$ roots, and the roots will consist of the primitive $s_k$ and their negatives. Therefore the primitive $s_k$ are Galois conjugates. When 4|N, $T_N[x]$ will be odd with degree $\varphi(N)/2$, and there is no longer a guarantee that the primitive $s_k$ will be among the roots, but it is easy to restore the missing symmetry by including the negative of each root. This can be done using $T_N[x]$ in conjunction with $T_N[-x]$. In Lemma 8 we will use the Root Lemma below to show that these statements are valid.

**Lemma 7** (Root Lemma): If $s_k = \tan(k\pi/N)$ exists then the roots of the minimal polynomial for $\tan(k\pi/N)$ are $s_{jk}$ where j is any integer such that (j,N) = 1. But when 4|N, j must also be 1 mod 4.

(Note that k does not have to be positive. We will only need the cases of k = 1 and k = -1 below.)

**Example** (N = 8): The primitive $s_k$ of N = 8 are $s_1 = \sqrt{2} - 1$ and $s_3 = \sqrt{2} + 1$ as shown below in black. It is clear that both $T_8[x]$ and $T_8[-x]$ are needed to yield the two primitive roots (but as generators $\mathbb{Q}(s_1) = \mathbb{Q}(-s_1)$, and they are both degree 2, so either would suffice for generating the primitive roots.) In the Root Lemma, setting k = 1 will yield $T_8[x]$ in blue. Since j must be 1 mod 4, the only possible j values are j = 1 and j = 5 = -3 mod 8, so the indices are kj = 1 and kj = -3.

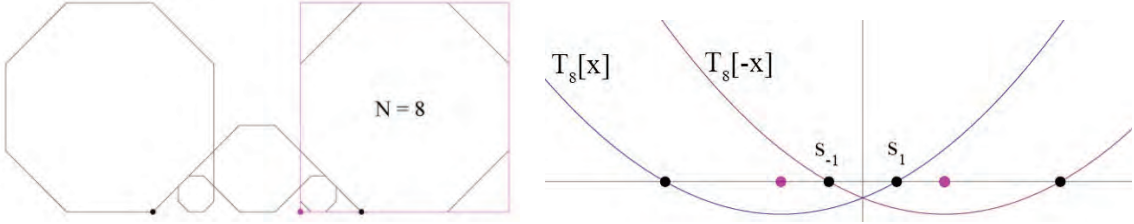

Normally setting k = -1 would yield the same roots, but when 4|N, it will yield the matching negative roots for $T_8[x]$. Setting k = -1 gives $T_8[-x]$ in red, and once again the only possible j values are j = 1 and j = 5, so the indices are kj = -1 and kj = -5 = 3 mod 8. This yields the two remaining roots. Therefore the full set of 'primitive' indices is {-3,-1,1,3}. The magenta star[2] corresponds to $s_2 = \tan(2\pi/8) = 1$ and the Basis Lemma below shows that it must be in $\mathbb{Q}(s_1)$.

**Lemma 8**: The primitive $s_k$ of a regular N-gon are the positive Galois conjugates of $s_1 = \tan(\pi/N)$ or $-s_1 = \tan(-\pi/N)$. There are always $\varphi(n)/2$ of these positive conjugates.

**Proof**: Note that when k = 1 or k = -1 in the Root Lemma, the matching $s_1 = \tan(\pi/N)$ and $s_{-1} = \tan(-\pi/N)$ always exist.

(i) Suppose that 4 does not divide N. Since j is mod N, we can assume j < N. The Root Lemma says that the set of Galois conjugates of $s_1$ is
$$G_1 = \{s_{jk} : k = 1 \text{ and } kj < N, (j,N) = 1\} \supset S = \{s_k : 1 \leq k < N/2, (k,N) = 1\}$$
so every primitive $s_k$ is a Galois conjugate of $s_1$. When (N,j) = 1, then (N,-j) = (N,N-j) is also 1, so every $s_{jk}$ in $G_1$ has a matching $s_{-jk}$. This implies that setting k = -1 would yield the same roots, so $T_N[x] = T_N[-x]$. The set S corresponds to the positive indices with j < N/2 – so the primitive roots are positive, and there are $\varphi(N)/2$ such roots.

(ii) Suppose that 4|N. The Root Lemma says that the Galois conjugates of $s_1$ and $-s_1$ are obtained by setting k = 1 or k = -1, so the combined conjugates are
$$G_2 = \{k = 1 \text{ or } k = -1 \text{ and } s_{jk} : kj < N, (j,N) = 1 \text{ and } ((j, 4) = 1 \text{ or } (-j,4) = 1)\}$$
Since N is twice even (j,N) = 1 implies j = 1 or 3 mod 4, so these can be dropped above
$$G_2 = \{k = 1 \text{ or } k = -1 \text{ and } s_{jk} : kj < N, (j,N) = 1\} \supset S = \{s_k : 1 \leq k < N/2, (k,N) = 1\}$$
Now the primitive $s_k$ on the right are spread between the two cases of k = 1 and k = -1 but there will always be an equal number of positive and negative indices since any valid j for k = 1 will remain valid for k = -1. In fact $T_N[x]$ (k = 1) = $-T_N[-x]$( k = -1). There are $\varphi(N)/2$ combined roots.
□

**Corollary**: For any regular N-gon, $\mathbb{Q}(s_1)$ generates all the primitive $s_k$ and is always degree $\varphi(N)/2$.

**Proof**: If N is not divisible by 4, Lemma 8 says that all the primitive $s_k$ are Galois conjugates of $s_1$, so they are in $\mathbb{Q}(s_1)$ (along with their negatives). If N is divisible by 4, the primitive $s_k$ are roots of $s_1$ or $-s_1$, and $\mathbb{Q}(s_1) = \mathbb{Q}(-s_1)$. Therefore $\mathbb{Q}(s_1)$ is always degree $\varphi(N)/2$. □

Note that this does not imply that the primitive $s_k$ are independent – just that the set $\{1, s_1, s_1^2, \ldots\}$ is independent. Also $s_1$ is not in $\mathbb{Q}(\zeta_N)$ unless $4|N$, so we will often use the 'surrogate' $is_1$.

**Example**: N = 15 has $\varphi(15) = 8$ and there are 7 'star' points $s_k$ with 4 primitive.

**NumberFieldSignature[Tan[Pi/15]]** yields {8,0}. This is an 'overkill' because the 8 real roots are of the form $\pm s_k$ for (k,n) = 1. Lemma 8 guarantees that all the primitive roots will be in $\mathbb{Q}(s_1)$, and the Basis Lemma below takes this one step further to show that all the degenerate $s_k$ are linear combinations of the primitive $s_k$, so they are also in $\mathbb{Q}(s_1)$. Here is an example with $s_5 = \tan(5\pi/15) = \sqrt{3}$, but in order to work entirely inside $\mathbb{Q}(\zeta_{15})$, we will solve for $is_5$ using $\mathbb{Q}(is_1)$,

**AlgebraicNumberPolynomial[ToNumberField[I*Tan[Pi/15], I*Tan[Pi/15]] ,x]**

yields $\dfrac{395x}{32} + \dfrac{3075x^3}{32} + \dfrac{2297x^5}{32} + \dfrac{25x^7}{32} = i\sqrt{3}$ where $x = i\tan(\pi/15)$. Since these are odd powers of $ix$, the $i$ terms will appear as +/- and hence $i$ can be cancelled to yield

$$\frac{395}{32}Tan[\frac{\pi}{15}] - \frac{3075}{32}Tan[\frac{\pi}{15}]^3 + \frac{2297}{32}Tan[\frac{\pi}{15}]^5 - \frac{25}{32}Tan[\frac{\pi}{15}]^7 = \sqrt{3}$$

Of course this same result can be obtained directly using $\tan(\pi/15)$ as generator – but in general it is more efficient to do calculations inside $\mathbb{Q}(\zeta_{15})$.

The minimal polynomial for $\tan(\pi/15)$ is **T[x_] := MinimalPolynomial[Tan[Pi/15]][x]** $= x^8 - 92x^6 + 134x^4 - 28x^2 + 1$. Using FirstFamily.nb with npoints = 15, Mathematica automatically finds the First Family, star polygons and star points. The primitive star points correspond to the units in $\mathbb{Z}_{15}$. To find them use the GCD function or **K = Table[k*DirichletCharacter[k, 1, 15], {k, 1, 15}];** followed by **K = Delete Cases[K, 0]** Then K = {1,2,4,7,8,11,13,14}. These will be the k-indices of the roots of T[x] – so the primitive $s_k$ will have indices 1,2,4,7 and we just show indices 1,2 and 4 on the right below and also the matching left side $s_k$ with indices -1(14),-2(13),-4(11) The complete set of conjugates is **S = Table[{Tan[K[[j]]*Pi/15],0} {j,1,8}].** The degenerate $s_k$ are shown in magenta below.

**Show[Plot[T[x], {x, -1.3, 1.3}], ListPlot[{S,T} PlotStyle -> {Black,Magenta}]}]** (on the right)

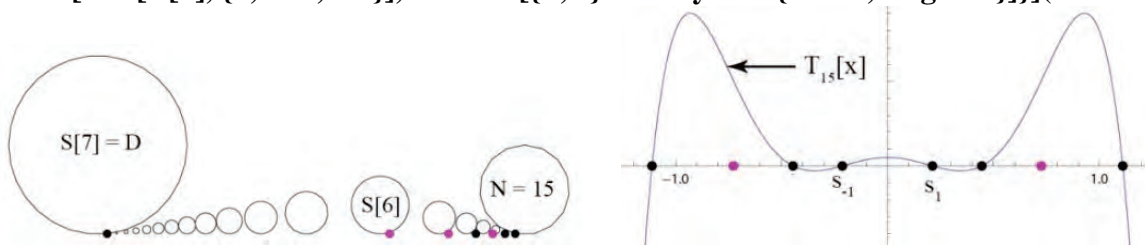

Note: If we set h = 1 for the First Family, the first coordinate of the star points will match the right side, but the second coordinate on the right is 0 by convention – while the left side is –h.

**Lemma 9** (Basis Lemma) For a given N-gon we have defined $s_k = \tan(k\pi/N)$ for $1 \leq k < N/2$. These 'star points' are classified as 'primitive' if $(k,N) = 1$ and 'degenerate' otherwise. Here we show that every degenerate $s_k$ can be written as a linear combination of the primitive $s_i$. (For convenience we show the equivalent result that every degenerate $is_k$ is a linear combination of the primitive $is_k$.)

**Proof**: Since $i\tan(\theta) = (e^{i2\theta}-1)/(e^{i2\theta}+1)$, $is_k = (e^{2k\pi i/N}-1)/(e^{2k\pi i/N}+1) = (\zeta^k-1)/(\zeta^k+1)$ where $\zeta = \exp(2\pi i/N)$, so $is_k \in \mathbb{Q}(\zeta)$ and the only difference between the primitive and degenerate $s_k$ is whether $(k,N) = 1$. (Note that this does not show that the $s_k$ are in $\mathbb{Q}(\zeta)$ because $i \in \mathbb{Q}(\zeta)$ only when $4|N$.)

The $\varphi(N)$ primitive roots of $\mathbb{Q}(\zeta)$ can be grouped in pairs to generate $\mathbb{Q}(\zeta)^+ = \mathbb{Q}(\zeta+\zeta^{-1})$. This is called the maximal totally real subfield of $\mathbb{Q}(\zeta)$. It plays an important role in algebraic number theory and $\mathbb{Q}(\zeta)$ is always a quadratic extension of $\mathbb{Q}(\zeta)^+$. Since $\mathbb{Q}(\zeta)^+ = \mathbb{Q}(\zeta) \cap \mathbb{R}$, the $is_k$ all lie in the complement of $\mathbb{Q}(\zeta)^+$, $\mathbb{Q}(\zeta) \cap \mathbb{R} i$, which is sometimes written as $\mathbb{Q}(\zeta)^-$.

The Complexity Theorem shows that $\mathbb{Q}(\zeta)^+$ has dimension $\varphi(N)/2$ over $\mathbb{Q}$, so $\mathbb{Q}(\zeta)^-$ also has dimension $\varphi(N)/2$ – as a vector space. Since $is_k \in \mathbb{Q}(\zeta)^-$ for $1 \leq k < N/2$, and the set S of primitive $is_k$ is independent with $|S| = \varphi(N)/2$, they form a basis for $\mathbb{Q}(\zeta)^-$. □

A similar argument can be applied to the cotangent to obtain the 'dual' result that every 'degenerate' $c_k = \cot(k\pi/n)$ can be written as a linear combination of the primitive $c_k$. This is done in [G2].

It should be clear that $is_1 = i\tan(\pi/N)$ can actually generate both $\mathbb{Q}(\zeta)^+$ and $\mathbb{Q}(\zeta)^-$ - using the even and odd powers respectively. This means that $\tan^2(\pi/N)$ will generate $\mathbb{Q}(\zeta)^+$ as we implied earlier. For $\mathbb{Q}(\zeta)^-$ a basis could be $\{is_1^k\}$ for k odd and $k < \varphi(N)$ as shown above for $N = 15$, but the primitive $is_k$ are a more practical choice of basis.

As a vector space, every element of $\mathbb{Q}(\zeta)^-$ has the form:

$$\sum_{1 \leq k < N/2, (k,n) = 1} a_i is_k = \sum_{1 \leq k < N/2, (k,n) = 1} a_i i Tan[k\pi/N] \text{ for } a_i \in \mathbb{Q}$$

Therefore if $M = N/d$, $\mathbb{Q}(\zeta_N) \supset \mathbb{Q}(\zeta_M)$ and there must be rational $a_i$ such that

$$\tan(j\pi/M) = \tan(jd\pi/N) = \sum_{1 \leq k < N/2, (k,n) = 1} a_i s_k = \sum_{1 \leq k < N/2, (k,n) = 1} a_i Tan[k\pi/N]$$

Note that since the $a_i$ can be negative, the indices k can be chosen negative (mod N) if desired. In general there are no elementary formulas that will yield these coefficients, so no one knows what $\mathbb{Q}(\zeta)^-$ really looks like for large N  Most results are derived from Dirichlet's $\mathcal{L}$-series result that $\mathcal{L}(s, \mathcal{X}) = \sum_1^\infty \frac{\mathcal{A}(n)}{n^s}$ is not zero for $s = 1$ and $\mathcal{X}$ a principle character mod N. This series is naturally more applicable to the cot case than the tan series, so these are the strongest results. Theorem 1 of [G2] gives a manageable formula for the coefficients in the cotangent case (see the second example below). The problem is that the known tan results depend on relating the character sets of tan and cot. Theorem 2 of [G2] -presents a "sort of closed formula" for the tan coefficients. For more on this issue see [Ch], [G1], [G2], [B] and [Ha]

**Example**: For a case like $N = 18$, known summation formula such as those found in [P] can be used to find the coefficients. For example $\tan(\pi/18) + \tan(\pi/18 + \pi/3) + \tan(\pi/18 + 2\pi/3) = \tan(6\pi/18)$. Note that all the factors on the left side are primitive in $N = 18$, with k values 1, 7 and 13 (-5), so $\tan(\pi/18) + \tan(7\pi/18) - \tan(5\pi/18) = \sqrt{3}$. By exploiting the connection between even and odd $s_k$ for N twice-odd, this can be written as $\cot(4\pi/9) + \cot(\pi/9) - \cot(2\pi/9) = \sqrt{3}$. But this does not immediately solve the problem of writing this degenerate $s_3 = \tan(3\pi/9)$ of $N = 9$ in terms of the primitive $s_k$ of $N = 9$. More on this below.

**Example**: $N = 1001$ and divisor $d = 13$ (from [G2]). This is a cotangent example and the goal is to find degenerate $\cot(j13\pi/1001) = \cot(j\pi/77)$ using the $\varphi(1001)/2 = 360$ primitive $s_k$ of $N = 1001$. The 360 coefficients will depend on j in a 'mod 13' fashion, so when $k = \pm j \bmod 13$, the matching $a_k$ have absolute value given by $a(j)$ where $a(1), a(2)..,a(6)$ are given by:

$a(1) = 3 \cdot 9 \cdot 1951 \cdot 6007/M$, $a(2) = 190871389/2M$, $a(3) = 19 \cdot 275929/M$, $a(4) = 19 \cdot 53 \cdot 103 \cdot 131/M$, $a(5) = 3^4 \cdot 31 \cdot 1951/2M$, $a(6) = 4373 \cdot 27793/2M$ – where $M = 5^2 \cdot 13^2 \cdot 61 \cdot 181 \cdot 1117$.

Patterns of this kind persist for power-series solutions and typically the p|N cases are hard: **AlgebraicNumberPolynomial[ToNumberField[Tan[Pi/11],Tan[Pi/77]],x]** will (eventually) return 30 coefficients which form 15 +/- pairs as above, each with common denominators of $M = 2^{53}$, 2M, 4M or 8M. As expected, the numerators are all divisible by 7. For example the coefficient of $x^{29}$ is $7^4 \cdot 11 \cdot 2351 \cdot 994087 \cdot 6038425783052745701/M$.

These results about star points will enable us to prove non-trivial results about the (canonical) scales for a regular N-gon. Recall that the scales are defined as $\{s_1/s_k: s_k = \tan(k\pi/N)$ for $1 \leq k < N\}$.

**Definition**: The *primitive scales* of a regular N-gon are the scales of the form $s_1/s_k$ for $s_k$ primitive.

Since $s_1$ is always primitive, the primitive scales are in 1-1 correspondence with the primitive star points so there are $\varphi(N)/2$ primitive scales. Even though the primitive $s_k$ are independent, this does not imply that the primitive scales are independent, and we will prove this below using the fact that the dual star points are linearly independent as shown by Siegel and Chowla.

In Section 1, we defined the dual star points of a regular N-gon to be $\{r_k = \cot(k\pi/N), 1 \leq k < N/2\}$ with the matching *dual scales* of the form $r_1/r_k = s_k/s_1$ and we noted that when $s_k$ is primitive, then $1/s_k$ is primitive with respect to cot, because the k indices are identical between tan and cot. Therefore the dual primitive scales are simply the inverse of the primitive scales.

**Theorem 2:** (complexity of scales): For any regular N-gon, the set of primitive scales is linearly independent and forms a basis for all the (canonical) scales. Therefore the complexity of the scales is $\varphi(N)/2$.

**Proof:** The scales of a regular N-gon are $\{t_k = s_1/s_k$ for $1 \leq k < N/2\}$ and the primitive scales are $T = \{t_k, (k,N) = 1\}$ so $|T| = \varphi(N)/2$

(i) To show that the primitive scales are linearly independent, suppose that $\sum_{1 \leq k < N/2} a_k t_k = 0$ with $(k,N) = 1$ and $a_i \in \mathbb{Q}$. Then $\frac{1}{s_1} \sum_{1 \leq k < N/2} a_k t_k = \sum_{1 \leq k < N/2} a_i r_k = 0$ contrary to the independence of the primitive dual star points $\{r_k, (k,N) = 1\}$.

(ii) To show that these primitive scales form a basis for all scales, suppose that $t_j = s_1/s_j$ has $s_j$ degenerate. Then the dual star point $r_j = 1/s_j$ is also degenerate but it can be written as $r_j = \sum_{1 \leq k < N/2} a_i r_k$ for $r_k$ primitive dual star points. Therefore $t_j = s_1 \sum_{1 \leq k < N/2} a_i r_k = \sum_{1 \leq k < N/2} a_i t_k$ for $(k,N) = 1$. □

**Corollary** (complexity of dual scales) For any regular N-gon, the set of primitive dual scales is linearly independent and forms a basis for all the dual scales. Therefore the complexity of the dual scales is $\varphi(N)/2$. (The proof is the dual of the proof Theorem 2 – exchanging star points with duals.)

Since the scales are of the form $is_1/is_k$ and $is_k \in \mathbb{Q}(\zeta)$, the scales are in $\mathbb{Q}(\zeta)$ and they are real so they are in the maximal real subfield $\mathbb{Q}(\zeta)^+$ - which has degree $\varphi(N)/2$. Therefore the primitive scales form a basis for $\mathbb{Q}(\zeta)^+$ and the dual scales are also in $\mathbb{Q}(\zeta)^+$, so primitive scales are units in $\mathbb{Q}(\zeta)^+$. This means that they form a group which has finite index in the full unit group. More on this later.

We noted earlier that the traditional generator of $\mathbb{Q}(\zeta)^+$ is $\lambda_N = 2\cos(2\pi/N)$ and based on the half-angle formula of Diophantus

$$\cos(2\pi/N) = \frac{1 - \tan^2(\pi/N)}{1 + \tan^2(\pi/N)}$$

it follows that $\tan^2(\pi/N)$ will also serve a generator for $\mathbb{Q}(\zeta)^+$ and when N is even $\tan^2(\pi/N)$ is GenScale[N], so it is a unit when GenScale[N] is primitive.

**Lemma 10**: GenScale[N] is primitive except when N is 2 mod 4.

**Proof**: By definition scale[j] is primitive iff $s_j$ is primitive and this is true iff $(j,N) = 1$. By definition GenScale[N] = scale[<N/2>] so it is primitive iff $(<N/2>,N) = 1$. Recall that <N/2> is the largest integer $k < N/2$, so $N-1 \leq 2k < N$. Since $(N,N+1) = 1$, $(<N/2>,N)$ is either 1 or 2.
(i) When N is even <N/2> = (N-2)/2 so $(\frac{N-2}{2}, N) = (\frac{2M-2}{2}, 2M) = (M-1, 2M)$. If N is twice-odd, M-1 is even, so $(<N/2>,N) = 2$. When N is twice-even, M-1 is odd, so $(<N/2>,N) = 1$.
(ii) When N is odd, <N/2> = (N-1)/2 is even, so $(<N/2>,N) = 1$. □

Therefore by Theorem 2, GenScale[N] is a unit in $\mathbb{Q}(\zeta)^+$, except when N = 2 mod 4 and in his case GenScale[N/2] is a unit – and of course $\mathbb{Q}(\zeta_N) \approx \mathbb{Q}(\zeta_{N/2})$. In addition the Twice-odd Lemma shows that GenScale[N/2] is the penultimate scale of N, so it is a natural choice as surrogate for GenScale[N] when N is 2 mod 4. (Many authors simply omit this case from consideration, but we believe that that there is something to be learned by including this case.)

**Corollary**: For any regular N-gon, the primitive scales are a unit basis for $\mathbb{Q}(\zeta)^+$, and GenScale[N] is a unit generator for $\mathbb{Q}(\zeta)^+$, except when N is 2 mod 4, in which case GenScale[N/2] is a unit generator.

**The Scaling Conjecture for Regular N-gons**: For any regular N-gon, all of the scales needed to describe the dynamics of the outer-billiards map can be chosen from $\mathbb{Q}(\zeta)^+$, and hence these scales are $\mathbb{Q}$-linear combinations of the canonical primitive scales determined by the First Families. These primitive scales are algebraic units and among them are GenScale[N] (or GenScale[N/2] if N is 2 mod N). These will always be generators for $\mathbb{Q}(\zeta)^+$ (See Appendix C for a more specific version of this conjecture – with examples.)

Progressing from $\tan(2\pi/N) \to \tan(\pi/N) \to i\tan(\pi/N)$ reduces the cases from 3 to 2 to 1 and in the last case, $i\tan(\pi/N)$ can be split into its real and imaginary parts to generate $\mathbb{Q}(\zeta)^+$ and the complement $\mathbb{Q}(\zeta)^-$.

**Example**: N = 27 has <N/2> = 13 and scale[13] is GenScale[27]. Since $\varphi(27) = 18$, there are 9 primitive scales and the maximal real subfield $\mathbb{Q}(\zeta)^+$ is degree 9.

**NumberFieldSignature[Tan[Pi/27]^2]** = **NumberFieldSignature[GenScale[27]]**={9,0} and $\tan^2(\pi/27)$ and GenScale[27] are both integral generators of $\mathbb{Q}(\zeta)^+$, but only GenScale[27] = $\tan(\pi/27)\cdot\tan(\pi/54)$ is a unit, so it would be a natural choice of generator when N is odd. To make a comparison, we will solve for degenerate scale[9] = $\tan(\pi/27)/\tan(9\pi/27) = \tan(\pi/27)/\sqrt{3}$ usng both of these generators. (To find scale[9] in terms of $\tan^2(\pi/27)$ or $\tan(\pi/27)\cdot\tan(2\pi/27)$ sounds easy, but it is not because $\tan(\pi/27)$ is not in $\mathbb{Q}(\zeta_{27})^+$ (or even in $\mathbb{Q}(\zeta_{27})$)).

**AlgebraicNumberPolynomial[ToNumberField[scale[9],Tan[Pi/27]^2],x]** =

$$\frac{6435}{32768} - \frac{88311 x}{8192} + \frac{420353 x^2}{4096} - \frac{2535207 x^3}{8192} + \frac{5678187 x^4}{16384} - \frac{3843943 x^5}{24576} + \frac{53727 x^6}{2048} - \frac{11253 x^7}{8192} + \frac{143 x^8}{32768}$$

**AlgebraicNumberPolynomial[ToNumberField[scale[9],GenerationScale[27]],x]** =

$$\frac{35}{384} - \frac{671 x}{192} + \frac{171 x^2}{64} + \frac{871 x^3}{64} + \frac{109 x^4}{16} - \frac{253 x^5}{192} - \frac{241 x^6}{192} - \frac{25 x^7}{192} + \frac{5 x^8}{384}$$

As indicated in the Basis Lemma, $i\tan(\pi/27)$ is in the vector-space complement of $\mathbb{Q}(\zeta_{27})^+$, which we call $\mathbb{Q}(\zeta_{27})^-$, along with all the other 'imaginary' star points $is_k$, but $i\tan(\pi/27)$ it is still capable of generating scale[9] - in the guise of $\tan^2(\pi/27)$ (with proper sign changes).

**AlgebraicNumberPolynomial[ToNumberField[scale[9],I*Tan[Pi/27]],x]** =

$$\frac{6435}{32768} + \frac{88311 x^2}{8192} + \frac{420353 x^4}{4096} + \frac{2535207 x^6}{8192} + \frac{5678187 x^8}{16384} + \frac{3843943 x^{10}}{24576} + \frac{53727 x^{12}}{2048} + \frac{11253 x^{14}}{8192} + \frac{143 x^{16}}{32768}$$

**Example:** N = 3,4 and 6 have 'linear' complexity so they have just one primitive star point (star[1]) and hence one 'primitive' scale – which is the identity. In terms of cyclotomic fields, these all have $\varphi(N) = 2$ so $\mathbb{Q}(\zeta)$ must be of the form $\mathbb{Q}(\sqrt{d})$. It is easy to find d for N = 3 since $\zeta_3 = (i\sqrt{3} -1)/2$, Therefor $\mathbb{Q}(\zeta_3) = \mathbb{Q}(i\sqrt{3})$ and this works for N = 6 also since $\zeta_6 = (i\sqrt{3} + 1)/2$. For N = 4, $\mathbb{Q}(\zeta_4) = \mathbb{Q}(i)$.

**Example**: N = 6 shown here has primitive $s_1 = \tan(\pi/6) = \sqrt{3}/3$ and degenerate $s_2 = \tan(2\pi/6) = \sqrt{3}$, so in the Basis Lemma, $a_1 = 1/3$. Since star[1] = $\{-1/\sqrt{3},-1\}$, scale[2] = 1/3, so the identity scale will suffice. The minimal polynomial for $\tan[\pi/6] = T_6[x] = 3x^2 -1$ as shown on the right. This is not monic over $\mathbb{Z}$ so $\tan(\pi/6)$ (and GenScale[6] = $\tan^2(\pi/6)$ ) are not algebraic integers.

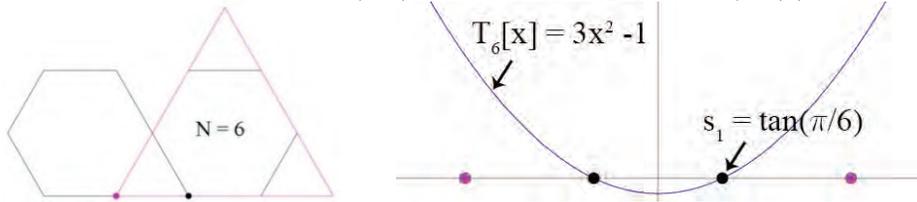

**Example**: N = 8: $s_1 = \tan(\pi/8)$ and $s_3 = \tan(3\pi/8)$ are primitive with values = $\sqrt{2} -1$ & $\sqrt{2} +1$. The degenerate $s_2$ is $\tan(2\pi/8) = 1$, so in the Basis Lemma, $(-1/2)s_1 + (1/2)s_3 = s_2$. As in the case of N = 12, cofunction symmetry yields scale[2] = $\sqrt{\text{scale}[3]}$ and the corresponding linear relationship is scale[2] = (1/2)(scale[1] – scale[3]) = (1/2)(1 – GenScale[8])

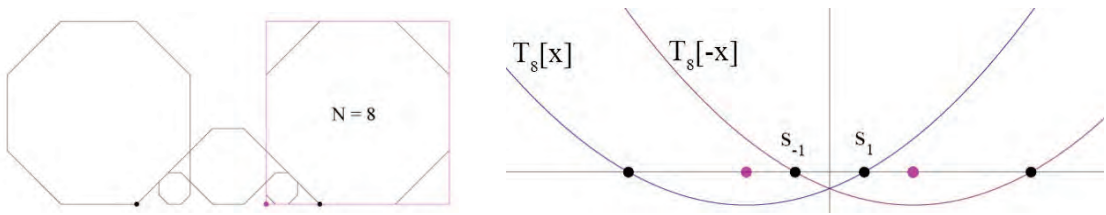

The full set of negative and positive primitive indices is $\{-3,-1,1,3\}$ and these are the prime residue classes mod 8, since $\{-3,-1,1,3\} = \{7,5,1,3\}$. The matching Galois group is isomorphic to $\mathbb{Z}^x_8 = \{1,3,5,7\}$. This is the first N with a non-cyclic Galois group and there are three subgroups of order 2. The three quadratic extensions are $\mathbb{Q}(i)$, $\mathbb{Q}(i\sqrt{2}) = \mathbb{Q}(\zeta_8)$ and $\mathbb{Q}(\sqrt{2}) = \mathbb{Q}(\zeta)^+$.

Like all 'quadratic' N-gons, N = 8 has only one non trivial scale - which is scale[3] = GenScale[8] = $(1-\sqrt{2})^2$. Renormalization yields a return-time ('temporal') scaling of 9, so the fractal dimension of the singularity set W is $-\text{Ln}[9]/\text{Ln}[\text{GenScale}[8]] \approx 1.246477$. See [S2].

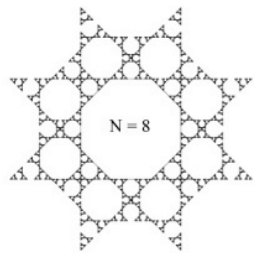

**Example**: For N = 9, the $s_k$ are {$\tan(\pi/9)$, $\tan(2\pi/9)$, $\tan(3\pi/9)$, $\tan(4\pi/9)$} and they are all primitive except for $s_3 = \tan(3\pi/9) = \sqrt{3}$ shown in magenta. In terms the Basis Lemma, there must be a rational solution to $a_1s_1 + a_2s_2 + a_3s_4 = \sqrt{3}$ and we noted earlier that such a solution does exist with the primitive $s_k$ of N = 18 so $\cot(\pi/9) - \cot(2\pi/9) + \cot(4\pi/9) = \sqrt{3}$. Therefore it is easy to derive the fact that $\tan(\pi/9) - \tan(2\pi/9) + \tan(4\pi/9) = 3\sqrt{3}$, so $(1/3)(s_1 - s_2 + s_4) = s_3$.

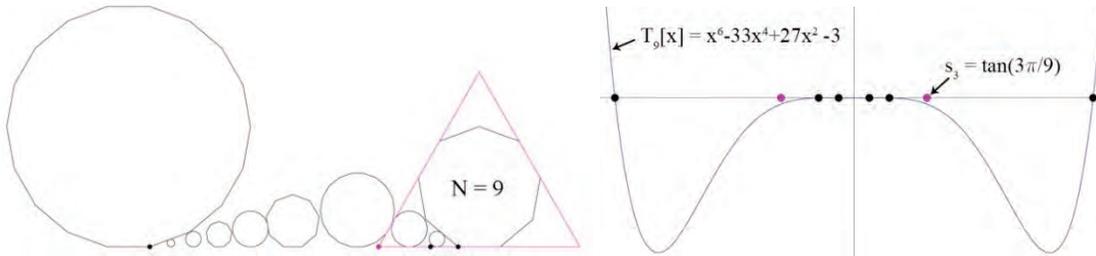

To get a linear relationship for the scales, the cotangent relationship above is perfect since $\sqrt{3} = 3/s_3$. Therefore scale[3] = (1/3)(scale[1] - scale[2] + scale[4]) = (1/3)(1 - scale[2] + GenScale[9]). (It is also true that scale[3]$^2$ = (1- 2scale[2])/3.) The important is issue is that N = 9 has just two non-trivial scales – just like N = 7. Dynamically speaking, N = 9 and N = 7 are quite similar even though N = 7 has no degenerate scales and no mutations.

Since 9 is the power of an odd prime, the Galois group $G_9$ must be cyclic of order $(3-1)(3^{2-1})$ so every divisor of 6 defines a subgroup. As always, the elements of $G_9$ are in 1-1 correspondence with the negative and positive indices of the primitive star points - namely {-4,-2,-1,1,2,4} = {5,7,8,1,2,4}. The generators are 2 and 5, and the proper subgroups are {1,4,7} and {1,8}.

**Example**: (N = 12) The $s_k$ are { 2- $\sqrt{3}$, 1/$\sqrt{3}$, 1, $\sqrt{3}$, 2 +$\sqrt{3}$}. These are shown on the right below with the primitive $s_1$ and $s_5$ in black. It should be clear that the degenerate $s_2$, $s_3$, $s_4$ (corresponding to N=6,4 and 3) are linear combinations of the conjugate pair $s_1$ and $s_5$.

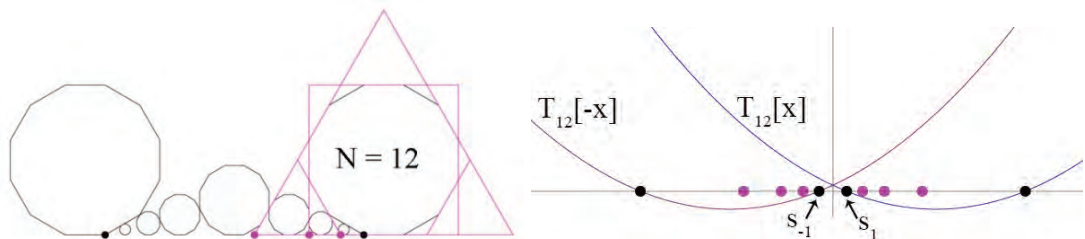

Since 4|N, the primitive $s_k$ are the combined positive roots of $T_{12}[x] = x^2-4x + 1$ in blue and $T_{12}[-x] = x^2 + 4x + 1$ in red. (In the case of N = 12, it would be sufficient to use just $T_{12}[x]$, but using $T_{12}[-x]$ will yield the full range of negative and positive roots.) With $T_{12}[x]$ and k = 1, the Root Lemma says that the only legal j values are 1 and 5, and both are less than N/2, so $T_{12}[x]$ generates the primitive $s_1$ and $s_5$ while $T_{12}[-x]$ generates the matching $s_{-1}$ and $s_{-5}$.

The full set of conjugate indices is {-5,-1,1,5}= {7,11,1,5} and these indices can serve as the prime residue classes mod 12 so the Galois group is isomorphic to $\mathbb{Z}^{\times}_{12} = \{1,5,7,11\}$. The matching automorphisms are {$\sigma_1, \sigma_5, \sigma_7, \sigma_{11}$} where $\sigma_k: \zeta \to \zeta^k$. This group is not cyclic. It is generated by the pair {$\sigma_5, \sigma_7$} and is in fact the Klein-4 group. Therefore N = 12 has three quadratic extensions - which correspond to the three degenerate star points above.

The quadratic extensions of N = 12 are $\mathbb{Q}(\sqrt{-3})$, $\mathbb{Q}(i)$ and $\mathbb{Q}(\sqrt{3})$. The maximal real subfield is $\mathbb{Q}^+ = \mathbb{Q}(\sqrt{3})$ and adjoining $i$ to $\mathbb{Q}^+$, yields $\mathbb{Q}(\zeta) = \mathbb{Q}(\sqrt{3}, i)$. Therefore N = 12 has the 'same' subfield structure as N = 8 with $\sqrt{3}$ in place of $\sqrt{2}$ - and this structure replicates itself with $\mathbb{Q}(\zeta_{24}) = \mathbb{Q}(\sqrt{2}, \sqrt{3}, i)$. Based on this fact, it would seem reasonable that the dynamics of N = 12 should be related to N = 8, and the dynamics of N = 24 should partially reflect both of these.

Of course $\sqrt{2}$ and $\sqrt{3}$ are the rotation parameters $2\cos(2\pi/8)$ and $2\cos(2\pi/12)$. As noted earlier, $\lambda_N = 2\cos(2\pi/N)$ always generates the maximal real subfield $\mathbb{Q}^+_N$. Below is a series of plots in the range $\sqrt{2}$ to $\sqrt{3}$ which show how the outer billiards web for N = 8 can be smoothly transformed into the web for N = 12 as the scaling in $\mathbb{Q}^+_N$ varies with rotation angle $\theta$.

These plots are based on the Digital Filter map described in Appendix B. This map mimics the outer-billiards map but it allows any $\theta$ in $(0, \pi/2]$ (and hence variable step-size). Here the angular decrements are $2\pi/120$, stating with N = 8 at $\theta = 2\pi(15/120)$, so the 4$^{th}$ plot will have $\theta = 2\pi(12/120)$, which gives N = 10 at $(\sqrt{5} + 1)/2$. In a sense this shows the full 'quadratic' range.

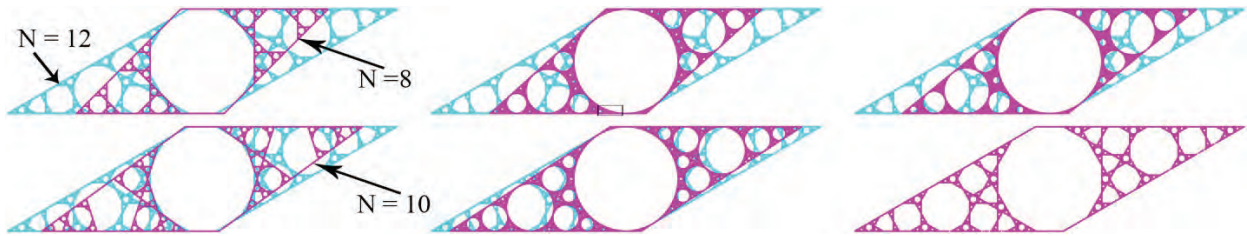

Some of these plots do not correspond to the dynamics of any known regular N-gon. For example the second plot has $\theta = 2\pi(14/120) = 2\pi(7/60)$. This would correspond to N = 60, but the 'step-size' would have to be 7 instead of 1 as in the traditional outer billiards map. For N = 60, the Df map allows 15 distinct step sizes, but step-15 would replicate N = 4. Step-14 is interesting because it reproduces the dynamics local to the central S[28] tile of N = 60. The rest are a mystery and have no obvious correlation with the step-1 outer billiards map of N = 60. (There have been no published studies of these multiple step-size cases, but possibly such a study could yield a better understanding of the step-1 case.)

The web plot on the left below is the same as final magenta image above. Since N = 12 has quadratic complexity there is only one non-trivial primitive scale, namely GenScale[12]. There appears to be uniform geometric and temporal scaling which persists for all 'generations' even though there are mutations in S[2] and S[3] (and the matching DS[2], DS[3]).

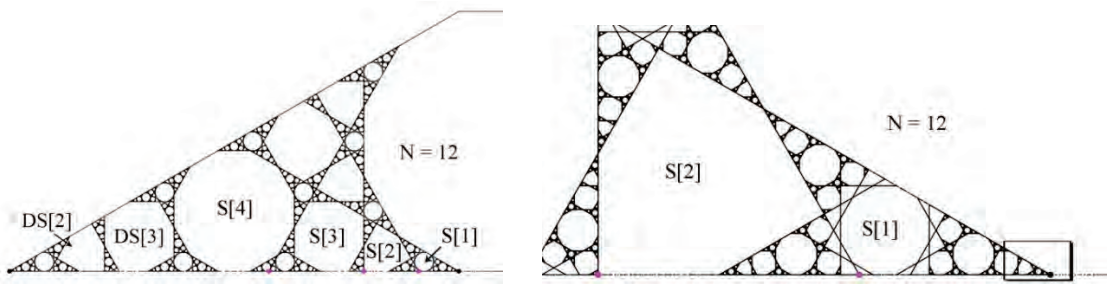

The mutations in S[2] and S[3] correspond to the 'degenerate' star points star[2] and star[3] with $(k,N) > 1$. But star[4] is also degenerate, and S[4] is not mutated. We have no definitive criteria for these mutations – but a necessary condition for the mutation of an S[k] appears to be the 'degeneracy' condition, $(k,N) > 1$, because this guarantees that the period of S[k] will be shorter than the canonical period N, and the local web development is related to this period. This degeneracy condition can also be applied to the star points of D and the DS[k] tiles, or the First Family of any regular tile.

The enlargement on the right above, hints at structure on all scales. When N is twice-even, there are no M-D pairs to serve as templates for future generations. Except for mutations, all the First Family tiles are just scaled copies of N. However the S[1] and S[2] tiles often foster extended generations with S[1] as the new 'matriarch'.

This makes sense because the S[1] tile is scaled by GenScale[N] relative to N and D. (When N is odd this scaling still exists but S[1] is a 2N-gon.) For N = 12 above, S[1] is a perfect regular dodecagon scaled by GenScale[12] ≈ 0.0717968 relative to N = 12. (This S[1] is unique in that it has its own invariant local 'web' and matching star polygons- so it is locally 'rational' and not quadratic. Since the matching star point is $s_2$ – which is star[1] of N = 6, this may be a remnant of the rationality of N = 6.)

If S[1] is the 'matriarch' of the next generation, the matching 'D' tile is the mutated S[2]. This poses no problem for future generations because mutations can propagate in the same fashion as unmutated tiles- so the next generation S[2] will be an exact copy of S[2] - scaled by GenScale[12].

In this chain converging to star[1], the tiles of generation[k] will be scaled by GenScale[12]$^k$ relative to generation[0] which is the First Family. In general the individual tiles in each generation may not be self-similar, but every generation will have an S[1] and S[2] in their canonical positions. For all 'quadratic' regular polygons the chains will eventually have perfect self-similarity. For N = 5, 8 and 10 this self-similarity begins with the First Family, but for N =12, it appears to begin with generation[1] shown on the right above.

The matching 'temporal' scaling of these tiles appears to be 27 in the limit. For example the S[1] tiles in the chain have periods 12, 420, 14148, 387252,.. with ratios which approach 27. The S[2] tiles also have periods which approach this same limit. This would yield a fractal dimension of :

$$-\text{Ln}[27]/\text{Ln}[\text{GenScale}[12]] = -3\text{Ln}[3]/2\text{Ln}[\text{Tan}[\pi/12]] \approx 1.2513$$

By comparison, the fractal dimensions of N = 8 is known to be

$$-\text{Ln}[9]/\text{Ln}[\text{GenScale}[8]] = -2\text{Ln}[3]/2\text{Ln}[\text{Tan}[\pi/8]] \approx 1.24648$$

This 3/2 similarity is probably a consequence of the fact noted above - that the cyclotomic fields are related by $\sqrt{3}$ and $\sqrt{2}$, so $\text{Tan}[\pi/12] = 2-\sqrt{3}$ and $\text{Tan}[\pi/8] = \sqrt{2} -1$. This might help to explain the 'unusual' period 27 temporal scaling for N = 12.

For odd N, GenScale[N] relates the $s_1$ terms of N and 2N, as shown below for N = 5, with temporal limit N+1

- Ln[6]/Ln[GenScale[5]] = -(Ln[2]+ Ln[3]) / (Ln[Tan[π/5] + Ln[Tan[π/10]])

= -(Ln[2]+ Ln[3]) / (Ln$\sqrt{5-2\sqrt{5}}$ + Ln$\sqrt{1-\frac{2}{\sqrt{5}}}$ ) ≈ 1.24114

It seems reasonable that these three 'quadratic' fractal dimensions are increasing with N, but for cubic complexity and higher, there are no obvious measures of dimension – since there are multiple primitive scales. However if each of these N-gons has a well-defined spectrum of dimensions, there will be a maximal Hausdorff dimension, and we conjecture that these maximal Hausdorff dimensions will approach 2 as N increases.

**Note**: At this time there is no theory that attempts to relate the 'in-situ' dynamics of a polygon with the 'in-vitro' dynamics – with the polygon as the generator. The maximal D tile for a regular N-gon is an important exception. It has the same local web whether it is at the origin or simply a First Family member, S[<N/2>]. The Df Theorem in [H2] provides some isolated examples – such as the case of N = 60 mentioned earlier. But these involve iterations with multiple steps and this process is poorly understood.

Most of the natural mutations which occur in the First Family are equilateral and not equiangular as shown on the left below for S[3] of N = 12. They preserve half the dihedral symmetry of the regular case, but the in-vitro dynamics of these 'woven' polygons shows imperfect rings of D tiles and little hope of bounded dynamics. By rotating the internal factors of a regular N-gon as shown on the right, the angles are preserved and it still has half the dihedral symmetry. These 'Riffle' mutations are better behaved. See [H2],[Ho] and [S3].

| The S[3] tile of N = 12 is an equilateral octagon as shown below. This is the typical form of a 'mutated' tile. We call it a 'woven' octagon because it has two embedded squares with different radii. Its dihedral symmetry group is $D_4$ instead of $D_8$ for the regular octagon, | The 'Riffle' octagon shown below is formed by rotating one of the embedded squares. It is equiangular and also has dihedral symmetry group $D_4$. We call these 'semi-regular' polygons and they are part of our quest to find classes of non-regular polygons with bounded dynamics. It has been conjectured that this class is measure zero. |
|---|---|
| 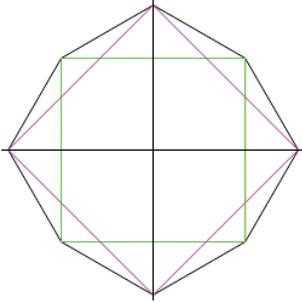 | 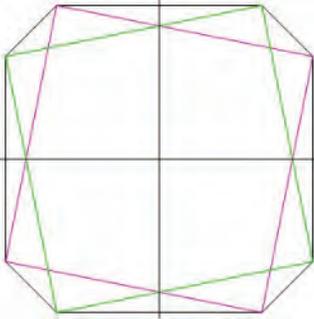 |

# Appendix A: Webs and Generalized Star Polygons

Since the star polygon $\{N,\langle N/2\rangle\}$ for a regular polygon consists of 'extended' edges of N, it is part of the 'singularity' set W for the outer billiards map $\tau$. We refer to W as the 'web'. The web at level-k, is defined to be those points where the outer billiards map $\tau^k$ (or its inverse) is not defined. For any piecewise affine map, the limiting web is very important because the complement defines regions (tiles) where the dynamics are always well-defined. Each S[k] of the First Family for N defines a step-k orbit of $\tau$, so they all have period N or N/k when (k,N) >1. This means that these tiles arise early in the web generation process and hence they can serve as templates for the web evolution.

It is not difficult to show that the region defined by the star polygon $\{N,\langle N/2\rangle\}$ becomes invariant in fewer than $\langle N/2\rangle$ iterations of the web. This means that this 'inner star' region can be used as a template to study the local web evolution (which is conjugate to the global evolution when N is regular).

Below are the first few iterations of the local web W for N = 14 – showing levels 1,2,3,4,5 and 10. Since the star polygon $\{14,6\}$ is the level-1 (local) web for N = 14, the subsequent webs are all generated from $\{14,6\}$. Therefore they can be regarded as 'generalized' star polygons of N = 14. In general, this web structure appears to be fractal or multi fractal and it is a curious fact that this fractal structure arises from a mapping which is discontinuous but definitely not 'chaotic' - yet it closely resembles the Poincare cross-sections which arise from chaotic mappings with positive Lyapunov exponents.

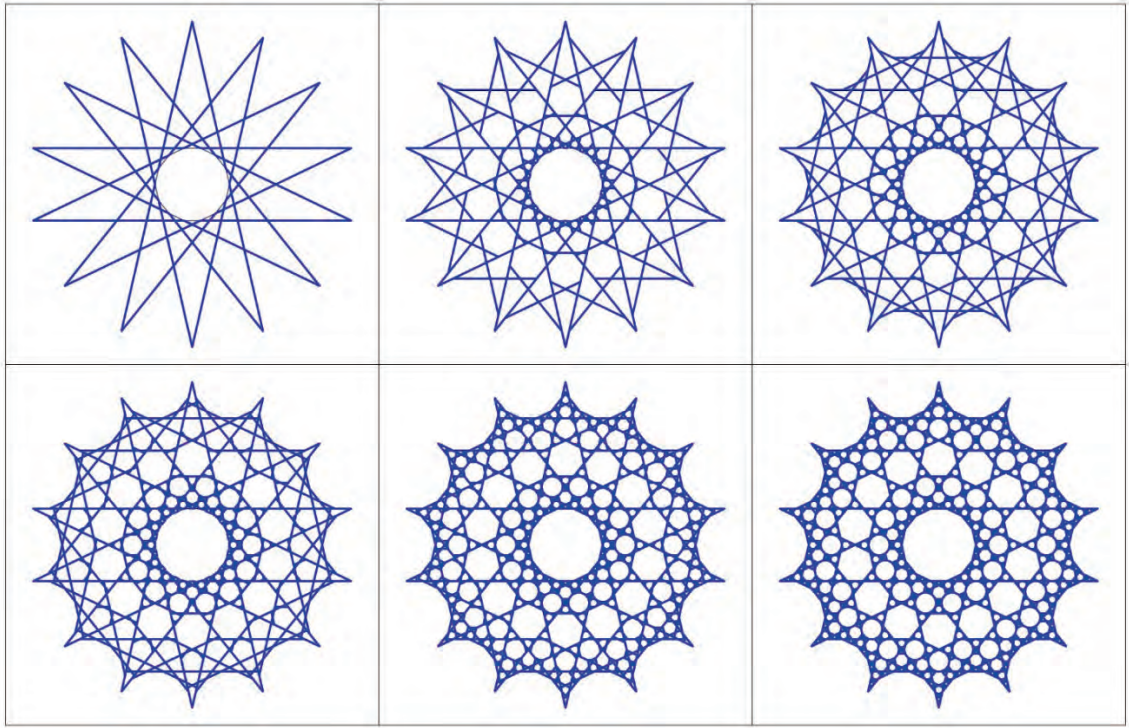

The symmetry of the N-even family implies that this template can be further reduced to half of the magenta rhombus shown below. This is what we call the First Family template. The Digital Filter map (Df) mimics the outer billiards map on this 'torus' and provides a manifold increase in 'space' and 'time' efficiency. The image on the right is a level 1000 Df web which takes just a few line of code and generates this web in a few seconds.

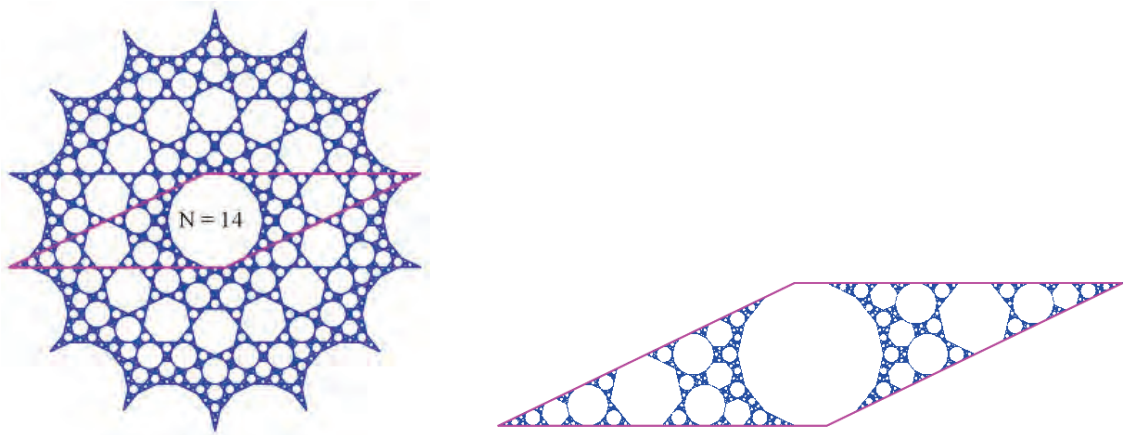

Note that there are 7 of these overlapping rhombi, and the overlaps cancel so that each rhombus contains exactly 1/7 of the tiles. This is only true for N even. Below is the embedding of this rhombus in the star polygon for N = 14.

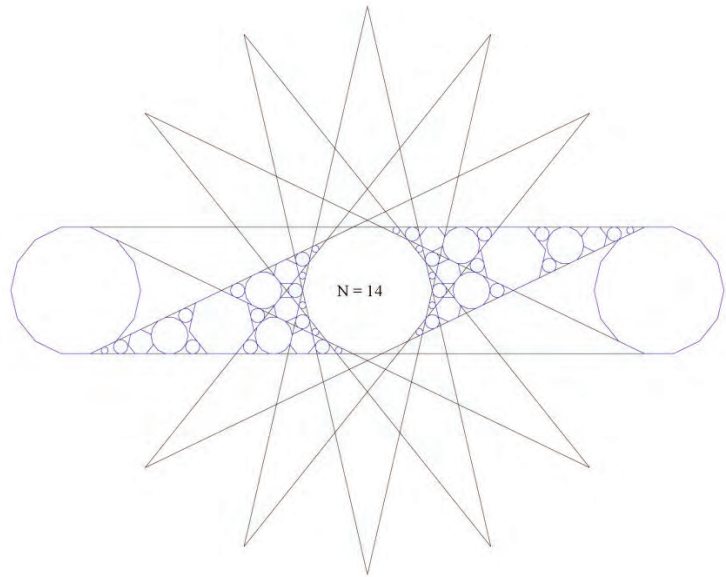

The Digital Filter map and the outer billiards map are examples of piecewise rational rotations. These maps will be discussed below. See also [D], [A], [H1], and Appendix F of [H2].

## Appendix B: Scaling for Rational Rotations

The Digital Filter map is an example of an *affine piecewise rotation*. These are two-dimensional analogs of interval exchange maps- where rotations typically play the part of exchanges. The linear form for an affine piecewise rotation by $\theta$ can be written as:

$$T\begin{bmatrix} x \\ y \end{bmatrix} = \begin{bmatrix} 0 & -1 \\ 1 & 2\cos\theta \end{bmatrix}\begin{bmatrix} x \\ y \end{bmatrix} \quad \text{or in complex form } T[z] = \rho z \text{ where } |\rho| = 1$$

Typically the phase space X is partitioned unto a finite number of mutually disjoint 'atoms' $A_i$ and $T_i$ acts on atom $A_i$ as a rotation possibly followed by a translation. The combined map T: X →X is defined as $T(z) = T_i(z)$ iff $z \in A_i$. When the map T is bijective, the image of a partition is a partition. When the rotation angle $\theta$ is a rational multiple of $\pi$, T is called a *rational rotation* or a polygonal rotation.

**Example**: The Digital Filter map has linear form and matching Jordan normal form given by

$$A = \begin{bmatrix} 0 & 1 \\ -1 & 2\cos\theta \end{bmatrix} \sim \begin{bmatrix} \cos\theta & \sin\theta \\ -\sin\theta & \cos\theta \end{bmatrix}$$

The space in which Df operates is $[-1,1)^2$ so it is a map on a 2-Torus: $T^2 \to T^2$. Like all affine piecewise rotations, $\text{Det}[A] = 1$, so it preserves area and is called a symplectic map. Symplectic maps on tori have been an area of interest to mathematicians and physicists since Henri Poincare (1854-1912) realized their value in the analysis of conservative (Hamiltonian) systems. In physics these are sometimes called 'kicked' Hamiltonians.

The eigenvalues of A are complex of the form $\lambda$ and $1/\lambda$ with unit absolute value so $\lambda = e^{2\pi i \theta}$ and A represents a rotation. But this is an 'elliptical' rotation which can be conjugated to a pure rotation. When studying the dynamics of maps based on matrices such as A, if the trace $2\cos\theta$ is the solution to a polynomial equation of low degree, there are computational advantages (exact arithmetic) to leaving A in its original form instead of its conjugate form.

To reproduce the (local) outer billiards dynamics of a regular N-gon set $\theta = 2\pi/N$, and then $\lambda_N = 2\cos\theta$ is the trace of A and also the generator of the maximal real subfield of $\mathbb{Q}(\zeta_N)$ and its minimal polynomial is always monic of degree $\varphi(N)/2$. This simplifies algebraic analysis, but the 'underflow' and 'overflow' conditions of the toral map necessitates the use of a sawtooth function such as $f(z) = \text{Mod}[z+1,2]-1$. (In 1997 Peter Ashwin [A] showed that the digital filter map is equivalent to a sawtooth version of the Standard Map – which is the classical model for Hamiltonian 'chaos'.)

**Definition**: The Digital Filter map Df: $[-1,1)^2 \to [-1,1)^2$ is defined as

Df[{x,y}]:={y, $f(-x + 2y\cos\theta)$} where $f(z) = \text{Mod}[z+1,2]-1$

For any given value of the trace $2\cos\theta$, Df is an affine piecewise rotation . The rotations can be irrational or rational but by symmetry the rational rotations need N to be even.

**Example**: $N = 14$, $\theta = 2\pi/14$. The three 'atoms' A, B and C are shown below. Region A is the 'overflow' region and C the 'underflow' – with B the 'linear' region. The Df map applies a (clockwise) elliptical rotation of $\theta = 2\pi/14$ to each region and the sawtooth nonlinearity $f$ provides the corresponding translation – which is vertical by -2, 0 and + 2 respectively for A, B and C as shown on the right below.

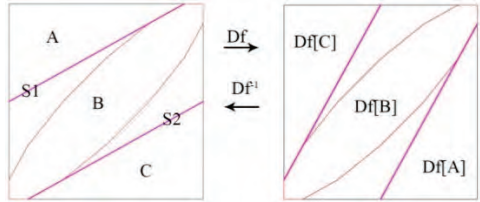

The 'speratices' S1 and S2 are extended edges of N, so the atoms mimic the star structure of $N = 14$. For the outer billiards map the web is formed by iterating these extended edges under $\tau$. Here it is only necessary to itereate S1 and S2 under Df. The resulting level 1000 web is shown below in magenta along with its rectified version in blue. This is a prefect copy of the local web for $N = 14$ – and hence also a perfect copy of the local $N = 7$ web.

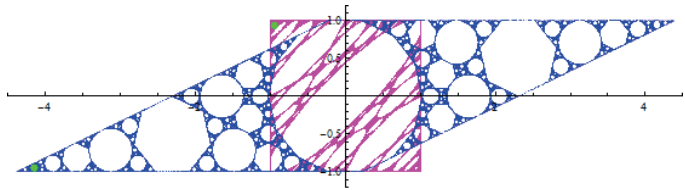

There are computational advantages to using triangular regions for the atoms and it may be possible to implement the outer billiards map with a three triangle map such as the 'dart' shown here. Mappings such as this were studied by Adler, Kitchens and Tresser [AKT (2001)]

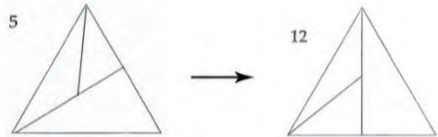

A recent study by X. Bressaud and G. Poggiaspalla [BP] used computer analysis to categorize the possible bijective polygonal piecewise isometries with a phase space consisting of two or three triangles. For two triangles there is only one case which yields non-trivial dynamics. This is what they call the 'tower' case.

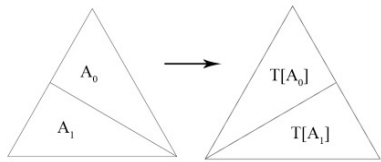

There are only two ways to implement a reflection such as this using orientation preserving rotations and A. Goetz [Go1] found them using rotations by $\pi/5$ and $\pi/7$, but the most natural phase space for $\pi/7$ is the three triangle case which we will look at below.

**Example**: The two triangle $\pi/5$ case of [Go1] is based on $N = 10$ as shown below. It can be implemented in complex form using $a = \zeta_{10} = \cos(\pi/5) + i\sin(\pi/5)$. The triangles $A_0$ and $A_1$ have vertices $\{0, a^2 + a^4 + a^6, -1\}$ and $\{0, -1, a^6\}$. The piecewise rotations and corresponding translations are $T_0(z) := a^4 z + a^2 + a^4 + a^6$ and $T_1(z) = a^6 z + a^6$ (so the rotations are $\pm (\pi - \pi/5)$). The combined map is $T(z) = T_0(z)$ if $\text{Im}[z] > 0$, otherwise $T_1(z)$. The first iteration is shown below along with the residual set obtained by iterating a non-periodic point.

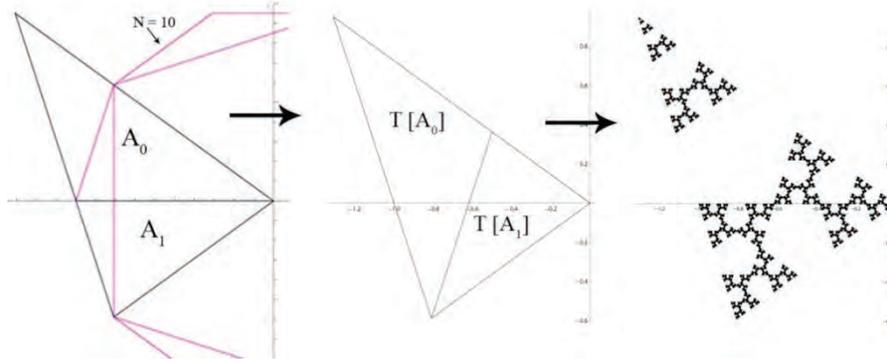

The geometric scale factor is $k = 1/\lambda_{10}$ where $\lambda_{10} = 2\cos(2\pi/10)$ is the Golden Ratio $\gamma$. The limiting singularity set consists solely of pentagons and the 'temporal' scaling is 2, so the Hausdorff dimension is $\text{Log}[2]/\text{Log}[k] \approx 1.4404$. By comparison, $N = 5$ (and $N = 10$) have a scale factor of GenScale[5] – which is $k^3$. The temporal scaling for $N = 5$ is 6 so the Hausdorff dimension can be written as $-(\text{Log}[2] + \text{Log}[3])/3\text{Log}[k] \approx 1.241$.

With a phase space based on three triangles Bressaud and Poggiaspalla found 810 cases involving rational rotations– including 258 'cubic' cases based on rotations by $\pi/7$. One of these three-atom $\pi/7$ solutions is an extension of the example above. It has been studied by Goetz, Poggiaspalla, Kahng, Lowenstein, Vivaldi and Kapustov. This example is presented below.

**Example:** Set $a = \zeta_{14} = \cos(\pi/7) + i\sin(\pi/7)$ and define triangles $A_0$ and $A_1$ and $A_2$ with vertices $\{0, a^5-1,, -1\}$, $\{0, -1, -a^3\}$ & $\{0, -a^3, -a^3+a^2\}$. The corresponding transformations are $T_0(z) = za^6 + a^5 -1$, $T_1(z) = -az - a^4 + a^5 - a^6 -1$ and $T_2(z) = za^6 - a^3$. Then $T[z] = T_i(z)$ iff $z \in A_i$. (Note that the rotations are $\pm (\pi - \pi/7)$, so swapping $\pi/7$ and $\pi/5$ will make this example the same as the example above when restricted to $A_0$ and $A_1$.)

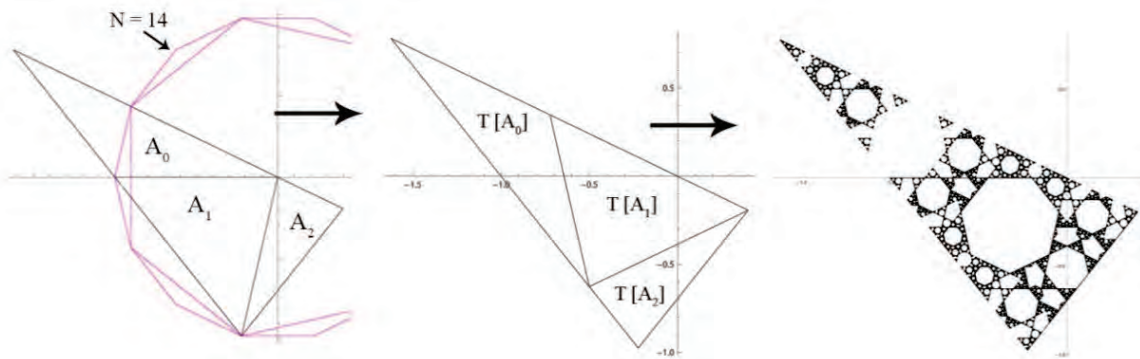

The large heptagons $H_0$ and $H_1$ have centers which are the fixed points of $A_0$ and $A_1$. The $\pi/7$ 'tower' shown here appears to be self-similar beginning with $H_1[1]$. Goetz and Poggiasapalla [GP] used symbolic complex analysis to aid in the analysis of these sequences of $H_0$'s and $H_1$'s.

They determined that the sequences share a geometric scale factor of $\upsilon = 4\sin^2(\pi/14)$ (which is $\lambda_{14} \cdot$ GenScale[7] where $\lambda_{14} = 2\cos(\pi/7)$ is the trace of their rotation matrix). The temporal scale factors of these sequences are $(2^{2n+1}+1)/3$ and $(4^{n+1}-1)/3$, so the temporal scaling approaches 4 in either case - one from the top and the other from the bottom. This is similar to the First Family scaling of $N = 7$ and $N = 14$, but there are significant differences in the detailed scaling.

In LKV],Lowenstein et al modified this map so that it had a 'recursive tiling' property that enabled them to determine the scaling sequences and hence the spectrum of dimensions with maximal Hausdorff dimension $\approx 1.652$. More on this below.

These examples are all affine piecewise rotation so they share the same linear model and all are based on rational rotations with $\lambda_N = 2\cos(2\pi/N)$ as the trace of the rotational matrix. In practice the most common scenario is for N to be odd with rotations based on $\lambda_{2N} = 2\cos(2\pi/2N)$ and scaling terms in the ring of integers $\mathbb{Z}[\lambda_{2N}]$.

We have shown that GenScale[N] or GenScale[N/2] are natural choices of generators for $\mathbb{Q}^+_N$ - alongside $\lambda_N = 2\cos(2\pi/N)$ and below we show how these generators are related:

**Lemma 11**: When N is even GenScale[N] $= \tan^2(\pi/N) = \dfrac{4 - 4Cos^2(2\pi/2N)}{4Cos^2(2\pi/2N)} = \dfrac{4 - \lambda_{2N}^2}{\lambda_{2N}^2}$

If N is odd GenScale[N] $= \tan(\pi/N) \cdot \tan(\pi/2N) =$

$$\frac{\sin(\pi/N)}{\cos(\pi/N)}\left(\frac{1-\cos(\pi/N)}{\sin(\pi/N)}\right) = \frac{1-\cos(\pi/N)}{\cos(\pi/N)} = \frac{2 - \lambda_{2N}}{\lambda_{2N}}$$

Therefore for N even, $\lambda_{2N} = \dfrac{2}{\sqrt{\text{GenScale}[N] + 1}}$ and for N odd $\lambda_{2N} = \dfrac{2}{\text{GenScale}[N] + 1}$ □

For a cyclotomic field $\mathbb{Q}(\zeta_N)$, the matching ring of integers is $\mathbb{Z}[\zeta]$ where $\zeta$ is a unit generator for $\mathbb{Q}(\zeta_N)$ so it is no surprise that the ring of integers in $\mathbb{Q}(\zeta_N)^+$ is $\mathbb{Z}[\lambda_N]$ and this will remain true for any integral generator of $\mathbb{Q}(\zeta_N)^+$.

Dirichlet's Unit Theorem, says that the group of units in a number field K is finitely generated with rank given by $r = \{r_1 + r_2 - 1\}$ where $r_1$ and $r_2$ are the number of real and pairs of complex embeddings. For a cyclotomic field $\mathbb{Q}(\zeta_N)$ there are no real roots or embeddings and there are $\varphi(N)/2$ complex pairs of embeddings, so r is always $\varphi(N)/2 - 1$. The maximal real subfield $\mathbb{Q}(\zeta_N)^+$ has the same rank with real roots instead of complex.

Therefore the unit groups $U_N$ and $U_N^+$ have the same rank. Since they are finitely generated abelian groups – they will have the form $T \oplus \mathbb{Z}^r$ where T is the torsion subgroup. When K is $\mathbb{Q}(\zeta_N)$, T will be the group of roots of unity because they have finite order, so every element of $U_N$ will have the form $\zeta e_1^{j_1} e_2^{j_2} \ldots e_r^{j_r}$ with the $e_i$ units in $\mathbb{Z}[\zeta]$ and the $j_i$ integers. The $e_i$ are called the 'fundamental units' and they form a (multiplicative) basis for $U_N$. All of this remains true for $\mathbb{Q}(\zeta_N)^+$ but the torsion subgroup is just $\{\pm 1\}$. Dirichlet's result does not address the issue of how $U_N$ and $U_N^+$ are related and that will be discussed below.

**Example**: N = 14: Mathematica will generate fundamental units for $\mathbb{Q}(\zeta)$ or $\mathbb{Q}(\zeta)^+$ - although the results are not unique. **NumberFieldSignature** will return $\{r_1, r_2\} = \{3,0\}$ and $\{0,3\}$ respectively for $\lambda_{14}$ and $\zeta_{14}$, so r = 2 in both cases.

**NumberFieldFundamentalUnits[Exp[2IPi/14]]** yields $e_1 = \zeta + \zeta^5$ and $e_2 = 1 - \zeta + \zeta^2 - \zeta^3 - \zeta^5$
And we will see below how these match up with units in $U_N^+$.
**NumberFieldFundamentalUnits[2Cos[2 Pi/14]]** yields $1/\lambda_{14} - 1$ and $2 - \lambda_{14}^2$ with product $1/\lambda_{14}$

In studies of rotations by $\pi/7$ from [LKV], the authors chose similar fundamental units of $1/\lambda_{14}$ and $1 - \lambda_{14}$ which are found in Appendix B of [C]. These can be written in terms of GenScale[7] as: $\frac{\text{GenScale}[7]+1}{2}$ and $\frac{\text{GenScale}[7]-1}{\text{GenScale}[7]+1}$ with product $\frac{\text{GenScale}[7]-1}{2}$

We have shown that the primitive scales for any regular N-gon form a unit basis for $\mathbb{Q}(\zeta)^+$. For N = 14, the primitive scales are scale[1] = 1, scale[3] and scale[5] = GenScale[7]. It is easy to write $\lambda_{14}$ as a linear combination of any of these scales since $\lambda_{14}$ is itself a unit. We noted earlier that $\lambda_{14} = 4\sin^2(\pi/14)/\text{GenScale}[7]$ and $\lambda_{14}$ is also scale[3]·$\sin(3\pi/14)/2\sin^2(\pi/14)$. Clearly any primitive scale can serve as a 'surrogate' unit generator for $\mathbb{Q}(\zeta)^+$ but this does not imply that the primitive scales are 'fundamental' units. However, based on the examples above with N = 14, there may be efficient algorithms for generating the fundamental units using the primitive scales.

If Mathematica is asked to find the fundamental units using scale[3] or scale[5] – the results will be identical: $\{\lambda_{14}^2 - 2 \text{ and } \lambda_{14}\} = \{(5-\text{GenScale}[7])^2/4 \text{ and } (2/(\text{GenScale}[7]) + 1)\}$ - and these are essentially the same as the results shown above using $\lambda_{14}$ as generator - or the fundamental units used in [LKV].

Theorem 4.12 in [W] shows a close relationship between the units in any CM field such as $\mathbb{Q}(\zeta_N)$ and the matching maximal real subfield $\mathbb{Q}(\zeta_N)^+$ - namely if $U_N$ is the unit group in $\mathbb{Q}(\zeta_N)$ and $U_N^+$ is the unit group in $\mathbb{Q}(\zeta_N)^+$, then $[U_N: W \cdot U_N^+] = 1$ or 2 where W is the group of roots of unity in $\mathbb{Q}(\zeta_N)$. This follows because any CM field is a totally imaginary quadratic extension of a real subfield, so complex conjugation is well defined and fixes the real subfield. Therefore for any unit $e$ in $\mathbb{Q}(\zeta_N)$, $|\bar{e}|$ is a unit in $\mathbb{Q}(\zeta_N)^+$ and $e/|\bar{e}|$ is a unit with norm 1 so $\pm$ a root of unity.

**Example**: Note that $e = 1 - \zeta_N$ will be a unit in $\mathbb{Q}(\zeta_N)$ whenever N is not a prime power, so for N = 14, $e/|\bar{e}|$ is of the form $\pm \zeta^k$ – and here the sign is negative and k = 3. Clearly the index is 2 and it is easy to see that the positive sign for k = 3 will occur when e is the 'fundamental' unit $\zeta + \zeta^5$.

Since the units of $\mathbb{Q}(\zeta_N)$ and $\mathbb{Q}(\zeta_N)^+$ only differ by roots of unity, the basic units are the 'same', and hence the ratio of the regulators is $2^r$ for N a prime power (two complex roots for each real root) – but if N is not a prime power, the index is 2 and the ratio is reduced by 2. These issues are important in the computation of class numbers. Kummer's criteria is based on a comparison of the data and structure of $\mathbb{Q}(\zeta_N)$ and $\mathbb{Q}(\zeta_N)^+$ which involves 'dividing out' their Dedekind zeta functions so that the ratio of the class numbers can be reduced to their Bernoulli numbers.

**Appendix C – The Tile-Scaling Conjecture for Regular Polygons**

**Definition**: For a regular N-gon , we define a *'canonical regular tile'* to be any regular tile which can be scaled (relative to N) by an element of the scaling field $S_N = \mathbb{Q}_N^+$.

**Example**: Every regular N-gon has a matching maximal D tile and this tile is always canonical because when is even, hD/hN = 1 and when N is odd, hD/hN = ScaleSwap[N,2N] which is in $S_N$ because $\mathbb{Q}(\zeta_N) = \mathbb{Q}(\zeta_{2N})$. The remaining First Family members are canonical because their heights are given by hS[k] = hD·GenScale[N]/scale[k] which is in $S_N$.

For 'quadratic' regular N-gons, like N = 5 and N = 8, all tiles are 'canonical regular tiles' because they are either First Family members or First Family members scaled by GenScale$^k$. This is no longer true for N-gons with cubic complexity of higher.

**Definition**: For a regular N-gon , we define a *'canonical non-regular tile'* to any non-regular tile whose sides can be scaled (relative to N) by an element of the scaling field $S_N$

**Example:** All 'mutated' First Family tiles are equilateral non-regular tiles which are canonical because their edges are the difference of two star points. We will demonstrate this with N = 9.

**The Tile-Scaling Conjecture for Regular Polygons**: All tiles which arise in the outer-billiards dynamics of a regular N-gon are canonical.

**Example**: N = 18 (and N = 9) have cubic complexity along with N= 7 and N = 14. Therefore there are 3 primitive scales – but one is the trivial scale[1] = 1, so the dynamics of N = 9 and N = 7 are based on just the two non-trivial primitive scales – namely scale[2] and scale[4] for N = 9 and scale[2] and scale[3] for N = 7. This seems to imply that the dynamics in both cases are a mixture of 'quadratic –type' self-similar dynamics (when the overlap of scales is suppressed) and 'cubic-type' complex dynamics when the two scales interact. This dichotomy can be seen in the second generation below for N = 18 – with D[1] playing the role as the new D..

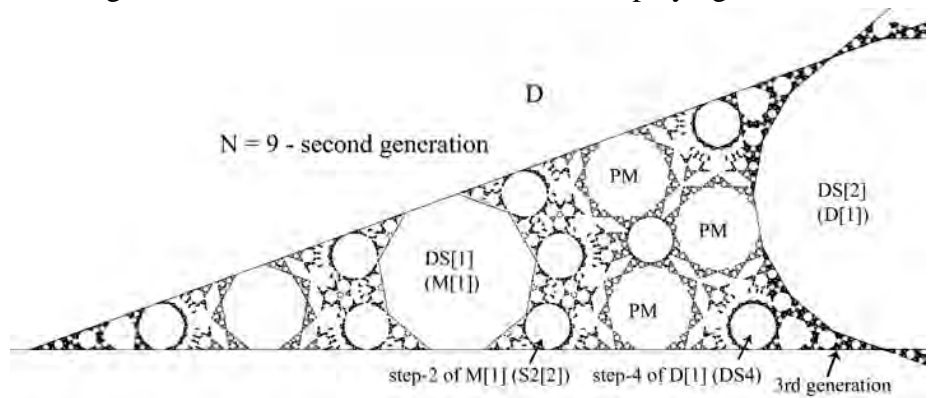

The dynamics on the edges of D[1] are very 'canonical' but this breaks down in the vicinity of DS4 so the PM tiles (and the elongated octagons) have no obvious relationship with the First Family of D[1]. But we will give evidence to support our conjecture that all the tiles shown here are 'canonical'. In particular we will show that the mutated DS3, the PM tile and the elongated octagons are all 'canonical'. (For convenience we will work inside N = 18 (with height h), so the star points are of the form star[k] = ±h·{$s_k$,1} where $s_k$= Tan[kπ/18] for {k,1,18}]).

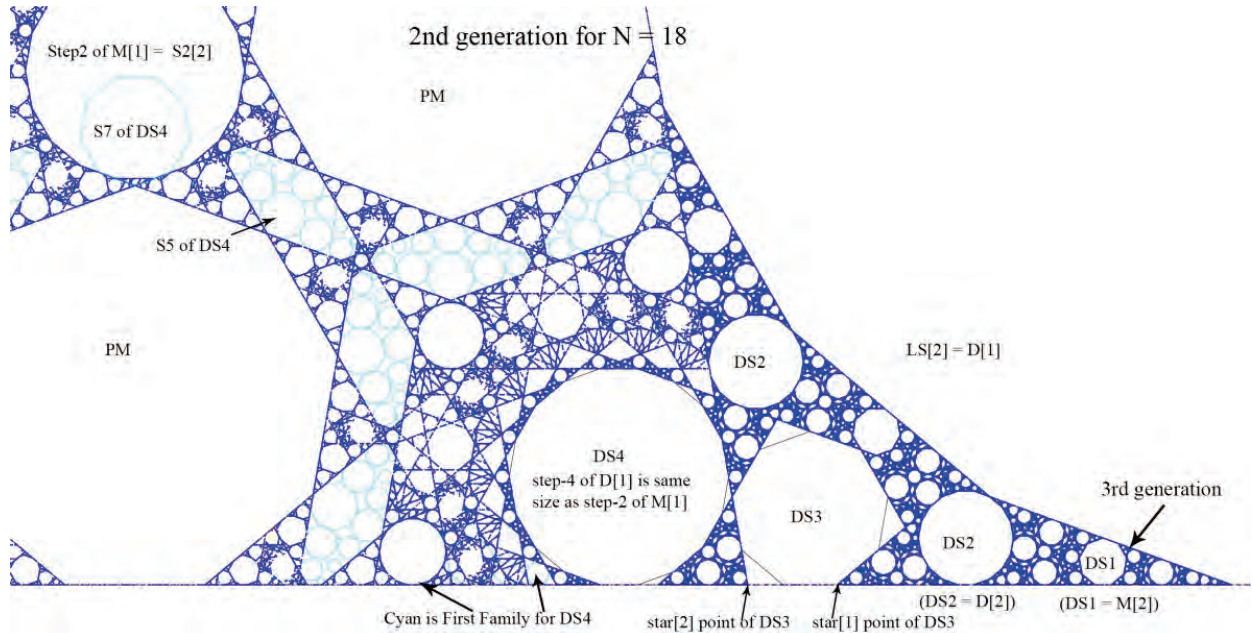

(i) To show that the 'mutated' DS3 shown above is 'canonical'.

Note: The mutations shown here in DS3 and DS4 are typical 'extended edge' mutations which result from incomplete local webs. These webs are the 'singularity sets' formed from mapping the extended edges of N = 18 under $\tau^{-1}$. In the First Family of any regular N-gon, the S[k] tiles will normally have period N, but when GCD[k,N] >1, this period can be shortened and this can result in a non-regular tile with some edges longer than the regular case. These mutations can occur in any generation. The DS3 shown above is a 2$^{nd}$ generation scaled copy of the mutated S[3] tile of the First Family. This mutation occurred because the points in S[3] have period 6 instead of the canonical period 18. (The DS4 tile shown here is also based on a scaled copy of S[4], but S[4] itself is not mutated - even though it has period 9. We do not have sufficient conditions for these mutations to occur – but the resulting tile will always be equilateral with side length scaled by an element of the scaling field $S_N$ , so it will be a 'canonical' non-regular tile.

The side of the mutated DS3 is the distance between star[2] and star[1] of the underlying regular 9-gon. (it is irrelevant which side is used for star points – since they are symmetric.) So

sDS3 = star[1][[1]] + star[2][[1]] (of DS3) = hDS3·($s_1$ + $s_2$) = hDS3·(Tan[π/9] + Tan[2π/9])

Finding hDS3 is easy because back in the First Family for N = 18, any S[k] tile has hS[k] = h·GenScale[18]/scale[k], so hS[3] = h·GenScale[18]/scale[3] and the scale to go from N = 18 to S[2] (or LS[2]) is GenScale/scale[2], so the side of DS3 is

sDS3 = h·(GenScale[18]$^2$/(scale[3]·scale[2]))·(Tan[Pi/9] + Tan[2Pi/9]) ≈ 0.00786020·h

By symmetry all mutated N-gons are equilateral and retain half the dihedral symmetry of the un-mutated tile so sDS3 represents an arbitrary side of DS33 and this formula gives this side to arbitrary accuracy – as long as h is known.

However our goal here is to show that sDS3 can be scaled relative to N = 18 (or N = 9) by an element in the scaling field $\mathbb{Q}^+$. Since the un-mutated DS3 is a regular 9-gon, it would be appropriate to scale it relative to the side of N = 9 – which is sN9 = h*2*Tan[π/9]. This is just a 'gender change' from N = 18, so the height will be the same as N = 18, so the h's will cancel.

DS3scale = sDS3/sN9 = (1/2)·Cot[π/9]·(GenScale[18]$^2$/(scale[3]·scale[2]))·(Tan[Pi/9]+Tan[2Pi/9])

It would be easy to expand the right-side to show that it is in the scaling field $S_{18}$ (or $S_9$), but the lazy approach is to let Mathematica do the work. Our convention for twice-odds is to use $S_{N/2}$ instead of $S_N$, and we can do this inside N = 18 because GenScale[9] = scale[7] of N = 18.

**AlgebraicNumberPolynomial[ToNumberField[DS3scale,scale[7]],x]** gives

$$-\frac{7}{12}+\frac{55x}{6}+\frac{17x^2}{12} \quad \text{so DS3scale} = -\frac{7}{12}+\frac{55}{6}Tan[\frac{\pi}{18}]Tan[\frac{\pi}{9}]+\frac{17}{12}Tan[\frac{\pi}{18}]^2Tan[\frac{\pi}{9}]^2 \approx .01079786$$

Of course this new scale is independent of the dimensions of N = 9. It matches our calculations above since sDS3 = DS3Scale*2hTan[π/9] = 0.00786020063·h

Below is a vector plot of the mutated DS3:

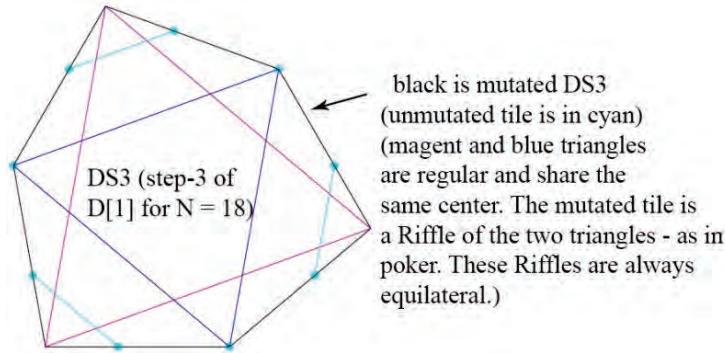

DS3 (step-3 of D[1] for N = 18)

black is mutated DS3 (unmutated tile is in cyan) (magent and blue triangles are regular and share the same center. The mutated tile is a Riffle of the two triangles - as in poker. These Riffles are always equilateral.)

This is what we call a 'woven' or Riffle polygon because it consists of two regular triangles with a common center and different radii. For DS4, the resulting Riffle consists of two regular hexagons. Below is the mutation of S[2] of N = 20 - which consists of two regular pentagons. (From the standpoint of dynamics, N = 10 and N = 20 are virtually unrelated.)

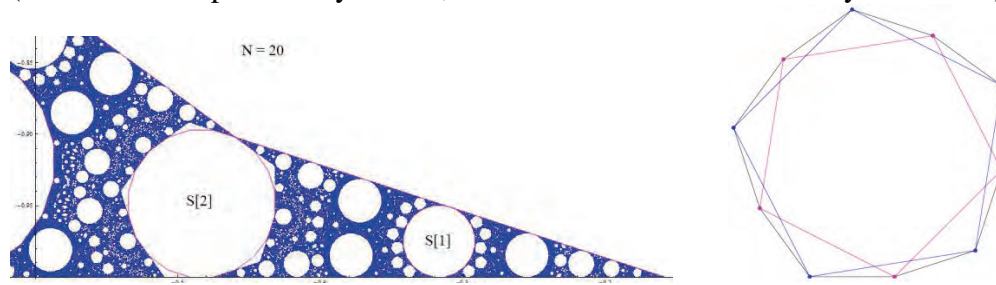

(ii) The scaling for the PM tiles of N = 18 is slightly more difficult – because the two star points that define sPM are in different tiles (but within the same 'family'). Here the underlying First Family is relative to the step-4 tile of D[1] – which we call DS4. It can be observed in the vector plot above that sPM is the distance between $s_1$ vertices of S7 and S5 of DS4. We will do this calculation in the First Family of N = 18 – and then scale the result to return to the 2$^{nd}$ generation.

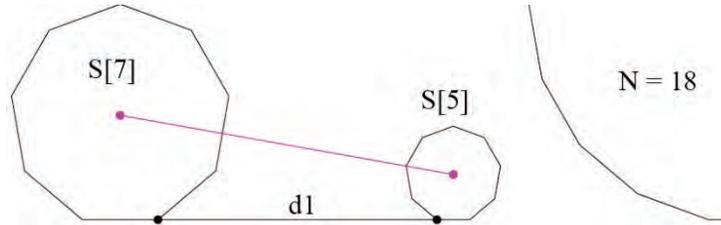

The distance d1 shown here is the side of PM scaled up to the First Family for N = 18. We will find this distance using the fact that the horizontal distance between centers of S[7] and S[5] is $h[Cot[\frac{2\pi}{9}] - Cot[\frac{\pi}{9}]]$. Therefore d1 = $h[Cot[\frac{2\pi}{9}] - Cot[\pi/9]]$ - (sum of $s_1$ points of S[7] and

S[5])] = $h[Cot[\frac{2\pi}{9}] - Cot[\pi/9]]$ - [Tan[π/9]·h·GenScale[18]($\frac{1}{scale[5]} + \frac{1}{scale[7]}$)]

So d1 = $h[Cot[\frac{\pi}{9}] - Cot[\frac{2\pi}{9}] - Tan[\frac{\pi}{9}](Cot[\frac{\pi}{9}]Tan[\frac{\pi}{18}] + Cot[\frac{2\pi}{9}]Tan[\frac{\pi}{18}])] \approx 1.30291 \cdot h$

Back in the 2$^{nd}$ generation, sPM = d1*GenScale[18]$^2$/(scale[4]·scale[2]) so

sPM = $h[Tan[\frac{\pi}{18}]^2 Tan[\frac{\pi}{9}](Cot[\frac{\pi}{9}] - Cot[\frac{2\pi}{9}] - (Cot[\frac{\pi}{9}]Tan[\frac{\pi}{18}] + Cot[\frac{2\pi}{9}]Tan[\frac{\pi}{18}])Tan[\frac{\pi}{9}])Tan[\frac{2\pi}{9}] \approx$

0.0123717824965·h. This is exact as long as h is known, but what we want is a scale – and the natural choice is to scale PM relative to the side of N = 9 (using the same h as N = 18), so PMscale = sPM/sN9 = sPM/2hTan[π/9] ≈ 0.0169956  (the h's cancel)

So PMscale =
$\frac{1}{2}Cot[\frac{\pi}{9}](Tan[\frac{\pi}{18}]^2 Tan[\frac{\pi}{9}](Cot[\frac{\pi}{9}] - Cot[\frac{2\pi}{9}] - (Cot[\frac{\pi}{9}]Tan[\frac{\pi}{18}] + Cot[\frac{2\pi}{9}]Tan[\frac{\pi}{18}])Tan[\frac{\pi}{9}])Tan[\frac{2\pi}{9}])$

To check that this new scale is in the scaling field $S_9$:

**AlgebraicNumberPolynomial[ToNumberField[PMscale,GenScale[9],x]** gives

$\frac{11}{3} - \frac{169x}{3} - \frac{25x^2}{3}$ so PMscale = $\frac{11}{3} - \frac{169}{3}Tan[\frac{\pi}{18}]Tan[\frac{\pi}{9}] - \frac{25}{3}Tan[\frac{\pi}{18}]^2 Tan[\frac{\pi}{9}]^2 \approx 0.0169956$

Note that M[1] is the 'matriarch' of this 2$^{nd}$ generation and M[1] is scaled by GenScale[18]= $Tan[\frac{\pi}{18}]^2 \approx .0310912$, so sPM/sM[1] ≈ 0.546637 and this is close to the .455927 ratio for the PM's of N = 7 – which we will study later.

(iii) Now that the exact parameters of the PM tiles are known, it is possible to study the non-regular octagons which surround these PM tiles. In the vector plot below the small S2[3] tiles are also S[2] relative to the M[2] and this implies that they are S[2]'s relative to the First Family of DS4. We will use this family to define the parameters of Dx - and hence show that all the Dx tiles are 'canonical'. (It is an optical illusion that Dx looks shorter than the neighboring tiles.)

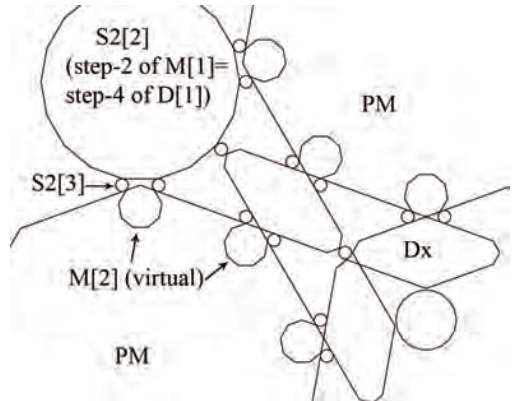

Note that the long edges of Dx are formed from extended edges of PM . This is not unusual. On a larger scale, the edges of the S2[2] tile shown above, are also extended edges of the surrounding PM tiles, so there is a degree of self-similarity between S2[2] and the next-generation S2[3] tiles.

The virtual embedded tiles are formed from rational rotations ($k\pi/18$) of the First Family tiles of DS4. Every vertex of Dx is either a vertex of a canonical tile or a star point of a canonical tile so this is indeed a 'canonical' tile.

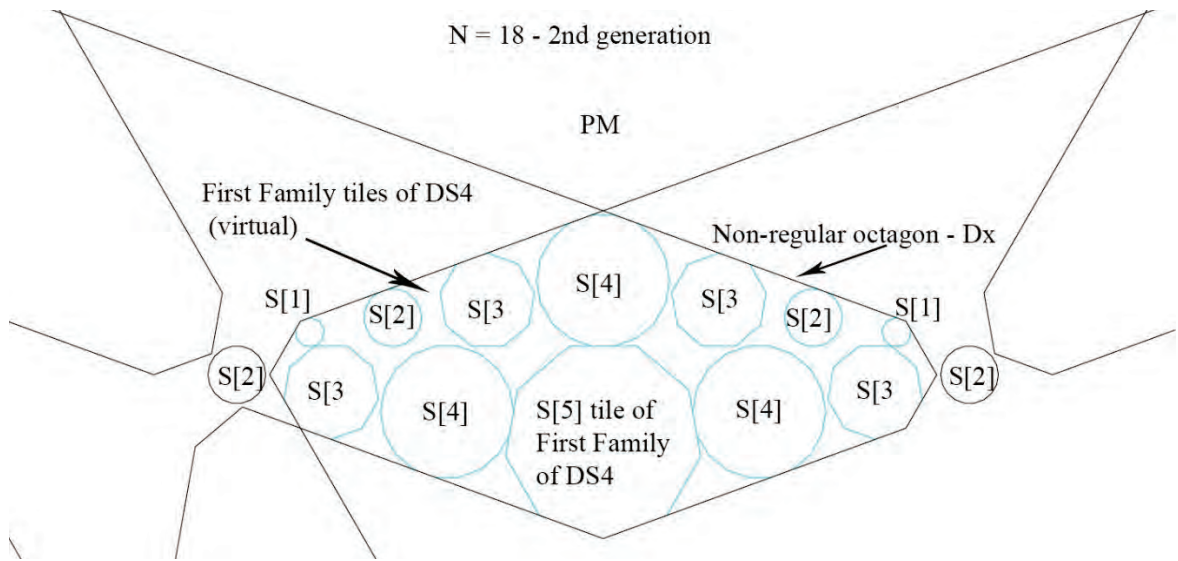

The two external S[2] tiles are really S2[3]'s as discussed above, but relative to M[2] and the First Family of DS4, they are S[2] tiles. They appear to be 'suspended' but the web is always connected, and these tiles do indeed contribute star[2] points to define the tip of the octagon (but the local web is very complex because of the unusual interaction between the scaling of S[1], S[2] and S[3] of DS4. That is why this region has traditionally been known as the Small Hadron Collider.

The long edges of Dx clearly depend on the distance between S[5] and S[3], so this mimics the geometry of the PM tiles. The calculations are very similar to those for PM, and the surprising result is that back in the first generation, the edge length of Dx is simply $2h\tan(\pi/9)$ – which is the side of N = 9 circumscribed about N = 18 (so the h's match).

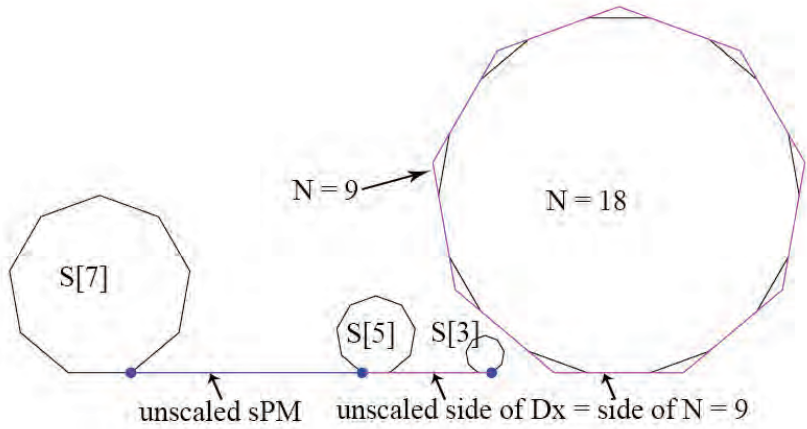

The short edges of Dx span the gap between star[3] and star[7] of S2[3] as shown below

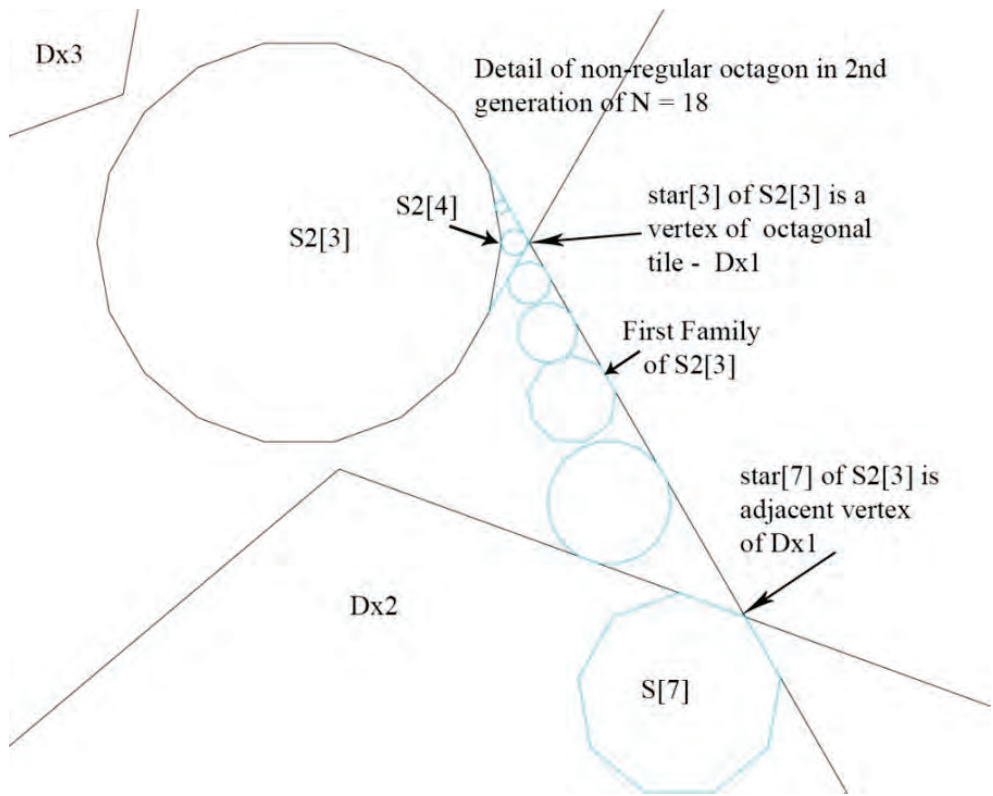

Note that the gap between S3[3] and Dx can be defined by a 'geometric' sequence of S2[k] tiles, so it is $2 \cdot rS2[4][\text{GenScale} + \text{GenScale}^2 + ...]$ where $rS2[4] = rD\text{GenScale}^4/\text{scale}[2]$. Of course this gap can also be found using star[3] of S3[3] – which is a vertex of Dx.

**Example: N= 7**: The PM tiles of N = 9, get their name from the PortalM tiles that play a major role in the second generation for N = 7. As the First Family structure breaks down in the 2$^{nd}$ generation, these PM tiles can be regarded as 'surrogate' S[2] tiles of M[1] for N = 7 and 'surrogate' S[3]'s for N = 9.

As with N = 9, the scale of these PM tiles has no obvious relationship with the First Family scaling of N = 7, but in both cases, there are geometric relationship between tiles that can be used to derive their scale. For both N = 9 and N = 7, they are about half the size of the M[1] tiles.

Since these PM tiles first occur in the second generation, there is no PM[1], so we have traditionally called them PM[2] tiles – but for N = 9 we called them simply PM tiles – so we apologize for the nomenclature. These PM tiles play an important role in the second generation (and all even generations) for N = 7 so these generations are known as 'Portal generations' because the PM tiles have unique local dynamics which resemble a Baroque doorway.

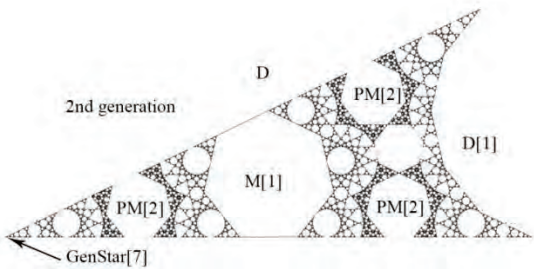
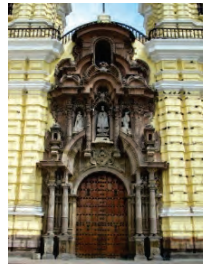

To find the dimensions of PM[2] we will work entirely within N = 7, so star points (and all horizontal distances) will be scaled by h = height of N = 7. (Typically we will use a radius of 1, so h will be $Cos[\pi/7]$ – which is in $S_7$). But our default convention for scaling of tiles is with respect to hD, so N= 14 plays an implicit role. However h and hD share the same side length, so hD/h =·ScaleSwap[7,14] = GenScale[7]/GenScale[14] = 2+ GenScale[7].

Below is a vector plot of this region showing that the PM tiles likely originate on the edges of DS3. There is an alternating sequence of PM[k] tiles and DS3[k] tiles converging to the star[3] point of D[1]. From the left side there is an orthogonal sequence of D[k] tiles converging to the same point- but these tiles alternate real and virtual.

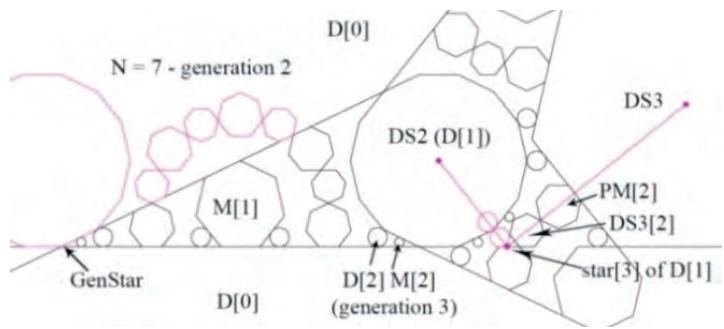

In the 'normal' evolution of families for regular polygons, the step-2 tile of D plays the role the 2nd generation D tile – which we call D[1]. Note that D[1] is also step-2 of DS3 on the right. This is unique to N = 7 and these two tiles share very complex dynamics. The star[3] point of D[1] lies on the top edge of a PM[2] tile and this is the edge that we will calculate.

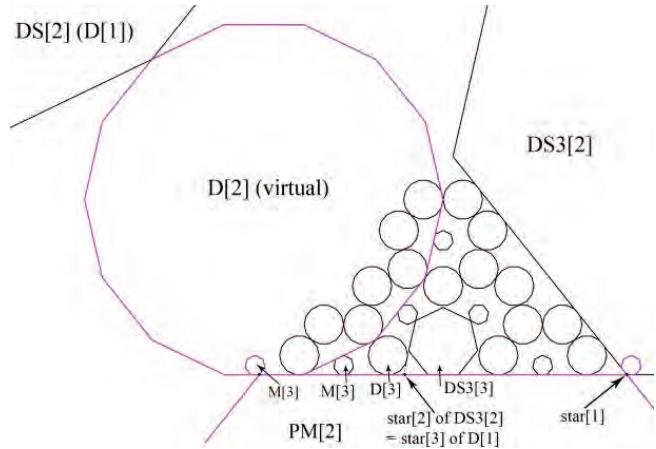

The side of PM[2] is the distance between the star[1] points of the two magenta M[3]'s. All the tiles in between are canonical, so it is an easy matter to find this distance in the First Generation and then rescale by $GenScale^2$ to return to the second generation.

Below is what we call the Short Family for N = 7 – because DS3 can be regarded as the 'matriarch' with DS2 as the step-2 tile. But this short family structure can only be continued on the left if DS2 is regarded as D[1] and **not** a step-2 of DS3 – because these two roles are incompatible.

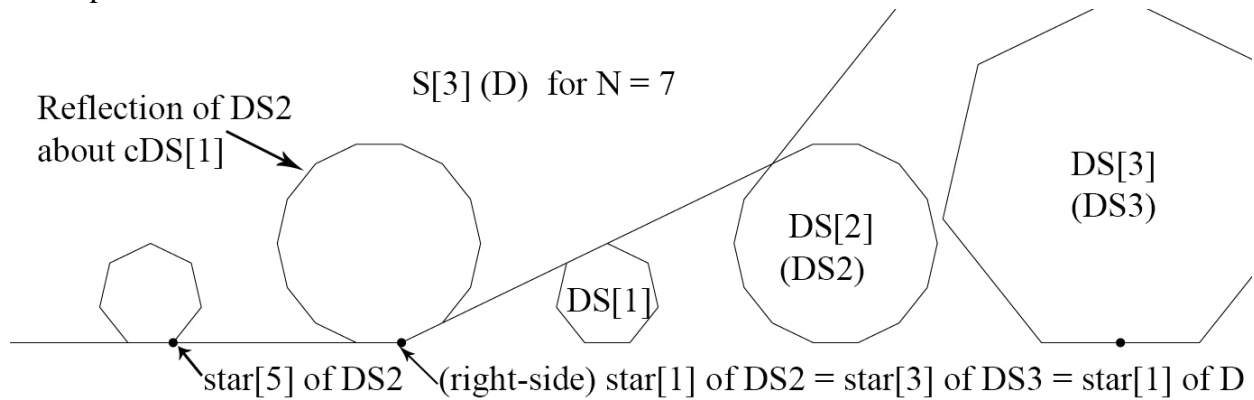

The desired distance is twice the distance shown here between far left and right black markers. The segments above are both differences of star points of known tiles so they can be found in exact form. From left to right:

D1 = star[5][[1]] + star[1][[1]] of DS2 = $(s_5 + s_1)$*hDS2 = (Tan[5*Pi/14] + Tan[Pi/14])*hDS2

D2 = star[3][[1]] of DS3 = $s_3$*hDS3 = Tan[3*Pi/7]*hDS3

So the side of PM[2] scaled up to the First Family is sPMFF = 2h[D1+D2].  (Back in the second generation this is sPM[2] = sPMFF*GenScale$^2$

By definition, these heights are:

hDS[2 ] = hD*GenScale[14]/scale[2] and  hDS[3 ] = hD*GenScale[14]/scale[3] so inside N= 7:

hDS[2] =h( $Sin[\frac{\pi}{7}]Tan[\frac{\pi}{7}]$) and  hDS[3] = h( $Sin[\frac{\pi}{7}]Tan[\frac{3\pi}{14}]$ )

so sPM= h(2·GenScale$^2$·Sin[π/7]·((Tan[5π/14]+Tan[π/14])·Tan[π/7]+Tan[3π/14]·Tan[3π/7])) ≈ 0.0482670611128757742970·h

Therefore PMscale = sPM/sN7 = sPM/(2hTan[π/7] =

Cot[π/7]·GenScale$^2$·Sin[π/7]·((Tan[5π/14]+Tan[π/14])·Tan[π/7]+Tan[3π/14]·Tan[3π/7]))

Which simplifies to  $Sin[\frac{\pi}{7}]Tan[\frac{\pi}{14}]^2 Tan[\frac{\pi}{7}](1+Tan[\frac{\pi}{14}]Tan[\frac{\pi}{7}]+Cot[\frac{\pi}{14}]Tan[\frac{3\pi}{14}])$

This is clearly in the scaling field S$_7$:

**AlgebraicNumberPolynomial[ToNumberField[PMscale,GenScale],x]**  gives

$\frac{1}{8}-\frac{3x}{4}+\frac{5x^2}{8}$  so PMscale = $\frac{1}{8}-\frac{3}{4}Tan[\frac{\pi}{14}]Tan[\frac{\pi}{7}]+\frac{5}{8}Tan[\frac{\pi}{14}]^2 Tan[\frac{\pi}{7}]^2$ ≈ 0.0501137925752

And  therefore sPM = $2hTan[\frac{\pi}{7}]\left[\frac{1}{8}-\frac{3}{4}Tan[\frac{\pi}{14}]Tan[\frac{\pi}{7}]+\frac{5}{8}Tan[\frac{\pi}{14}]^2 Tan[\frac{\pi}{7}]^2\right]$ ≈ 0.0482670611·h

It makes sense to compare PMscale with the other scales. Our interest is a comparison of PM with M[1] – which is the 'matriarch' of this 2$^{nd}$ generation. Since M[1] is scaled by GenScale relative to N = 7 we want the ratio of PMscale and GenScale

PMscale/GenScale =  $Sin[\frac{\pi}{7}]Tan[\frac{\pi}{14}](1+Tan[\frac{\pi}{14}]Tan[\frac{\pi}{7}]+Cot[\frac{\pi}{14}]Tan[\frac{3\pi}{14}])$ ≈ 0.455926999989

Of course sPM/sM[1] is this same ratio so PM[2] is less than half the size of M[1].

These PM tiles define a new scale for N = 7 and it is likely that there are an infinite number of such scales defined by quadratic polynomials in GenScale. Since PM[2] is canonical, PM[2]·GenScale$^k$ is also canonical for all k and these PM tiles seem to exist at least for even k. There are also 'S2-scaled' PM tiles with side (or height) scaled by scale[2] relative to the normal PM's. We call these PMS2 tiles. (The 'main-line' scaling is scale[3] (GenScale) – which is also known as S1 scaling.) These PMS2 tiles occur naturally because any regular tile such as S2 can support a local First Family and these families are identical to the 'main-line' families except for the extra scale[2] factor. Of course this means there are $PM(S2)^k$ tiles for any k.

**Example**: Below is a portion of a 4$^{th}$ generation S2 family with 'matriarch' MS2[3]. The associated PMS2[4] has side sPM[2]·GenScale$^2$/scale[2] ≈ 0.00151843·h.

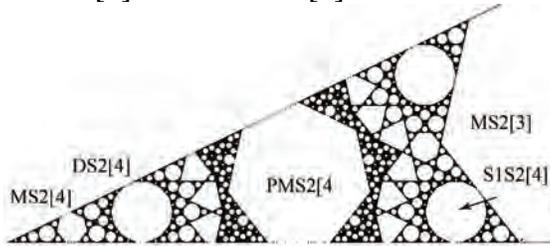

Below is a graphic showing a 'normal' PM[4] with radius rPM[2]·GenScale$^2$. This is the second tiles in what we believe is an infinite chain of PM[2k] tiles converging in a self-similar fashion to the star[3] point of D[1]. For odd generations, these PM tiles alternate with DS3[k] tiles. Converging on the left are D[k] tiles alternating real and virtual (so the small D[4] shown below is virtual). As indicated earlier, all of these tiles exist in the top edge of a PM[2]. This type of alternation is not unusual for N = 7.

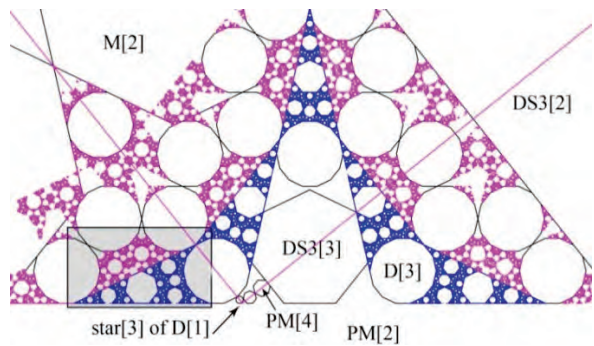

The shaded region is enlarged below in a vector plot – showing how the PM[4]'s are integrated into the 4$^{th}$ generation in place of S[2]'s of M[3]. This region shows clearly the dichotomy between the 'normal' (blue) dynamics of the 4$^{th}$ generation and the unpredictable (magenta) dynamics across the border.

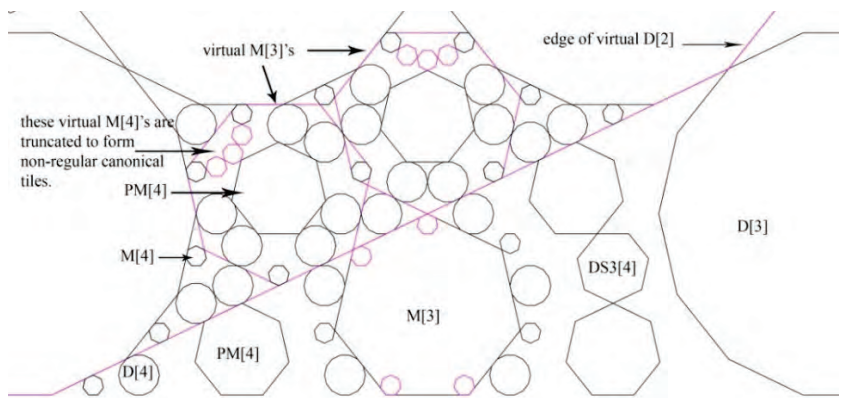

On the left above, note the proximity of D[4] and PM[4]. This is normal for an 'even' (Portal) generation. The problem is that D[4] is a 'main-line' tile with S1 (GenScale) scaling and this scaling is incompatible with PM scaling. Therefore the dynamics on the right side of D[4] are unpredictable while the region on the left of D[4] is a typical 5$^{th}$ generation with 'matriarch' M[4]. The self-similar 'towers' on the edges of PM[4] are a testament to the clash in scaling between D[4] and PM[4]. The situation is much worse when S2 scaling is introduced into the mix, as shown below. Note that the proximity of D[4] and PM[4] is unchanged from above, but now they are surrounded by S2[4] tiles – which are larger and more dominant because hS2[4]/hD[4] = 1/scale[2 ] ≈ 2.60388.

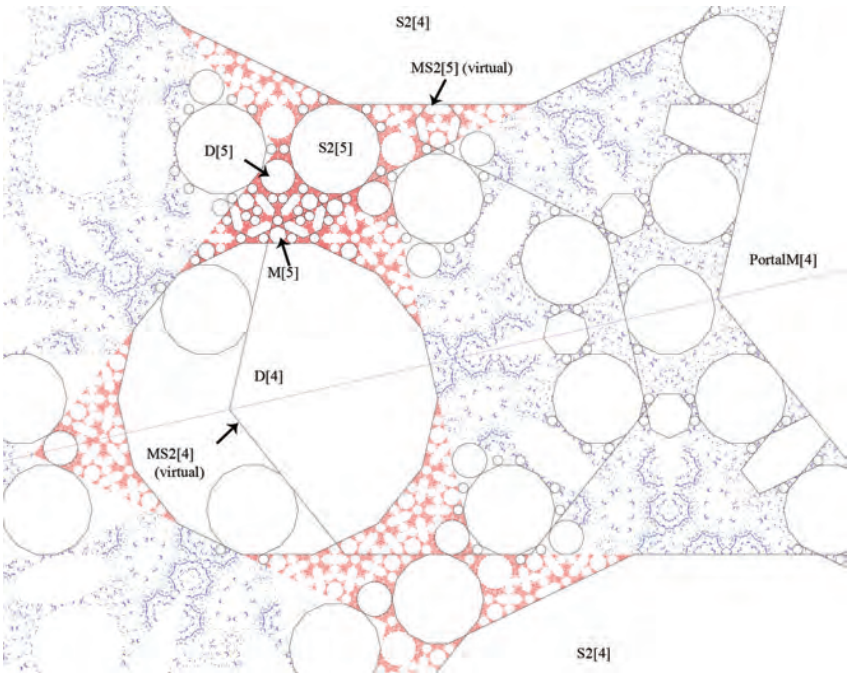

This region occurs at the top of S1 as shown below - and it is replicated at star[2] and star[1] of N =7. Since S1 is adjacent to S2, their scales interact in an unpredictable fashion and both families have PM tiles which cause even more turmoil. The light blue dynamics above are apparently the remnants of the PM 'towers' from a normal 4$^{th}$ generation.

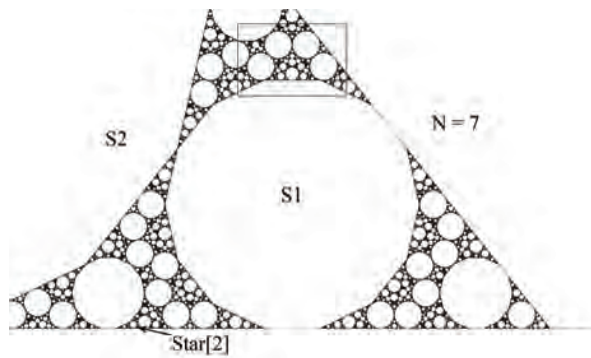

The invariant orange region is enlarged below. It shows the (one sided) interaction between S2 and S1 dynamics. D[4] is the lone 4[th] generation S1-scaled tile in this S2-dominated region – so D[4] is a 'stranger in a strange land'.

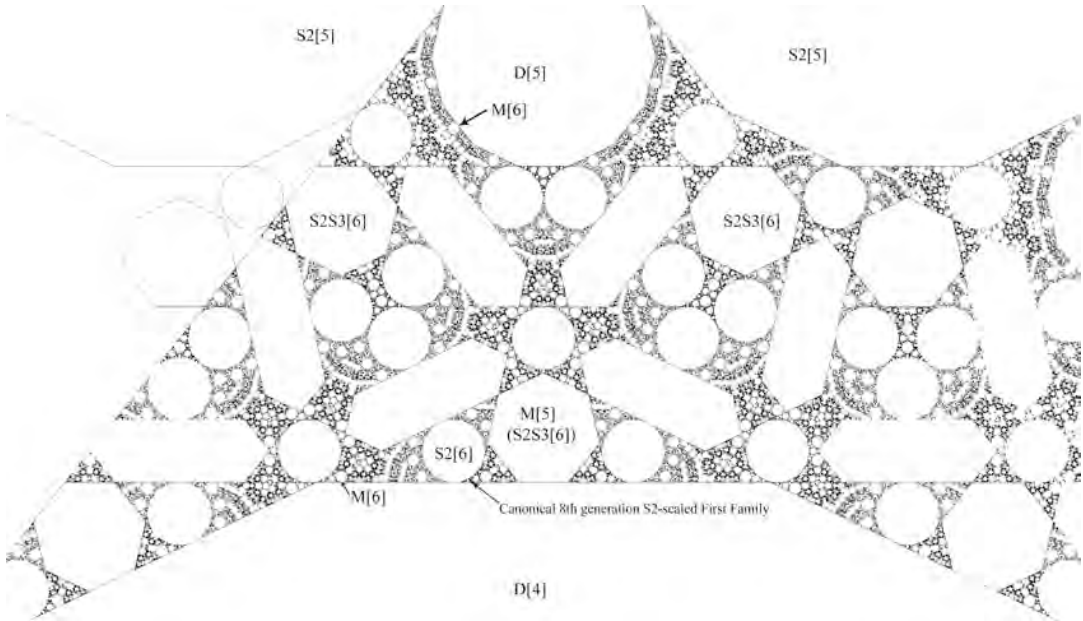

The small M[5] on the top edge of D[4] has no matching D[5] – so there is no hope of a main-line 5[th] generation here. The virtual M[4] outlined in the earlier plot, has a 'normal' relationship with the large S2[4] tiles, and as indicated earlier, it is not unusual for M tiles to have next generation M's at their vertices. These virtual M[5]'s play a large role in the large-scale dynamics at the top of D[4] and one of their 'progeny' is M[5] at the top of D[4].

The problem is that M[5] of D[4] is also a step-3 of the large S2[5] tiles and these roles are totally incompatible. Hence this region has geometry which has never been observed before, but it is our contention that all the tiles shown here are 'canonical' - even though we do not fully understand their evolution.

For example each of the elongated hexagons are clearly formed from two (virtual) S2[6] tiles which are linked by their extended edges, and hence they are trivially canonical. This type of linkage is similar to that found for $N = 9$ earlier and it is also common with $N = 11$. Below is the 'chaotic' region to the left of M[5] showing a normal S2-scaled 8[th] generation in blue.

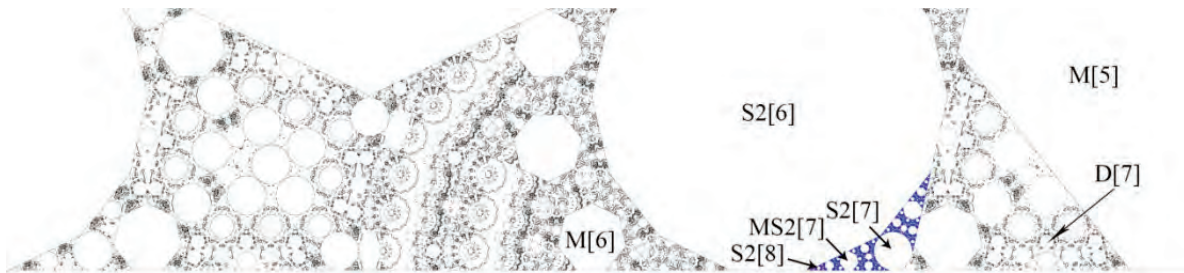

**Example**: N = 11

In [H2] we present an introduction to the dynamics of the 'quintic' polygon N = 11. Since this is a 4k+3 prime, there are no obvious signs of 'generation' scaling after the first two generations. In particular the 'D' chain contains only the minimal D and D[1] and there are no signs of further D[k] tiles. The 'M' chain survives for one more generation as M, M[1] and M[2] (they exist on alternating edges of D[1]). In general these D and M tiles co-exist with each other – so there seems little hope of extended generations. This means that any analysis will be of a very different nature than the examples above. It has taken a long time to develop analytic methods that can be applied to N = 11. It is our hope that these methods can be applied on a broader scope.

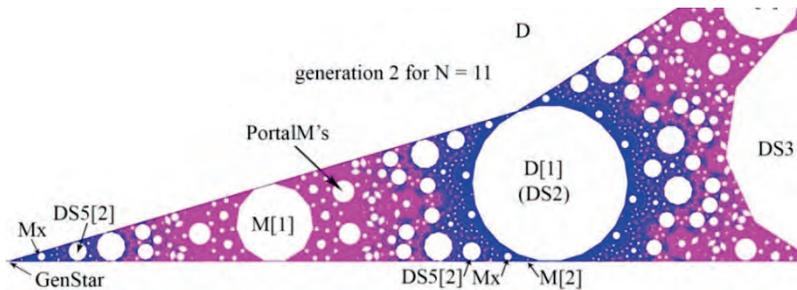

Below is a vector plot of this second generation region – showing what an ideal second generation would look like. It is easy to make plots like these because it simply involves scaling the FirstFamily of N = 11 by GenScale (scale[5]) and translating to the center of M[1] and then reflecting to get the cyan right side.

FirstFamily of M[1]-D[1] = **TranslationTransform[cM[1]]/@(FirstFamily*GenScale)**;

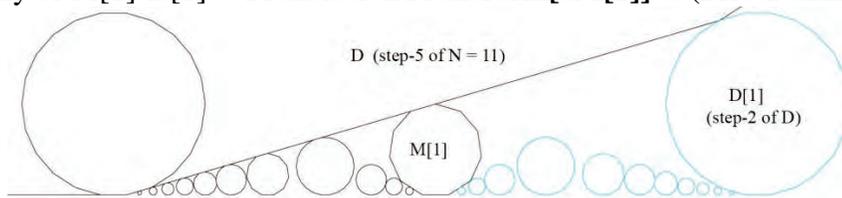

Comparing this with the plot above, it is clear that the reflective symmetry is preserved, but the only First Family survivors are D[1], M[1], M[2] on the edge of D[1] (too small to see here) and DS5[2] - and their symmetric counterparts at the GenStar point. Below is an enlargement of the region at the foot of D[1] – showing M[2], Mx and DS5[5] – which we will call simply DS5 (step 5 of D[1]). Since the 'mystery' tile Mx is 'conforming' with the bounds of the known cyan family, it would seem that there should be some connection between Mx and this family – but no connection was ever found. Here we show the rather bizarre genesis of Mx.

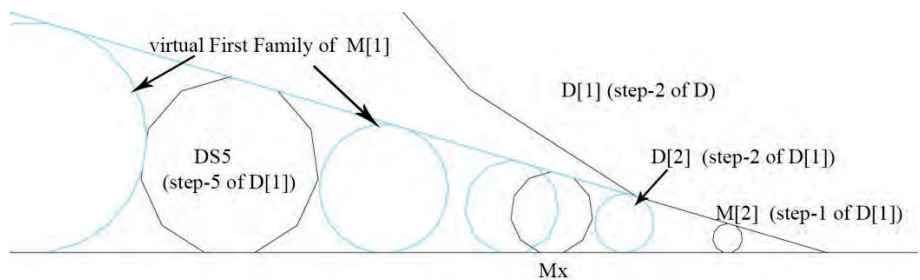

**Two-Star Lemma**: Any two star points of a regular tile determine the canonical scale of the tile.

**Proof**: The canonical scale of a regular n-gon tile P is determined by its height hP. Note that any star point of P is sufficient to determine the height (apothem), if the center is known, since by definition Tan[kPi/n] = star[k][[1]]/hP. Clearly two star points are sufficient to determine the relative displacement of either as follows:

Without loss of generality we can assume that the two given star points p1 and p2 lie on a common horizontal edge of P and that the 'indices' of each are known (by simply extending the edges of P). Suppose that p1 is star[k] and p2 is star[j]. The distance between them is the horizontal displacement - which we call c. In the coordinate system based on the center of P, this displacement is unchanged. If x and y are the coordinates of p1 and p2 relative to P, then
$$\frac{x}{Tan[\frac{kPi}{n}]} = \frac{x+c}{Tan[\frac{jPI}{n}]}$$ because both are equal to hP. Therefore

x = (c/Tan[jPi/n])(Cot[kPi/N]-Cot[[jPi/n]) and hP = x/Tan[kPi/n] = y/Tan[jPi/n]. □

For Mx , star5 is known because it is a vertex of D[1]. We will find the exact coordinates of the star[4] point of Mx – shown in blue below. To do this we will use the 'web' algorithm as described in [H2]. The extended horizontal edges of D[1] also lie on the extended edge of N = 11, so the segment H1 shown below is part of the level-0 web for N = 11. There are matching segments for all 11 D tiles and all of these segments together are iterated by the web algorithm. There is no guarantee that any of these 11 segments will map back to D[1], but since D[1] is a step-2 tile of D, and the D's map to each other, these 11 segments will map primarily to the D[1]'s. It takes 13 iterations of the web for H2 to arise as a symmetric version of H1. Technically H2 is not an image of H1 because all images of H1 must be parallel to H1 – but H2 is the image of an H1 in another D). We will use H2 to determine the strange 'offsets' that eventually determine the star[4] point of Mx – although other choices work also.

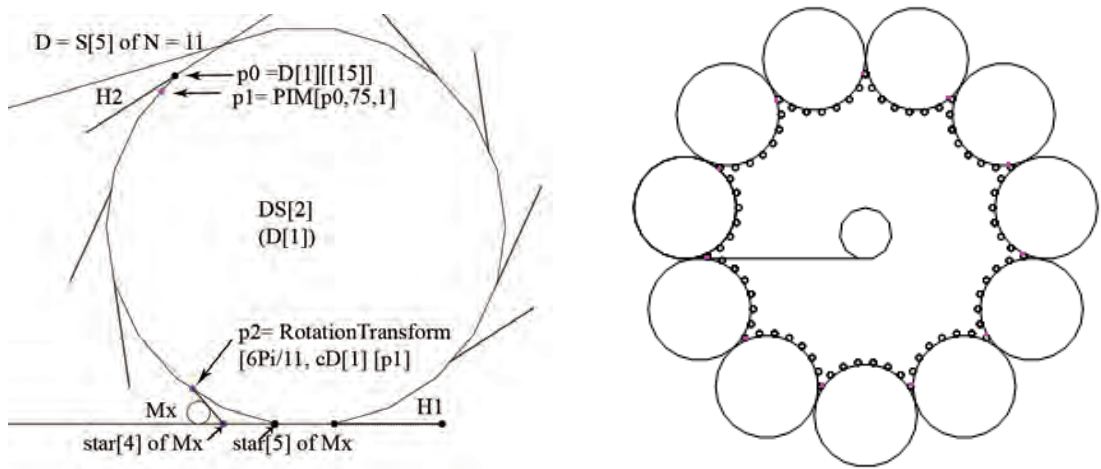

After 9*11 iterations of the web, H1 (and its matching segments) will yield the nine segments shown on the left above – including the crucial truncated segment that determines a side of Mx. All that we need is an 'exact' formula for p2, because p2 determines star[4] - since all edges must be parallel to edges of D[1]. And then star[4] and star[5] together determine Mx.

The steps to find p2 are outlined on the graphic above. The vertex p0 is known to arbitrary precision and we suspect that it maps to p1, but vertices technically have no image under $\tau$. However vertices can have one-sided limit points and we will show that the limit point of p0 is indeed the edge point p1. (So vertices can indeed map to edges.) To do this, pick a point on the segment H2 'close' to p0. It is not critical that the error be very small – because the web tells us that the interior portion of H2 does not split for at least 200 iterations – and hence all neighbors of p0 on this segment will temporarily share the same orbit.

Since $\tau(p) = 2c_k - p$ for some vertex $c_k$ of N = 11, $\tau^k(p) = (-1)^k p + 2Q$ where Q is an alternating sum of vertices. These vertices are known to arbitrary accuracy, so if the initial point p is known, the orbit of p is 'exact' – as long the sequence of 'corners' is known. This is what we call the 'corner sequence'. Mathematica tracks these sequences as part of our analysis and there are 'checks and balances' that can be used to determine the validity of such a sequence.

The center of D[1] has period 77 and all other points in D[1] have period 154. (Only centers can have odd period). The 'step sequence' of D[1] is the number of vertices advanced on each iteration (the first difference of corner indices). The step sequence for **all** points of D[1] is {5,5,5,5,5,5,4} (compared with {5} for the 'patent' D). So the step sequence of D[1] is period 7 and of course this divides the actual period.

It makes sense that the points adjacent to D[1] will have similar step sequences. The first 20 corner indices for a point close to p0 are {11,5,10,4,9,3,8,1,6,11,5,10,4,9,2,7,1,6,11,5} and it is clear they advance by {5,5,5,5,5,5,4} (mod 11), but of course this sequence will eventually break down. Clearly no web point, edge point or vertex point can have a periodic orbit.

So we will use an approximation of p0 to determine the corner sequence Q and then simply plug p0 into this sequence. This allows us to define an 'ideal' (one-sided) orbit for p0 - and the results can be checked in many ways – including the fact that the resulting hMx must be in the scaling field $S_{11}$.

Here are the calculations:

(i) rDS[2] = $Sec[\frac{\pi}{22}]Sin[\frac{\pi}{11}]Tan[\frac{\pi}{11}]$

cDS[2] = $\{-Cos[\frac{\pi}{11}]Cot[\frac{\pi}{22}](1-Tan[\frac{\pi}{11}]^2), -Cos[\frac{\pi}{11}](1-Tan[\frac{\pi}{11}]^2)\}$

(In the calculations to follow we will give just the first coordinate of the resulting points)

p0[[1]] = DS[2][[15]]][[1]] = $\{-Cos[\frac{\pi}{11}]Cot[\frac{\pi}{22}] + Sin[\frac{\pi}{11}](Cot[\frac{\pi}{22}] - Sec[\frac{\pi}{22}]Sin[\frac{5\pi}{22}])Tan[\frac{\pi}{11}]\}$

(ii) To get the step sequence for p0, we will use a point px that is on H2 and distance $10^{-8}$ from p0. Because we have been using step sequences to generate projections (see Projections in the PDF list at dynamicsofpolygons.org), our software is orientated toward the 'return map' $\tau^2$ which generates sequences of vectors of the form $c_k - c_j$. We will generate at least 150 points in the orbit of px to get the 75 displacements that define p1. (The 'suspected' number of iterations is determined by visually tracking the 'return map'). The indices of the corner sequences are generated by IND[px,150]. The first 20 indices are given above. Then P1= PIM[p0,75,1] takes the 150 indices and constructs the first 75 terms of the vector sequence to reproduce the actual 'orbit' using p0 as the 'surrogate' for px.

This is called the P1 'projection' of p0 - and it is a series of vector displacements based on the indices found above. (This is an example of an 'algebraic graph' as defined in [S2].) . The first element in P1 is p0. Our interest is in the 75$^{th}$ element of P1 so set p1= P1[[75]]. The simplified form of p1[[1]]=

$$-6Cos[\frac{3\pi}{22}] - Cos[\frac{\pi}{11}]Cot[\frac{\pi}{22}] + 4Sin[\frac{\pi}{11}] + 8Sin[\frac{2\pi}{11}] + Cot[\frac{\pi}{22}]Sin[\frac{\pi}{11}]Tan[\frac{\pi}{11}] - Sec[\frac{\pi}{22}]Sin[\frac{\pi}{11}]Sin[\frac{5\pi}{22}]Tan[\frac{\pi}{11}]$$

(iii) To find p2, rotate p1 in an 'exact' fashion – which is the default for Mathematica:

p2= **RotationTransform[6*Pi/11, cDS[2]][p1]**

$$p2[[1]] = \frac{(-1)^{15/22}(-1 - 4(-1)^{1/11} + 6(-1)^{3/11} - 7(-1)^{4/11} + 7(-1)^{6/11} - 6(-1)^{7/11} + 4(-1)^{9/11} + (-1)^{10/11})}{2(-1 + (-1)^{1/11})(1 + (-1)^{2/11})}$$

(the complex part will vanish)

(iv) The slope of this line segment determined by p2 is known exactly, so it is trivial to solve for star[4], but the exact result is somewhat messy: Star[4][[1]] =

$$\frac{1}{-Cos[\frac{5\pi}{22}] + Sin[\frac{2\pi}{11}]}(Cos[\frac{2\pi}{11}] - Sin[\frac{5\pi}{22}])(\frac{1 + 12(-1)^{1/11} + (-1)^{2/11} + (-1)^{3/11} - 3(-1)^{4/11} - 9(-1)^{5/11} - 2(-1)^{6/11} + 2(-1)^{7/11} + 9(-1)^{8/11} + 3(-1)^{9/11} - (-1)^{10/11}}{2 + 2(-1)^{2/11}}$$

$$-Cos[\frac{\pi}{11}] - Cos[\frac{\pi}{22}]Cos[\frac{\pi}{11}]Cot[\frac{\pi}{22}] - Cos[\frac{\pi}{11}]Sin[\frac{\pi}{22}] +$$

$$\frac{(-1)^{15/22}(-1 - 4(-1)^{1/11} + 6(-1)^{3/11} - 7(-1)^{4/11} + 7(-1)^{6/11} - 6(-1)^{7/11} + 4(-1)^{9/11} + (-1)^{10/11})(-Cos[\frac{5\pi}{22}] + Sin[\frac{2\pi}{11}])}{2(-1 + (-1)^{1/11})(1 + (-1)^{2/11})(Cos[\frac{2\pi}{11}] - Sin[\frac{5\pi}{22}])} +$$

$$Sin[\frac{\pi}{11}]Tan[\frac{\pi}{11}] + Cos[\frac{\pi}{22}]Cot[\frac{\pi}{22}]Sin[\frac{\pi}{11}]Tan[\frac{\pi}{11}] + Sin[\frac{\pi}{22}]Sin[\frac{\pi}{11}]Tan[\frac{\pi}{11}] +$$

$$Cos[\frac{\pi}{11}]Sin[\frac{\pi}{22}](1 - Tan[\frac{\pi}{11}]^2) - Cos[\frac{\pi}{22}](Sec[\frac{\pi}{22}]Sin[\frac{\pi}{11}]Tan[\frac{\pi}{11}] - Cos[\frac{\pi}{11}]Cot[\frac{\pi}{22}](1 - Tan[\frac{\pi}{11}]^2)))$$

≈ -6.1402943334830363571499067 (this assumes a radius 1 for N = 11). The second coordinate of star[4] is the same as the 'base' edge of N = 11 – namely minus the apothem.

Star points like these can be multiplied (or divided) by Tan[Pi/11] to yield elements in the scaling field $S_{11}$ - which is generated by GenScale (scale[5] of N = 11). This provides a check on the validity of the calculations so far and gives a more manageable formula for star[4]:

**AlgebraicNumberPolynomial[ToNumberField[star[4][[1]]*Tan[Pi/11],GenScale],y]** =

$-\frac{111}{32} + \frac{621y}{16} + \frac{133y^2}{8} - \frac{25y^3}{16} - \frac{13y^4}{32}$  Setting y = GenScale will yield star[4][[1]]*Tan[Pi/11].

(v) Finally the Two-Star Lemma will yield the coordinates of star[4] relative to Mx – and hence the scale of Mx: Set c = Star[4][[1]] – star[5][[1]]; and y1= Tan[4*Pi/11], y2 = Tan[5*Pi/11], then

MxStar[4][[1]] = $\frac{c}{y2}(\frac{1}{y1} - \frac{1}{y2})$. As expected this result is 'messy' but exact. But what we really want is hMx = -MxStar[[4]][[1]]/Tan[4*Pi/11] which is given below:

$$-\frac{40 - 4(-1)^{1/11} + 24(-1)^{2/11} - 44(-1)^{3/11} + 8(-1)^{4/11} - 36(-1)^{5/11} + 29(-1)^{6/11} - 5(-1)^{7/11} + 43(-1)^{8/11} - 19(-1)^{9/11} + 12(-1)^{10/11}}{2(1+(-1)^{2/11})(1-(-1)^{2/11}+(-1)^{6/11}+(-1)^{7/11})(Cos[\frac{5\pi}{22}] - Sin[\frac{2\pi}{11}])(Tan[\frac{\pi}{22}] - Tan[\frac{3\pi}{22}])}$$

≈ 0.0047466778294017318376548  (again the imaginary part vanishes)

This height is actually a scaled distance so it is in $S_{11}$ and Mathematica will yield a polynomial form for hMx:

**AlgebraicNumberPolynomial[ToNumberField[hMx,GenScale,z]=**

$\frac{21}{32} - \frac{241z}{16} - \frac{35z^2}{4} - \frac{3z^3}{16} + \frac{11z^4}{32}$  so hMx =

$\frac{21}{32} - \frac{241}{16}Tan[\frac{\pi}{22}]Tan[\frac{\pi}{11}] - \frac{35}{4}Tan[\frac{\pi}{22}]^2 Tan[\frac{\pi}{11}]^2 - \frac{3}{16}Tan[\frac{\pi}{22}]^3 Tan[\frac{\pi}{11}]^3 + \frac{11}{32}Tan[\frac{\pi}{22}]^4 Tan[\frac{\pi}{11}]^4$

To convert this to a scale, divide by hN11. Since we are assuming a radius of 1, hN11 = Cos[π/11], so MxScale = hMx/hN11 = hMx/Cos[π/11] ≈ 0.004947068879. This scale is in $S_{11}$ because 1/Cos[π/11] = 1 + GenScale.

This proves that Mx is a 'canonical' tile and the matching MxScale is in the scaling field $S_{11}$ so it can interact algebraically with the 5 existing scales, but we do not know if any of the resulting variants of Mx exist. It is trivial to construct an ideal First Family for Mx, but this family has no obvious 'hits' with existing tiles – it would have been a great surprise if there were any such relationships because there does not appear to be any clear-cut connections between Mx and the neighboring canonical tiles.

The offsets that we used here to define Mx, result from 'mutations' in the predicted orbits of the vertices of the D[1] tiles – where points temporarily fail to map to the edges of D[1] tiles. Therefore we can no longer predict their corner sequence and are forced to rely on Mathematica for the actual corner sequence.

Based on the nature of the GenScale polynomial above, it seems unlikely that there are simple algebraic relationships between the MxScale and the remaining scales of N = 11. We know that all the scales can be generated from GenScale but it would be challenging to find any meaningful relations between these polynomials and the Mx polynomial above.

Vertex 'offsets' such as these are clearly a part of the evolution of 'most' regular polygons – and it appears that the typical edge of a tile may code a significant amount of the overall dynamics. (Every extended edge of a regular polygon 'codes' the entire dynamics because the evolution of these points defines the singularity set.) In the case of p2 studied here, this appears to be just the first of an infinite sequence of 'offsets' - some of which can be observed in the web plot below.

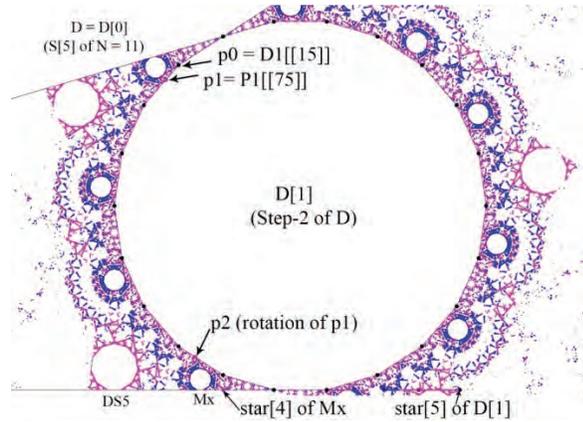

The blue points show part of the sphere of influence of Mx. They are generated by initial edge points between (right-side) star[3] and star[4] of D[1]. The magenta points are generated by neighboring intervals. In the limit there seems to be a level of interaction between intervals and it seems likely that there are points on the trailing edge of D[1] with 'locally dense' orbits.

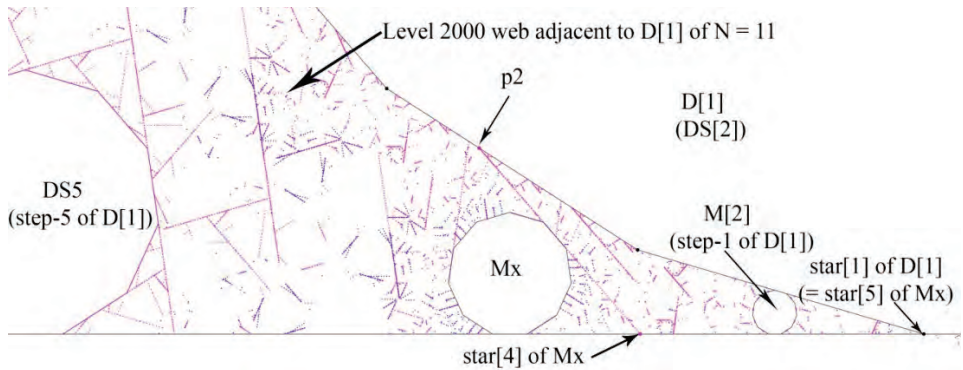

Below is an overlay of this web on an extended web that is generated to a level that varies between one and two billion iterations.

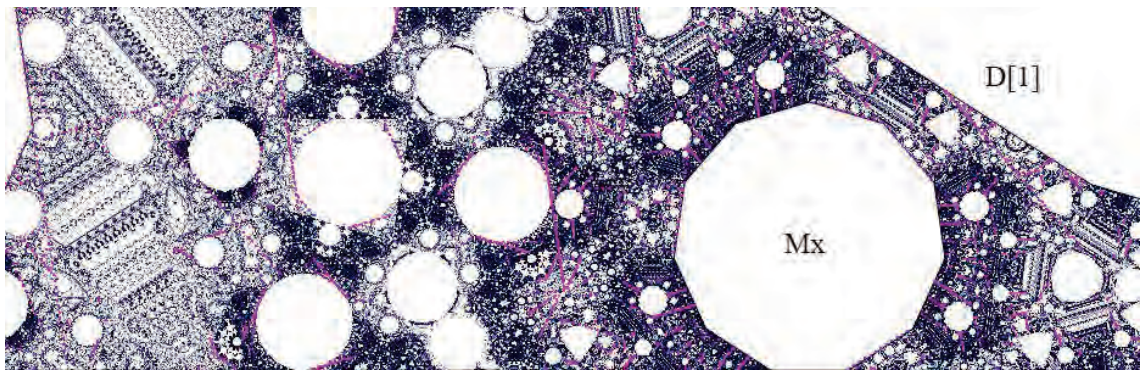

The scaling of Mx is just a small 'dent' in the armor of N = 11, but these methods show promise for further analysis. One practical benefit of knowing Mx exactly, is that it is possible to study the edge-dynamics at close range. These dynamics appear to be locally invariant and hence they resemble the strange invariant 'islands' observed in the region surrounding Mx. It may be that these invariant islands reappear at smaller scales in the dense region close to Mx – like a local 'Mandelbrot set'.

**Summary**:

(i) For a given regular N-gon, we have defined the First Family of 'tiles' - with nucleus consisting of the S[k] tiles. Each S[k] tile corresponds to 'star' point $s_k = \tan(k\pi/N)$ and matching scale $s_1/s_k$. The primitive star points and scales satisfy $(k,N) = 1$, and based on the results of Siegel and Chowla, we have shown that the primitive star points and scales are independent with degree $\varphi(N)/2$.

(ii) These primitive scales lie in the maximal real subfield of the cyclotomic field $\mathbb{Q}(\zeta_N)$, which we denote as $\mathbb{Q}(\zeta_N)^+$. They have degree $\varphi(N)/2$ so they are a unit basis for this subfield – which is traditionally generated by $\lambda_N = 2\cos(2\pi/N)$. In all cases GenScale[N] or GenScale[N/2] can also serve as unit generators for $\mathbb{Q}(\zeta_N)^+$.

(iii) The primitive scales have the same complexity as the scales used by other investigators – since the dynamics and scaling of affine piecewise rational rotations are typically restricted to $\mathbb{Q}(\lambda)$ where $\lambda$ is the trace of the rotation matrix. Therefore it seems likely that the scaling needed to describe the dynamics of maps such as the outer-billiards map and Digital Filter map will lie in the maximal real subfield $\mathbb{Q}(\zeta_N)^+$ of the cyclotomic field $\mathbb{Q}(\zeta_N)$. The canonical scales described here may help to unify investigations in this area.

(iv) High in the 'wish' list are practical methods that can be used to relate the geometric and temporal scaling – maybe in some form of 'power' law – which in turn can be used to estimate the Hausdorff dimension of the simgulartity set. This was succesfully done in [LKV] using $\lambda_{14}$ and a "*finite order recursive tiling with scaling factors given by algebraic units*". But the recursive tiling required an extensive 'catalog' of scaling domains – which is typically not available.

(v) The recent Appendix C describes practical methods for applying the scaling results described in this paper. It shows that there is reason to believe that all tiles that arise from the outer billiards map on a regular polyon, are 'canonical' - and hence can be scaled by elements of the scaling field $S_N$.

**Links**:

(i) The author's web site at DynamicsOfPolygons.org is devoted to the outer billiards map and related maps from the perspective of a non-professional. This file is available there as FistFamilies.pdf – just click on PDFs.

(ii) A Mathematica notebook called FirstFamily.nb will generate the First Family and related star polygons for any regular polygon. It is also a full-fledged outer billiards notebook which works for all regular polygons. This notebook includes the Digital Filter map (which is only applicable for N even). The default height is 1 - to make it compatible with the Digital Filter map. For investigations without the Digital Filter map, it may be preferable to use one of the notebooks described below – with the more natural convention of radius 1. Of course the scaling is independent of these choices, so it is an easy matter to mix the conventions.

 (iii) Outer Billiards notebooks for all convex polygons (radius 1 convention for regular cases). There are four cases: Nodd, NTwiceOdd, NTwiceEven and Nonregular.

(iv) The open source PARI software at pari.math.u-bordeaux.fr has impressive facilities for computer algebra and algebraic number theory. There is an excellent introduction to Galois Theory and PARI in *Fields and Galois Theory* by J.S Milne – which is available at www.jmilne.org/math/

(v) We also recommend the open source (GPL- GNU) Sage software which was originally devised in 2005 by William Stein at U.C. San Diego, as a Python and C++ based library. It currently has hundreds of developers around the world and an extensive library of routines for number theory, algebra and geometry. Sage runs on an Oracle Virtual Machine which will install on almost any operating system.